\documentclass[11pt]{article}

\usepackage[T1,T2A]{fontenc}
\usepackage[utf8]{inputenc}
\usepackage{amsfonts}
\usepackage{amssymb}
\usepackage{epsfig}
\usepackage{tikz}
\usepackage{float}
\usepackage{theorem}

\usepackage{color}

\newtheorem{theorem}{Theorem}[section]
\newtheorem{lemma}[theorem]{Lemma}
\newtheorem{corollary}[theorem]{Corollary}

{\theorembodyfont{\rmfamily} 

}

\title{Joining dessins together}
\date{\today}
\author{Gareth A.~Jones}

\begin{document}

\maketitle

\begin{abstract}
\noindent An operation of joining coset diagrams for a given group, introduced by Higman and developed by Conder in connection with Hurwitz groups, is reinterpreted and generalised as a connected sum operation on dessins d'enfants of a given type. A number of examples are given.
\end{abstract}

\medskip

\noindent{\bf Subject class:}  05C10, 14H57, 20B25, 30F10.

\medskip




\section{Introduction}\label{intro}

Around 1970 in Oxford, in a series of unpublished lectures on Generators and Relations, Graham Higman introduced a technique which he called `sewing coset diagrams together'. The basic idea was to combine two transitive permutation representations of a group $\Gamma$, of finite degrees $n_1$ and $n_2$, to give a transitive permutation representation of $\Gamma$ of degree $n_1+n_2$. In certain cases this could be done by taking permutation diagrams for these two representations, that is, coset diagrams for point-stabilisers with respect to certain generators for $\Gamma$, and joining them to obtain a larger connected diagram, also representing $\Gamma$, by replacing some pairs of fixed points of a generator of order 2 with 2-cycles. Under suitable conditions, this could be iterated to give permutation representations of $\Gamma$ of even larger degrees.

One aim for using this technique was to obtain relatively concise presentations for various finite groups, such as ${\rm PSL}_2(q)$ for certain values of $q$, and the sporadic simple group ${\rm J}_1$ of order 175,560 which had been discovered by Janko~\cite{Janko-65, Janko-66} a few years earlier. Another aim was to obtain finite quotients of certain finitely presented infinite groups, such as various triangle groups (including the modular and extended modular group), and Coxeter's groups $G^{p,q,r}$ (see~\cite[\S7.5]{CM}) with presentations
\[\langle A, B, C \mid A^p=B^q=C^r=(AB)^2=(BC)^2=(CA)^2=(ABC)^2=1\rangle.\]
Higman gave an outline argument, using this technique, that the alternating group ${\rm A}_n$ is a Hurwitz group, that is, a quotient of the triangle group\footnote{For consistency with dessin notation, we prefer $\Delta(3,2,7)$ to the isomorphic $\Delta(2,3,7)$.}
\[\Delta=\Delta(3,2,7)=\langle X, Y, Z\mid X^3=Y^2=Z^7=XYZ=1\rangle,\]
for all sufficiently large values of $n$. This was fully proved, extended and made more precise by Higman's student Marston Conder in his 1980 D.~Phil.~thesis and his first published paper~\cite{Conder-80}: he showed that ${\rm A}_n$ is a Hurwitz group for all $n\ge 168$, and also for a specific set of smaller values of $n$, ranging from 15 to 166. (See~\cite{Conder-10} for a useful recent survey of Hurwitz groups, also by Conder.)  He also applied this technique to realise all sufficiently large alternating and symmetric groups as quotients of the extended triangle group $\Delta[3,2,7]$.

This method of joining coset diagrams together was also applied by Wilson Stothers~\cite{Stothers-74, Stothers-77} in 1974 to the modular group $\Gamma={\rm PSL}_2(\mathbb Z)$, and in 1977 to the triangle group $\Delta$, in order to study the possible specifications (essentially signature, plus cusp-split in the case of $\Gamma$) for their subgroups of finite index. This work arose out of his Cambridge PhD thesis, written in 1972 under the supervision of Peter Swinnerton-Dyer. 

Conder's technique for finding alternating quotients of $\Delta$ was later extended by Pellegrini and Tamburini~\cite{PT}, who showed that the double cover $2.{\rm A}_n$ of ${\rm A}_n$ is a Hurwitz group for all $n\ge 231$. It was transferred by Lucchini, Tamburini and Wilson~\cite{LTW} from permutations to matrices in order to show that the groups ${\rm SL}_n(q)$ are all Hurwitz groups for $n\ge 287$, while Lucchini and Tamburini~\cite{LT} proved a similar result for various other families of classical groups. Conder's coset diagrams have also been used in~\cite{JP} to construct Beauville surfaces from pairs of Hurwitz dessins.

In~\cite{Conder-81}, Conder extended the technique to generalise his results in~\cite{Conder-80} to triangle groups of type $(3,2,k)$ for all $k\ge 7$; his diagrams have been used in~\cite{Jones-18a} to show that such groups have uncountably many maximal subgroups, and in~\cite{Jones-18b} to realise all countable groups as automorphism groups of maps and hypermaps of various types. By extending the joining operation to Dyck groups, Everitt~\cite{Everitt} showed that each non-elementary finitely generated Fuchsian group has all sufficiently large alternating groups as quotients.

The theory developed by Higman and Conder was entirely 1-dimensional: their basic tools were graphs, with edges implicitly labelled and directed. However, when the theory is applied to triangle groups (including extended triangle groups), as in many of Higman's examples, and in almost all of the others discussed above, the coset diagrams can be interpreted as maps on compact oriented surfaces (often the sphere), and therefore, following Grothendieck \cite{Gro}, as dessins d'enfants, that is, as `pictures' of algebraic curves defined over the field $\overline\mathbb Q$ of algebraic numbers, shown in the Appendix.
Moreover, their joining operations correspond to certain type-preserving connected sum operations on dessins, allowing Hurwitz dessins of arbitrary size and genus to be created from simple planar ingredients. Thus Conder's main result in~\cite{Conder-80} can be interpreted as using connected sums to construct, for each integer $n\ge 168$ (and many smaller $n$), a planar dessin with monodromy group ${\rm A}_n$, so that this group is realised as a Hurwitz group on its Galois cover, of genus $1+n!/168$.

The aim of this paper is to explain these reinterpretations, with illustrative examples, and to generalise the joining operations as connected sum operations on surfaces. In particular, whereas Higman and Conder based their technique on using the fixed points of a generator of order 2 of a triangle group, we shall extend this to fixed points of a generator of any order. 



\medskip

\noindent{\bf Acknowledgement} The author is grateful to Alexander Zvonkin for many very helpful comments.


\section{Preliminaries}\label{dessins}

We refer the reader to~\cite{GG, JW, LZ} for background on dessins. Here we summarise a few special aspects of dessins and other topics to be used later.

\subsection{Hurwitz dessins}\label{Hdessins}

Hurwitz showed that if $X$ is a compact Riemann surface of genus $g>1$ then $|{\rm Aut}\,X|\le84(g-1)$. A group $G$ is defined to be a {\em Hurwitz group\/} if $G\cong {\rm Aut}\,X$ for some surface $X$ attaining this bound, or equivalently, $G$ is a non-trivial finite quotient of the triangle group $\Delta:=\Delta(3,2,7)$. In this case there is a bijection, given by $X\cong\mathbb H/N$, between the isomorphism classes of such surfaces $X$ and the normal subgroups $N$ of $\Delta$ with $\Delta/N\cong G$.

Any transitive finite permutation representation $\theta:\Delta\to G\le {\rm Sym}\,\Omega$ of $\Delta$ gives a Hurwitz group $G$. Hence so does any subgroup $M$ of finite index in $\Delta$, where $\theta:\Delta\to G$ is the permutation representation of $\Delta$ on the cosets of $M$, and $N=\ker\theta$ is the core of $M$ in $\Delta$, the intersection of its conjugates.

In such cases $\theta$ or $M$ induces a {\em Hurwitz dessin\/} $\mathcal M$ on the Riemann surface $\mathbb H/M$; this is a bipartite map with the edges corresponding to the elements of $\Omega$, and the black and white vertices and faces corresponding to the cycles of the images $x, y$ and $z$ in $G$ of the standard generators of $\Delta$, so that $G$ acts as the monodromy group of $\mathcal M$. Similarly the subgroup $N$, which corresponds to the regular representation of $G$, induces a regular Hurwitz dessin $\mathcal N$ on the Hurwitz surface $\mathbb H/N$; this dessin $\mathcal N$ is the Galois cover, or minimal regular cover $\tilde{\mathcal M}$ of $\mathcal M$, with automorphism group ${\rm Aut}\,\mathcal N\cong G$.

If $|\Delta:M|=n$, and $x, y$ and $z$ (of orders $3, 2$ and $7$) have $\alpha, \beta$ and $\gamma$ fixed points on the cosets of $M$, then it follows from a result of Singerman~\cite{Sin} that $M$ has signature $\sigma=(g; 3^{[\alpha]}, 2^{[\beta]}, 7^{[\gamma]})$, where the Riemann--Hurwitz formula implies that $\mathcal M$ has genus
\begin{equation}\label{RH}
g=1+\frac{1}{84}(n-28\alpha-21\beta-36\gamma).
\end{equation}
This means that $M$ has a standard presentation as a Fuchsian group, with generators
\[A_i, B_i\;(i=1,\ldots, g), \;\; X_j\;(j=1,\ldots, \alpha), \;\; Y_k\;(k=1,\ldots,\beta),\;\; Z_l\;(l=1,\ldots,\gamma)\]
and defining relations
\[\prod_i[A_i,B_i]\cdot\prod_jX_j\cdot\prod_kY_k\cdot\prod_lZ_l=X_j^3=Y_k^2=Z_l^7=1.\]
We will also refer to $\sigma$ as the signature of the dessin $\mathcal M$ corresponding to $M$.

The values of $\alpha$, $\beta$ and $\gamma$ can be computed as follows. Any transitive permutation group $G$ can be regarded as permuting the cosets of a point-stabiliser $H$. A simple counting argument shows that the number of fixed points of an element $h\in G$ is $|h^G\cap H|.|C_G(h)|/|H|$, where $h^G$ and $C_G(h)$ denote the conjugacy class and centraliser of $h$ in $G$, so that $|h^G|.|C_G(h)|=|G|$ by the Orbit-Stabiliser Theorem.
Of course, since the periods $3, 2$ and $7$ of $\Delta$ are all prime, the parameters $n, \alpha, \beta$ and $\gamma$ uniquely determine the cycle-structures of $x, y$ and $z$, and hence the signatures of $M$ and $\mathcal D$.


\subsection{Jordan's Theorem}\label{Jordan}

Many Hurwitz dessins of degree $n$ have monodromy group $G$ isomorphic to the alternating group ${\rm A}_n$. To prove this in specific cases, one can often use the following theorem (see~\cite[Theorem~13.9]{Wie}):

\begin{theorem}[Jordan]\label{Jordanthm}
If $G$ is a primitive permutation group of degree $n$ containing a cycle of prime length $l\le n-3$ then $G\ge{\rm A}_n$.
\end{theorem}

We will also use the following extension of Jordan's Theorem~\cite{Jones-14}, where the primality condition is omitted. (The proof uses the classification of finite simple groups; this is the only part of this paper which depends upon it.)

\begin{theorem}\label{JJordanthm}
If $G$ is a primitive permutation group of degree $n$ containing a cycle of length $l\le n-3$ then $G\ge{\rm A}_n$.
\end{theorem}

In either case, if $G$ is the monodromy group of a Hurwitz dessin, then $G$ is a quotient of $\Delta$ and so must be perfect; thus $G\ne{\rm S}_n$ and hence $G={\rm A}_n$.


\subsection{From coset diagrams to dessins}

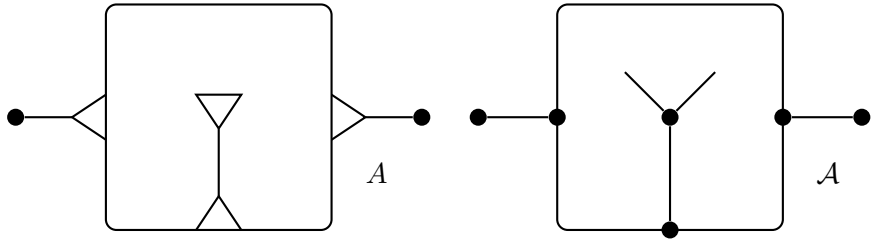
\begin{figure}[h!]
\begin{center}
\begin{tikzpicture}[scale=0.15, inner sep=0.8mm]

\draw [thick, rounded corners] (0,10) to (-10,10) to (-10,-10) to (10,-10) to (10,10) to (0,10);
\draw [thick] (10,2) to (13,0) to (10,-2);
\draw [thick] (-10,2) to (-13,0) to (-10,-2);
\node (c) at (18,0) [shape=circle, fill=black] {};
\node (d) at (-18,0) [shape=circle, fill=black] {};
\draw [thick] (13,0) to (c);
\draw [thick] (-13,0) to (d);
\draw [thick] (-2,-10) to (0,-7) to (2,-10);
\draw [thick] (0,-7) to (0,-1);
\draw [thick] (0,-1) to (-2,2) to (2,2) to (0,-1);
\node at (14,-5)  {$A$};


\draw [thick, rounded corners] (40,10) to (30,10) to (30,-10) to (50,-10) to (50,10) to (40,10);
\node (E) at (50,0) [shape=circle, fill=black] {};
\node (F) at (30,0) [shape=circle, fill=black] {};
\node (C) at (57,0) [shape=circle, fill=black] {};
\node (D) at (23,0) [shape=circle, fill=black] {};
\draw [thick] (C) to (E);
\draw [thick] (D) to (F);
\node (G) at (40,-10) [shape=circle, fill=black] {};
\node (H) at (40,0) [shape=circle, fill=black] {};
\draw [thick] (G) to (H);
\draw [thick] (36,4) to (H) to (44,4);
\node at (54,-5)  {$\mathcal A$};

\end{tikzpicture}

\end{center}
\caption{The diagram $A$ and dessin $\mathcal A$.}
\label{AandA}
\end{figure}

In order to construct permutation representations of $\Delta$, Higman and Conder used the vertices of small triangles to represent $3$-cycles of $X$, with their cyclic order given by the positive (anticlockwise) orientation of the page, and heavy dots for its fixed points; they used long edges to indicate $2$-cycles of $Y$. Any connected diagram formed in this way is a coset diagram for a point-stabiliser $M\le \Delta$. For example the left-hand part of Figure~\ref{AandA} shows the coset diagram $A$ from~\cite{Conder-80}; here $M$ has index $14$, and $X$ and $Y$ have cycle structures $1^{[2]},3^{[4]}$ and $1^{[2]},2^{[6]}$.

Any coset diagram $D$ for $\Delta$ can be converted into a Hurwitz dessin $\mathcal D$ by shrinking each triangle representing a $3$-cycle of $X$ to a single vertex, with the cyclic order giving the local orientation of the surface, and then adding a free edge for each fixed point of $Y$. Thus the directed edges of $\mathcal D$ correspond to the vertices of $D$, with the same actions of $X$ and $Y$. The faces of the dessin correspond to the cycles of $Z$. This transformation is illustrated in Figure~\ref{AandA}, which shows the spherical dessin $\mathcal A$ corresponding to the coset diagram $A$. The process can be reversed by truncating a Hurwitz dessin $\mathcal D$, so that each trivalent vertex is replaced with a small triangle, then removing the free edges and ignoring the faces, to give the coset diagram $D$. 

Conder's diagrams in~\cite{Conder-80} all have bilateral symmetry, so they can be regarded as coset diagrams for subgroups of the extended triangle group
\[\Delta[3,2,7]=\langle X, Y, T\mid X^3=Y^2=(XY)^7=(XT)^2=(YT)^2=1\rangle,\]
with $T$ permuting vertices by reflection in the vertical axis. The motivation for this was that Conder wanted to realise ${\rm S}_n$ and ${\rm A}_n$ as quotients of $\Delta[3,2,7]$ for all sufficiently large $n$. Here we do not have this ambition, so we will not restrict our diagrams and dessins to those with bilateral symmetry.

The dessins $\mathcal A,\ldots, \mathcal N$ corresponding to Conder's basic diagrams $A,\ldots, N$ are all drawn in~\cite{JP}, whereas only $\mathcal A$, $\mathcal B$, $\mathcal C$, $\mathcal F$ and $\mathcal G$ are used in this paper; their basic properties are listed in the Appendix.


\subsection{Macbeath--Hurwitz curves}\label{Macbeath}

Being perfect, every Hurwitz group is a covering of a non-abelian finite simple group, which is also a Hurwitz group. Up to any given finite order, most non-abelian finite simple groups have the form ${\rm PSL}_2(q)$ for prime-powers $q$, so it is natural that many of our examples will also have this form. 

\begin{theorem}[Macbeath~\cite{Macb}]\label{Macbthm}
${\rm PSL}_2(q)$ is a Hurwitz group if and only if one of the following conditions is satisfied:
\begin{enumerate}
\item $q=7$ or $q=p^3$ for a prime $p\equiv \pm 2$ or $\pm 3$ {\rm mod}~$(7)$, with one associated Hurwitz curve;
\item $q$ is a prime $p\equiv \pm 1$ {\rm mod}~$(7)$, with three associated Hurwitz curves distinguished by the choice of $\pm {\rm tr}(z)$.
\end{enumerate}
\end{theorem}

In case (1) the uniqueness of the corresponding regular dessin shows that it is defined over $\mathbb Q$. Streit~\cite{Streit} showed that in case~(2) the three dessins are defined over the real cyclotomic field $\mathbb Q(\cos(2\pi/7))={\mathbb Q}(e^{2\pi i/7})\cap{\mathbb R}$, and that for each $q$ they form an orbit under the Galois group ${\rm C}_3$ of that field. This orbit consists of one reflexible dessin and a chiral pair.

For odd $q$, ${\rm PSL}_2(q)$ has order $q(q^2-1)/2$, so the genus of the corresponding curve or curves is
\[g=1+\frac{q(q^2-1)}{84}.\]
On the other hand $|{\rm PSL}_2(8)|=504$, so in this case the curve has genus~$7$.


\section{Some small dessins}\label{small}

We now introduce some small Hurwitz dessins to be used later as examples.

\begin{figure}[h!]
\begin{center}
\begin{tikzpicture}[scale=0.4, inner sep=0.8mm]

\node (a) at (-5,0) [shape=circle, fill=black] {};
\node (b) at (-10,0) [shape=circle, fill=black] {};
\node (c) at (-12,2) [shape=circle, fill=black] {};
\draw [thick] (-3,2) to (a) to (-3,-2);
\draw [thick] (-12,-2) to (b) to (c);
\draw [thick] (c) to (b) to (a);

\node at (-7,-3)  {$\mathcal S$};


\node (a') at (5,0) [shape=circle, fill=black] {};
\node (b') at (10,0) [shape=circle, fill=black] {};
\node (c') at (12,2) [shape=circle, fill=black] {};
\draw [thick] (3,2) to (a') to (3,-2);
\draw [thick] (12,-2) to (b') to (c');
\draw [thick] (c') to (b') to (a');

\node at (7,-3)  {$\overline\mathcal S$};

\end{tikzpicture}

\end{center}
\caption{Dessins $\mathcal S$ and $\overline\mathcal S$.}
\label{dessinsSSbar}
\end{figure}
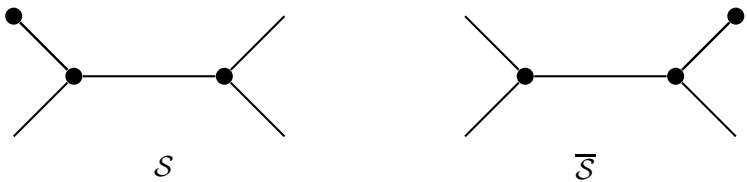

\medskip

\noindent{\bf Example} The dessin $\mathcal S$ and its mirror image $\overline\mathcal S$, shown in Figure~\ref{dessinsSSbar}, are the smallest Hurwitz dessins, having degree $n=7$. They both have genus $0$. In each case the monodromy group is ${\rm PGL}_3(2)$, in its natural action on the points or the lines of the Fano plane ${\mathbb P}^2({\mathbb F}_2)$; these are equivalent to the actions of the isomorphic group ${\rm PSL}_2(7)$ on its two conjugacy classes of subgroups $H\cong {\rm S}_4$. These dessins correspond to two conjugacy classes of subgroups $M$ of $\Delta$, of index $7$ and signature $(0;2,2,2,3)$.

\medskip

\begin{figure}[h!]
\begin{center}
\begin{tikzpicture}[scale=0.4, inner sep=0.8mm]

\draw [thick] (2,-2) arc (0:360:2);
\node (a) at (0.75,-2) [shape=circle, fill=black] {};
\node (a') at (-0.75,-2) [shape=circle, fill=black] {};
\node (b) at (0,0) [shape=circle, fill=black] {};
\node (c) at (0.75,-3.85) [shape=circle, fill=black] {};
\node (c') at (-0.75,-3.85) [shape=circle, fill=black] {};
\node (d) at (0,2) [shape=circle, fill=black] {};

\draw [thick] (c) to (a);
\draw [thick] (c') to (a');
\draw [thick] (b) to (d);
\draw [thick] (-1,3) to (d) to (1,3);

\node at (3.5,-3) {$\mathcal A$};


\draw [thick] (14,0) arc (0:360:3);
\node (A) at (8,0) [shape=circle, fill=black] {};
\node (A') at (14,0) [shape=circle, fill=black] {};
\node (B) at (11,-3) [shape=circle, fill=black] {};
\node (C) at (11,3) [shape=circle, fill=black] {};
\node (D) at (11,1) [shape=circle, fill=black] {};

\draw [thick] (6.5,0) to (A);
\draw [thick] (15.5,0) to (A');
\draw [thick] (B) to (11,-4.5);
\draw [thick] (C) to (D);
\draw [thick] (12,0) arc (0:360:1);

\node at (15,-3) {$\mathcal B$};


\draw [thick] (24,-2) arc (0:360:2);
\node (a2) at (24,-2) [shape=circle, fill=black] {};
\node (a2') at (20,-2) [shape=circle, fill=black] {};
\node (b2) at (22,0) [shape=circle, fill=black] {};
\node (c2) at (23,-0.3) [shape=circle, fill=black] {};
\node (c2') at (21,-0.3) [shape=circle, fill=black] {};
 
\draw [thick] (c2) to (23,-2);
\draw [thick] (c2') to (21,-2);
\draw [thick] (b2) to (22,1.5);

\node (d2) at (22,3) [shape=circle, fill=black] {};
\draw [thick, rounded corners] (a2) to (25,-2) to (25,3) to (19,3) to (19,-2) to (a2');

\node (e2) at (22,5) [shape=circle, fill=black] {};
\draw [thick] (d2) to (e2);
\draw [thick] (23,6) to (e2) to (21,6);

\node at (26,-3)  {$\mathcal C$};

\end{tikzpicture}

\end{center}
\caption{Dessins $\mathcal A$, $\mathcal B$ and $\mathcal C$.}
\label{dessinsABC}
\end{figure}
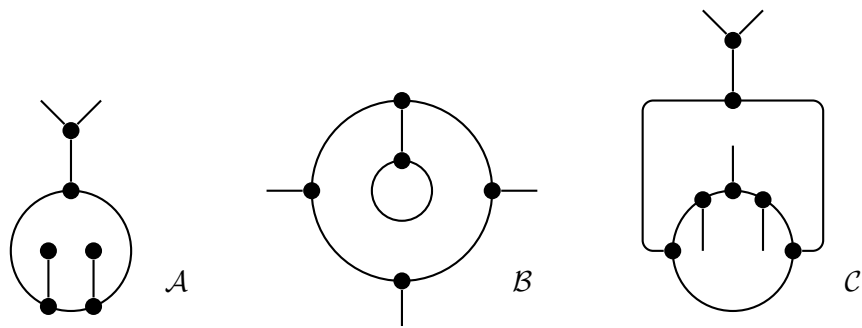

\noindent{\bf Example} Figure~\ref{dessinsABC} shows three Hurwitz dessins $\mathcal A$, $\mathcal B$ and $\mathcal C$, of genus $0$ and of degrees $14$, $15$ and $21$, based on Conder's coset diagrams $A$, $B$ and $C$ in~\cite{Conder-80}. Note that $\mathcal A$ has been redrawn from Figure~\ref{AandA} so that the free edges are now in the outer face; the reason for this will become clear in the next section, when we consider joining operations (see Figure~\ref{joining}).

The monodromy group of $\mathcal A$ is ${\rm PSL}_2(13)$ in its natural representation of degree 14 on the projective line ${\mathbb P}^1({\mathbb F}_{13})$; the point-stabilisers are the Sylow 13-normalisers $H\cong {\rm C}_{13}\rtimes {\rm C}_6$. One can take $x: t\mapsto 1/(1-t)$, $y: t\mapsto -t$ and $z: t\mapsto (t+1)/t$ as standard generators. This dessin corresponds to a conjugacy class of subgroups $M$ of index $14$ in $\Delta$; since $\mathcal A$ has genus $0$, and $x$ and $y$ each have two fixed points while $z$ has none, these subgroups $M$ are quadrilateral groups, with signature $(0; 2, 2, 3, 3)$.  We have $\mathcal A\cong\widetilde{\mathcal A}/H$ where $\widetilde{\mathcal A}$ is the Galois cover of $\mathcal A$, corresponding to the core $K$ of $M$ in $\Delta$; this is the unique reflexible dessin in the Galois orbit of three Macbeath--Hurwitz dessins of genus 14 with automorphism group  ${\rm PSL}_2(13)$. (See~\cite{JZ} for further details of this Galois orbit and the corresponding quotient dessins.)

In $\mathcal B$, the commutator $[x,y]$ has cycle structure $3,5,7$, so the monodromy group $G$ contains a $3$-cycle, a $5$-cycle and a $7$-cycle. A transitive group of degree $15$ with an element of order $7$ must be primitive (otherwise it would be contained in a wreath product ${\rm S}_5\wr{\rm S}_3$ or ${\rm S}_3\wr{\rm S}_5$, of order coprime to $7$), and a primitive group containing a $3$-cycle contains the alternating group (see~\cite[Theorem~13.3]{Wie}). Since $x, y$ and $z$ are even, $G={\rm A}_{15}$. (As Conder showed in~\cite{Conder-80}, this is the smallest alternating group which is a Hurwitz group.) In this case the corresponding subgroups $M$ of $\Delta$ have signature $(0; 2, 2, 2, 7)$.

The monodromy group of $\mathcal C$ is ${\rm PGL}_3(2)$, isomorphic to ${\rm PSL}_2(7)$, in its imprimitive action of degree $21$ on the flags (incident point-line pairs) of the Fano plane, or equivalently its action by conjugation on its involutions. The point-stabilisers are the Sylow $2$-subgroups $H\cong {\rm D}_4$. The corresponding subgroups $M$ of $\Delta$ have signature $(0; 2,2,2,2,2)$.

 
\section{Joins for Hurwitz dessins}\label{joins}

From now on, unless stated otherwise, all dessins will be Hurwitz dessins, that is, finite oriented hypermaps of type $(3,2,7)$. In this section we will consider the handles and joins introduced by Higman and used by Conder, reinterpreting them as operations on dessins.

\subsection{$y$-handles and $y$-joins}\label{kjoins}

\begin{figure}[h!]
\begin{center}
\begin{tikzpicture}[scale=0.5, inner sep=0.8mm]

\node (1) at (-8,3) [shape=circle, fill=black] {};
\draw [thick] (1) to (-8,1);
\draw [thick, dashed] (-8,1) to (-8,-1);
\draw [thick] (-10,5) to (1) to (-6,5);
\node at (-6,5.5) {$a$};
\node at (-10,5.5) {$b$};
\node at (-5.5,-1) {$k=1$};


\node (2) at (-2,2) [shape=circle, fill=black] {};
\node (3) at (2,2) [shape=circle, fill=black] {};
\draw [thick] (2) to (-2,5);
\draw [thick] (3) to (2,5);
\draw [thick] (-3,2) to (3,2);
\draw [thick, dashed] (-4,2) to (-3,2);
\draw [thick, dashed] (4,2) to (3,2);
\node at (2,5.5) {$a$};
\node at (-2,5.5) {$b$};
\node at (0,-1) {$k=2$};


\node (4) at (8,2) [shape=circle, fill=black] {};
\node (5) at (12,2) [shape=circle, fill=black] {};
\node (6) at (10,2) [shape=circle, fill=black] {};
\draw [thick] (4) to (8,5);
\draw [thick] (5) to (12,5);
\draw [thick] (7,2) to (13,2);
\draw [thick, dashed] (6,2) to (7,2);
\draw [thick, dashed] (14,2) to (13,2);
\draw [thick] (6) to (10,0);
\draw [thick, dashed] (10,0) to (10,-1);
\node at (12,5.5) {$a$};
\node at (8,5.5) {$b$};
\node at (8,-1) {$k=3$};

\end{tikzpicture}

\end{center}
\caption{$(k)$-handles for $k=1, 2, 3$.}
\label{k-handles}
\end{figure}
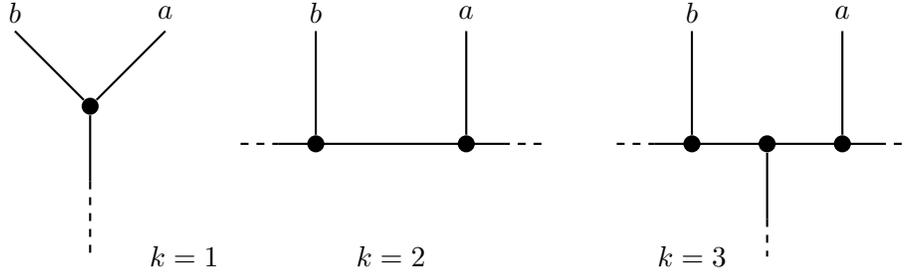

The joins used by Higman and Conder are based on fixed points of the generating involution $Y$  for $\Delta$ (or $X$ in their notation; the symbols are transposed here for more convenient explanation of links with dessins). Let $\mathcal D$ be a Hurwitz dessin with monodromy group
\[G=\langle x, y, z\mid x^3=y^2=z^7=xyz=\cdots=1\rangle.\]
This gives a faithful transitive permutation representation of $G$. As usual, cycles of $x$, $y$ and $z$ are represented as vertices, edges and faces. If $y$ has two fixed points (free edges) in the same face, we can uniquely label them $a$ and $b$ so that $b=ax$, $axyx$ or $axyxyx$, that is, $b=az^{1-k}x$ for some $k=1, 2$ or $3$; the ordered pair $(a,b)$ is then called a $(k)$-{\sl handle}. These are illustrated in Figure~\ref{k-handles}. If we do not wish to specify the value of $k$, we will sometimes call these $y$-{\sl handles}, to distinguish them from the $x$-handles to be defined later.

\noindent{\bf Example} In each of the dessins $\mathcal S$ and $\overline\mathcal S$ in Figure~\ref{dessinsSSbar}, $Y$ has three fixed points. One ordered pair of them forms a $(1)$-handle, and another forms a $(2)$-handle; there are no $(3)$-handles.

\medskip

\noindent{\bf Example} In the dessin $\mathcal A$ in Figure~\ref{dessinsABC}, $Y$ has two fixed points, forming a $(1)$-handle. In $\mathcal B$, $Y$ has three fixed points; two ordered pairs form $(2)$-handles, and a third forms a $(3)$-handle. In $\mathcal C$, $Y$ has five fixed points, providing a $(1)$-handle and a $(3)$-handle.

\begin{lemma}\label{handlecover}
Let $\mathcal D\to\overline{\mathcal D}$ be a covering of Hurwitz dessins. Then any $(k)$-handle in $\mathcal D$ projects onto a $(k)$-handle in $\overline{\mathcal D}$.
\end{lemma}

\noindent{\sl Proof.} If $ay=a$, $by=b$ and $b=az^{1-k}x$ in $\mathcal D$, then the images $\overline a$ and $\overline b$ of $a$ and $b$ in $\overline{\mathcal D}$ satisfy the corresponding equations. If $\overline a=\overline b$ then the covering reduces the number of sides of the face containing the handle $(a,b)$, so its image has just one side. Thus $Z$ fixes $\overline a$, and hence so does $\langle Y, Z\rangle=\Delta$, contradicting the transitivity of $\Delta$. Hence $\overline a\ne\overline b$, so $(\overline a,\overline b)$ is a $(k)$-handle. \hfill$\square$

\medskip

The following partial converse is obvious:

\begin{lemma}\label{lifting}
Let $\tilde{\mathcal D}\to{\mathcal D}$ be a $d$-sheeted covering of Hurwitz dessins. Then any $(k)$-handle $(a,b)$ in $\mathcal D$ lifts to $d$ disjoint $(k)$-handles in $\tilde{\mathcal D}$, provided the covering is unbranched over $a$ and $b$.  \hfill$\square$
\end{lemma}

\medskip

Dessins with arbitrarily many $(k)$-handles can be constructed by applying Lemma~\ref{lifting} to a dessin $\mathcal D$ with at least one $(k)$-handle $(a,b)$, where the covering $\tilde{\mathcal D}\to\mathcal D$ corresponds to a suitable index $d$ inclusion $\tilde M\le M$ of map subgroups. To avoid branching over $a$ and $b$ we require the elliptic generators of order $2$ of $M$ corresponding to the fixed points $a$ and $b$ of $Y$ to be elements of $\tilde{M}$. In all except a few small cases, such subgroups $\tilde M$ exist.

\medskip

\noindent{\bf Example} Let $\mathcal D=\mathcal C$, so that $M$ has signature $(0;2^{[5]})$. Let $Y_1$ and $Y_2$ be the elliptic generators corresponding to the $(1)$-handle, $Y_3$ and $Y_4$ those corresponding to the $(3)$-handle, and let $Y_5$ be the fifth, chosen so that $Y_1\ldots Y_5=1$. For any integer $d\ge 2$ define an epimorphism
\[\theta:M\to{\rm D}_d=\langle u, v\mid u^2=v^2=(uv)^d=1\rangle\]
by $Y_1, Y_2\mapsto u$, $Y_3, Y_4\mapsto v$, and $Y_5\mapsto 1$. Then $\tilde{M}=\theta^{-1}(\langle u\rangle)$ has index $d$ in $M$, and contains $Y_1$ and $Y_2$, so it corresponds to a covering of $\mathcal C$ with $d$ disjoint $(1)$-handles, while $\theta^{-1}(\langle v\rangle)$ corresponds to a covering with $d$ disjoint $(3)$-handles. Similarly, examples with arbitrarily many disjoint $(2)$-handles can be constructed as coverings of $\mathcal B$.

\medskip

\begin{figure}[h]
\begin{center}
\begin{tikzpicture}[scale=0.5, inner sep=0.8mm]


\draw [thick, dotted] (-3,10) arc (0:360:3);
\node (v2) at (-6,14) [shape=circle, fill=black] {};
\draw [thick] (-7,15) to (v2) to (-5,15);
\draw [thick] (v2) to (-6,13);

\node (a2) at (-4.3,15) {$a'$};
\node (b2) at (-7.6,15) {$b'$};
\node (D2) at (-6,10) {$\mathcal D_2$};


\draw [thick, dotted] (9,10) arc (0:360:3);
\node (v2') at (6,14) [shape=circle, fill=black] {};
\draw [thick] (v2') to (6,13);

\node (D2) at (6,10) {$\mathcal D_2$};


\draw [thick, dotted] (-3,0) arc (0:360:3);
\node (v1) at (-6,4) [shape=circle, fill=black] {};
\draw [thick] (-7,5) to (v1) to (-5,5);
\draw [thick] (v1) to (-6,3);

\node (a1) at (-4.3,5) {$a$};
\node (b1) at (-7.6,5) {$b$};
\node (D1) at (-6,0) {$\mathcal D_1$};


\draw [thick, dotted] (9,0) arc (0:360:3);
\node (v1') at (6,4) [shape=circle, fill=black] {};
\draw [thick] (v1') to (6,3);

\node (D1) at (6,0) {$\mathcal D_1$};


\draw [rounded corners, thick] (v1') to (10,4) to (10,14) to (v2');
\draw [rounded corners, thick] (v1') to (2,4) to (2,14) to (v2');

\end{tikzpicture}

\end{center}
\caption{Construction of $\mathcal D_1(1)\mathcal D_2$}
\label{joining}
\end{figure}
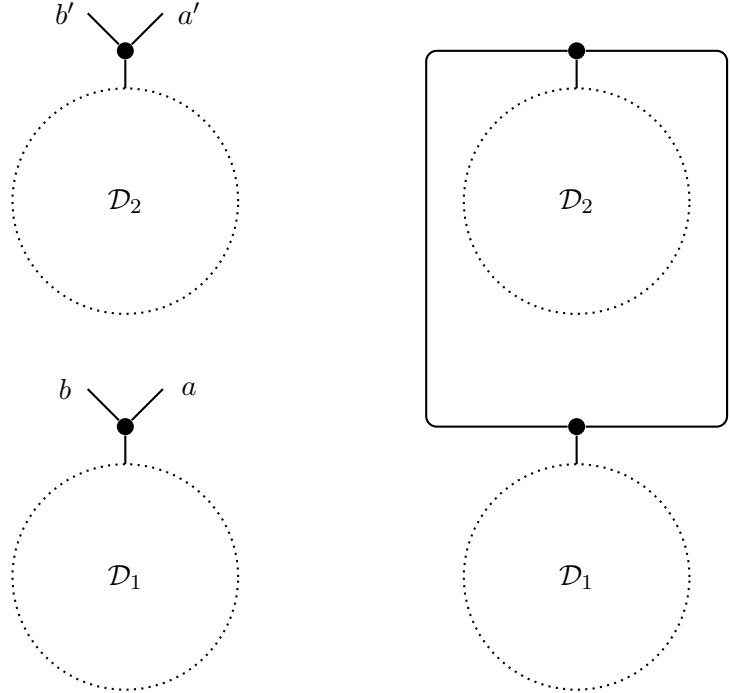

We now define joins of dessins. Suppose that Hurwitz dessins $\mathcal D_i\;(i=1, 2)$ of degree $d_i$ and genus $g_i$ have $(k)$-handles $(a_i, b_i)$ for the same value of $k$. We can form a $(k)$-{\sl join} $\mathcal D_1(k)\mathcal D_2$ by defining the cycles of $x$ and $y$ to be those they have on $\mathcal D_1$ or on $\mathcal D_2$, except that the four fixed points $a_i, b_i$ of $y$ are replaced with two 2-cycles $(a_1,a_2)$ and $(b_1,b_2)$.  (This is called a $k$-{\sl composition} in~\cite{Conder-80}, but we prefer to reserve the word `composition' for a different method of combining dessins, by composing their Bely\u\i\/ functions.) As with handles, if we do not wish to specify the value of $k$ we sometimes will call these $y$-{\sl joins}, to distinguish them from $x$-joins defined later. This operation is illustrated in Figure~\ref{joining} in the case $k=1$; the cases $k=2$ and $3$ are similar. Topologically this is a connected sum operation, where the two dessins are joined across cuts between the free ends of the half-edges representing the fixed points $a_i$ and $b_i$, with the cuts staying in the faces containing the two $k$-handles. If either dessin $\mathcal D_i$ has more than one $(k)$-handle, one may need to specify the choice of handles in order to define their $(k)$-join uniquely.

\medskip

The following simple result records the additivity properties of joins.

\begin{theorem}\label{additivity}
If Hurwitz dessins $\mathcal D_i\;(i=1,2)$ of degree $n_i$ and signature $(g_i;3^{[\alpha_i]},2^{[\beta_i]},7^{[\gamma_i]})$ have a $(k)$-join $\mathcal D$, then $\mathcal D$ has degree $n$ and signature  $(g;3^{[\alpha]},2^{[\beta]},7^{[\gamma]})$ where
\[n=n_1+n_2,\; g=g_1+g_2, \; \alpha=\alpha_1+\alpha_2, \; \beta=\beta_1+\beta_2-4\;\hbox{and}\; \gamma=\gamma_1+\gamma_2.\]
\end{theorem}

\noindent{\sl Proof.} The equations for $n, \alpha, \beta$ and $\gamma$ follow immediately from the definition of a join, and that for $g$ follows from applying the Riemann--Hurwitz formula~(\ref{RH}) to these three dessins. \hfill$\square$

\begin{figure}[h!]
\begin{center}
\begin{tikzpicture}[scale=0.5, inner sep=0.8mm]

\draw [thick] (2,11) arc (0:360:2);
\node (a1) at (0.75,11) [shape=circle, fill=black] {};
\node (a1') at (-0.75,11) [shape=circle, fill=black] {};
\node (b1) at (0,13) [shape=circle, fill=black] {};
\node (c1) at (0.75,9.15) [shape=circle, fill=black] {};
\node (c1') at (-0.75,9.15) [shape=circle, fill=black] {};
\node (d1) at (0,15) [shape=circle, fill=black] {};

\draw [thick] (c1) to (a1);
\draw [thick] (c1') to (a1');
\draw [thick] (b1) to (d1);
\draw [thick] (-1,16) to (d1) to (1,16);

\node at (3,14.5) {$\mathcal A$};


\draw [thick] (2,0) arc (0:360:2);
\node (a) at (2,0) [shape=circle, fill=black] {};
\node (a') at (-2,0) [shape=circle, fill=black] {};
\node (b) at (0,2) [shape=circle, fill=black] {};
\node (c) at (1,1.7) [shape=circle, fill=black] {};
\node (c') at (-1,1.7) [shape=circle, fill=black] {};
 
\draw [thick] (c) to (1,0);
\draw [thick] (c') to (-1,0);
\draw [thick] (b) to (0,3.5);

\node (d) at (0,5) [shape=circle, fill=black] {};
\draw [thick, rounded corners] (a) to (3,0) to (3,5) to (-3,5) to (-3,0) to (a');

\node (e) at (0,7) [shape=circle, fill=black] {};
\draw [thick] (d) to (e);
\draw [thick] (1,8) to (e) to (-1,8);

\node at (3,6.5)  {$\mathcal C$};


\draw [thick] (12,11) arc (0:360:2);
\node (A1) at (10.75,11) [shape=circle, fill=black] {};
\node (A1') at (9.25,11) [shape=circle, fill=black] {};
\node (B1) at (10,13) [shape=circle, fill=black] {};
\node (C1) at (10.75,9.15) [shape=circle, fill=black] {};
\node (C1') at (9.25,9.15) [shape=circle, fill=black] {};
\node (D1) at (10,15) [shape=circle, fill=black] {};

\draw [thick] (C1) to (A1);
\draw [thick] (C1') to (A1');
\draw [thick] (B1) to (D1);


\draw [thick] (12,0) arc (0:360:2);
\node (A) at (12,0) [shape=circle, fill=black] {};
\node (A') at (8,0) [shape=circle, fill=black] {};
\node (B) at (10,2) [shape=circle, fill=black] {};
\node (C) at (11,1.7) [shape=circle, fill=black] {};
\node (C') at (9,1.7) [shape=circle, fill=black] {};

\draw [thick] (C) to (11,0);
\draw [thick] (C') to (9,0);
\draw [thick] (B) to (10,3.5);

\node (D) at (10,5) [shape=circle, fill=black] {};
\draw [thick, rounded corners] (A) to (13,0) to (13,5) to (7,5) to (7,0) to (A');

\node (E) at (10,7) [shape=circle, fill=black] {};
\draw [thick] (D) to (E);

\draw [thick, rounded corners] (D1) to (7,15) to (7,7) to (13,7) to (13,15) to (D1);

\node at (15,6)  {$\mathcal A(1)\mathcal C$};

\end{tikzpicture}

\end{center}
\caption{Construction of $\mathcal A(1)\mathcal C$.}
\label{dessinA1C}
\end{figure}
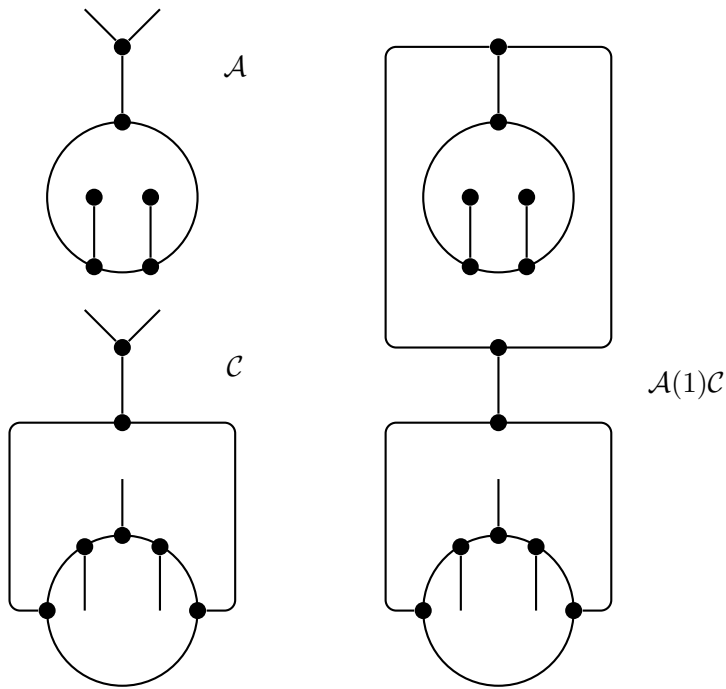

\medskip

\noindent{\bf Example.} The dessins $\mathcal A$ and $\mathcal C$ in Figure~\ref{dessinsABC} each have a unique $(1)$-handle, so we can form the join $\mathcal A(1)\mathcal C$ of degree $14+21=35$, as shown in Figure~\ref{dessinA1C}. This dessin has monodromy group $G={\rm A}_{35}$. To prove this, we first show that $G$ is primitive. If not, it has $a$ blocks of size $b$, so $G$ is embedded in the wreath product ${\rm S}_b\wr{\rm S}_a=({\rm S}_b)^a\rtimes{\rm S}_a$, where $ab=35$ and hence $\{a, b\}=\{5, 7\}$. We cannot have $a=5$ since $G$ would then act on the blocks as a Hurwitz subgroup of ${\rm S}_5$, so $a=7$ and $b=5$. This implies that the Sylow $7$-subgroups of $G$, which are mutually conjugate, are isomorphic to ${\rm C}_7$, so all elements of order $7$ have the same cycle structure. However, inspection of Figure~\ref{dessinA1C} shows that $w:=[x,y]$ has cycle structure $1^{[2]},2^{[2]},4^{[2]},21$, so $w^{12}$ has cycle structure $1^{[14]},7^{[3]}$, whereas for $z$ it is $7^{[5]}$. Thus $G$ is primitive. Now $w^4$ is a $21$-cycle, so $G\ge {\rm A}_{35}$ by Theorem~\ref{JJordanthm}, with equality since $G$ is perfect.
(This has been confirmed using GAP.) This example illustrates how the monodromy group of a join (in this case of order approximately $5\cdot17\times 10^{39}$) can be much larger than those of the two factors (here of orders $1092$ and $168$). (In principle, it could be smaller, but it is hard to think of an example.)

\medskip

\begin{figure}[h!]
\begin{center}
\begin{tikzpicture}[scale=0.4, inner sep=0.8mm]

\node (a) at (0,0) [shape=circle, fill=black] {};
\node (b) at (-5,0) [shape=circle, fill=black] {};
\node (c) at (-7,2) [shape=circle, fill=black] {};
\draw [thick] (-7,-2) to (b) to (c);
\draw [thick] (c) to (b) to (a);

\draw [thick] (5,0) arc (0:360:2.5);


\node (a') at (5,0) [shape=circle, fill=black] {};
\node (b') at (10,0) [shape=circle, fill=black] {};
\node (c') at (12,-2) [shape=circle, fill=black] {};
\draw [thick] (12,2) to (b') to (c');
\draw [thick] (c') to (b') to (a');

\end{tikzpicture}

\end{center}
\caption{Joining copies of $\mathcal S$ incorrectly.}
\label{wrongSS}
\end{figure}
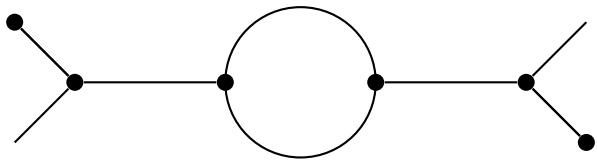

\noindent{\bf Warning.} When making a $y$-join, it is important that $a_1$ should be paired with $a_2$, and $b_1$ with $b_2$, rather than $a_1$ with $b_2$ and $a_2$ with $b_1$. For example, the dessin $\mathcal S$ in Figure~\ref{dessinsSSbar} has a $(1)$-handle, so we can form a Hurwitz dessin $\mathcal S(1)\mathcal S$ (as we will see in Figure~\ref{dessinsS*Sbar*}), but if we choose the wrong pairing for the fixed points of $y$ in the two copies of $\mathcal S$, as in Figure~\ref{wrongSS}, we obtain a dessin of type $(3,2,12)$ rather than $(3,2,7)$, with faces of valency $2$ and $12$; its monodromy group, of order $2688=2^7.3.7$, is solvable, so it cannot be a Hurwitz group. More generally, applying such `incorrect' pairings to Hurwitz dessins gives dessins which have quotients of the modular group
\[\Gamma={\rm PSL}_2(\mathbb Z)=\Delta(3,2,\infty)=\langle X, Y\mid X^3=Y^2=1\rangle,\]
rather than $\Delta(3,2,7)$, as monodromy groups. We will consider this in more detail in \S\ref{modular}, devoted to $\Gamma$.

It should be clear from the definition that joining is commutative and associative, in the sense that there are obvious isomorphisms $\mathcal D_1(k)\mathcal D_2\cong\mathcal D_2(k)\mathcal D_1$ and $\mathcal D_1(k)(\mathcal D_2(k')\mathcal D_3)\cong(\mathcal D_1(k)\mathcal D_2)(k')\mathcal D_3$ whenever these dessins are well-defined. In the case of associativity this requires the chosen $(k)$- and $(k')$-handles in $\mathcal D_2$ to be mutually disjoint. The following lemma ensures that, as long as we do not use copies of $\mathcal B$, $\mathcal S$ or $\overline\mathcal S$, distinct $y$-handles in the same dessin will be in distinct faces, so they will be mutually disjoint and can therefore be used independently of each other:

\begin{lemma}
The following conditions on a Hurwitz dessin $\mathcal D$ are equivalent:
\begin{enumerate}
\item $\mathcal D$ has two $y$-handles with a common fixed point;
\item $\mathcal D$ has two $y$-handles in the same face;
\item $\mathcal D$ is isomorphic to $\mathcal B$, $\mathcal S$ or $\overline\mathcal S$.
\end{enumerate}
\end{lemma}

\noindent{\sl Proof.} The implication $(1)\Rightarrow(2)$ is obvious, and $(3)\Rightarrow(1)$ is evident from Figures~\ref{dessinsSSbar} and \ref{dessinsABC}. The implication $(2)\Rightarrow(3)$ follows from straightforward case-by-case analysis, using the fact that a $(k)$-handle in a face uses four, three or four of its seven sides as $k=1, 2$ or $3$ (see Figure~\ref{k-handles}). \hfill$\square$


\subsection{The double cover $\mathcal D(k)\mathcal D$}

If a dessin $\mathcal D$ has a $(k)$-handle, one can form the join $\mathcal D^*=\mathcal D(k)\mathcal D$, a double cover of $\mathcal D$ branched over two points, namely those corresponding to the fixed points $a$ and $b$ of $Y$ in the chosen handle. 
There is an obvious automorphism of order $2$ of $\mathcal D^*$, interchanging the two copies of $\mathcal D$; the quotient is isomorphic to $\mathcal D$. The double covering $\mathcal D^*\to\mathcal D$ shows that the monodromy group $G^*$ of $\mathcal D^*$ must be imprimitive, with blocks of size $2$ permuted as the monodromy group $G$ of $\mathcal D$, so it is isomorphic to a subgroup of the wreath product ${\rm S}_2\wr G$.

If $\mathcal D$ has degree $n$ and signature $(g;3^{[\alpha]},2^{[\beta]},7^{[\gamma]})$, then it follows from Theorem~\ref{additivity} that $\mathcal D^*$ has degree $n^*=2n$ and signature
\[(g^*;3^{[\alpha^*]},2^{[\beta^*]},7^{[\gamma^*]})=(2g;3^{[2\alpha]},2^{[2\beta-4]},7^{[2\gamma]}).\]

Since $M^*$ has index 2 in $M$, it contains the subgroup $M'M^2$ generated by the commutators and squares of the elements of $M$. Now the standard presentation of $M$ (see \S\ref{Hdessins}) shows that $M/M'M^2$ is an elementary abelian $2$-group of rank $r=2g+\beta-1$, so $M$ has $2^r-1$ subgroups of index $2$, each corresponding to a double covering of $\mathcal D$. (For instance, $\mathcal D^*$ corresponds to the normal closure in $M$ of all its standard generators except the two involutions $Y_k$ corresponding to the fixed points of $Y$ in the chosen handle.) This immediately leads to the following:


\begin{lemma}\label{uniqueD*}
Suppose that a planar Hurwitz dessin $\mathcal D$ has exactly two fixed points for $Y$, forming a $(k)$-handle; then $\mathcal D^*$ is the unique Hurwitz dessin which is a double covering of $\mathcal D$. If, in addition, the point stabiliser $H$ in the monodromy group $G$ of $\mathcal D$ has a subgroup $H^*$ of index $2$, then $\mathcal D^*$ also has monodromy group $G$, with point stabiliser $H^*$. \hfill$\square$
\end{lemma}

However, if $g>0$ or $\beta>2$ then $\mathcal D^*$ is one of several double coverings of $\mathcal D$, and its monodromy group could be $G$ or a proper covering group of $G$ contained in ${\rm S}_2\wr G$ (as must happen if $H$ has no subgroup of index $2$).

\medskip



\begin{figure}[h!]
\begin{center}
\begin{tikzpicture}[scale=0.25, inner sep=0.8mm]

\draw [thick] (-5,5) arc (0:360:5);
\node (A) at (-10,10) [shape=circle, fill=black] {};
\node (B) at (-10,6) [shape=circle, fill=black] {};
\node (C) at (-10,0) [shape=circle, fill=black] {};
\node (D) at (-10,-4) [shape=circle, fill=black] {};
\node (E) at (-12,4) [shape=circle, fill=black] {};
\node (F) at (-12,-6) [shape=circle, fill=black] {};

\draw [thick] (A) to (B);
\draw [thick] (C) to (D);
\draw [thick] (E) to (B) to (-8,4);
\draw [thick] (F) to (D) to (-8,-6);

\node at (-4,-2) {$\mathcal S(1)\mathcal S$};


\draw [thick] (15,5) arc (0:360:5);
\node (A') at (10,10) [shape=circle, fill=black] {};
\node (B') at (10,6) [shape=circle, fill=black] {};
\node (C') at (10,0) [shape=circle, fill=black] {};
\node (D') at (10,-4) [shape=circle, fill=black] {};
\node (E') at (12,4) [shape=circle, fill=black] {};
\node (F') at (12,-6) [shape=circle, fill=black] {};

\draw [thick] (A') to (B');
\draw [thick] (C') to (D');
\draw [thick] (E') to (B') to (8,4);
\draw [thick] (F') to (D') to (8,-6);

\node at (16,-2) {$\overline\mathcal S(1)\overline\mathcal S$};


\draw [thick] (-5,-15) arc (0:360:5);
\node (a) at (-12,-10.5) [shape=circle, fill=black] {};
\node (b) at (-8,-10.5) [shape=circle, fill=black] {};
\node (c) at (-8,-15) [shape=circle, fill=black] {};
\node (d) at (-12,-19.5) [shape=circle, fill=black] {};
\node (e) at (-8,-19.5) [shape=circle, fill=black] {};
\node (f) at (-8,-24) [shape=circle, fill=black] {};

\draw [thick] (b) to (c);
\draw [thick] (e) to (f);
\draw [thick] (a) to (-12,-15);
\draw [thick] (d) to (-12,-24);

\node at (-4,-21.5) {$\mathcal S(2)\mathcal S$};


\draw [thick] (15,-15) arc (0:360:5);
\node (a') at (12,-10.5) [shape=circle, fill=black] {};
\node (b') at (8,-10.5) [shape=circle, fill=black] {};
\node (c') at (8,-15) [shape=circle, fill=black] {};
\node (d') at (12,-19.5) [shape=circle, fill=black] {};
\node (e') at (8,-19.5) [shape=circle, fill=black] {};
\node (f') at (8,-24) [shape=circle, fill=black] {};

\draw [thick] (b') to (c');
\draw [thick] (e') to (f');
\draw [thick] (a') to (12,-15);
\draw [thick] (d') to (12,-24);

\node at (16,-21.5) {$\overline\mathcal S(2)\overline\mathcal S$};

\end{tikzpicture}

\end{center}
\caption{The dessins $\mathcal S(k)\mathcal S$ and $\overline\mathcal S(k)\overline\mathcal S$.}
\label{dessinsS*Sbar*}
\end{figure}
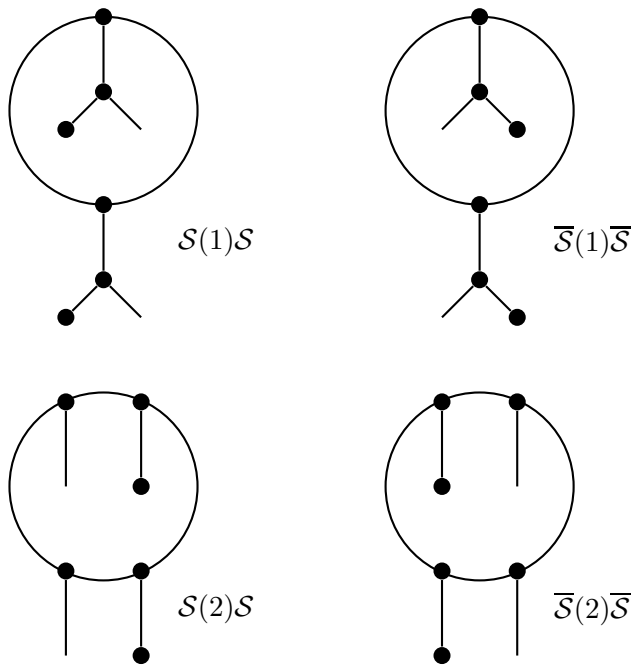

\noindent{\bf Example} Let $\mathcal D=\mathcal S$ or $\overline\mathcal S$ (see Figure~\ref{dessinsSSbar}), with monodromy group $G={\rm PGL}_3(2)$ acting on the points or lines of  the Fano plane, so that $H$ belongs to one of the two conjugacy classes of subgroups of $G$ isomorphic to ${\rm S}_4$. Then $g=0$ but $\beta=3$, so Lemma~\ref{uniqueD*} does not apply. Since $r=2$ there are three double coverings of $\mathcal D$. Taking $k=1$ we obtain the chiral pair of planar dessins $\mathcal S^*=\mathcal S(1)\mathcal S$ and $\overline\mathcal S^*=\overline\mathcal S(1)\overline\mathcal S$ in the top row of Figure~\ref{dessinsS*Sbar*}.  The monodromy group of each is the unique Hurwitz group $G^*$ of genus $17$, a non-split extension of an elementary abelian group of order $8$ by its automorphism group ${\rm GL}_3(2)\cong G$. There is a chiral pair of regular Hurwitz dessins associated with this group, and these two dessins are their quotients. 

Taking $k=2$ we obtain the dessins $\mathcal S(2)\mathcal S$ and $\overline\mathcal S(2)\overline\mathcal S$ in the bottom row of Figure~\ref{dessinsS*Sbar*}. They have the same monodromy group $G$ as $\mathcal S$ and $\overline\mathcal S$, but now in its imprimitive representations on the cosets of the unique subgroup $H^*$ of index $2$ in $H$, isomorphic to ${\rm A}_4$; the two conjugacy classes of such subgroups give the two representations.

The third double covering of $\mathcal D$ in this example does not arise from a $(k)$-join as defined earlier. However, like the coverings $\mathcal D(k)\mathcal D$ for $k=1, 2$ it is branched over two of the three fixed points of $y$ in $\mathcal D$, in this case satisfying $b=axyx^{-1}$ rather than $b=a(xy)^{k-1}x$. Figure~\ref{coversofSSbar} shows the chiral pair of Hurwitz dessins resulting from taking $\mathcal D=\mathcal S$ and $\overline\mathcal S$; their monodromy group is the Hurwitz group $G^*$ of genus $17$ associated earlier with $\mathcal D(1)\mathcal D$.

\begin{figure}[h!]
\begin{center}
\begin{tikzpicture}[scale=0.25, inner sep=0.8mm]

\draw [thick] (-5,-15) arc (0:360:5);
\node (a) at (-12,-10.5) [shape=circle, fill=black] {};
\node (b) at (-8,-10.5) [shape=circle, fill=black] {};
\node (c) at (-12,-6) [shape=circle, fill=black] {};
\node (d) at (-12,-19.5) [shape=circle, fill=black] {};
\node (e) at (-8,-19.5) [shape=circle, fill=black] {};
\node (f) at (-12,-15) [shape=circle, fill=black] {};

\draw [thick] (a) to (c);
\draw [thick] (d) to (f);
\draw [thick] (b) to (-8,-15);
\draw [thick] (e) to (-8,-24);


\draw [thick] (15,-15) arc (0:360:5);
\node (a') at (12,-10.5) [shape=circle, fill=black] {};
\node (b') at (8,-10.5) [shape=circle, fill=black] {};
\node (c') at (12,-6) [shape=circle, fill=black] {};
\node (d') at (12,-19.5) [shape=circle, fill=black] {};
\node (e') at (8,-19.5) [shape=circle, fill=black] {};
\node (f') at (12,-15) [shape=circle, fill=black] {};

\draw [thick] (a') to (c');
\draw [thick] (d') to (f');
\draw [thick] (b') to (8,-15);
\draw [thick] (e') to (8,-24);

\end{tikzpicture}

\end{center}
\caption{The third double coverings of $\mathcal S$ and $\overline\mathcal S$.}
\label{coversofSSbar}
\end{figure}
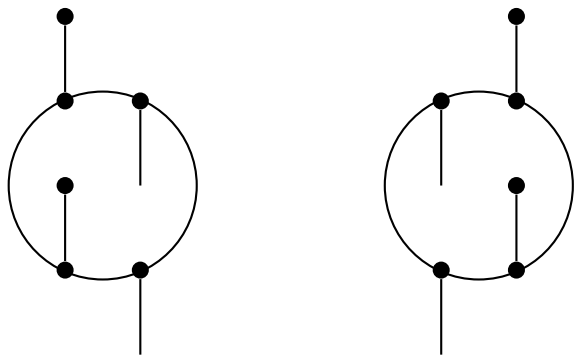


\begin{figure}[h!]
\begin{center}
\begin{tikzpicture}[scale=0.35, inner sep=0.8mm]

\draw [thick] (2,8) arc (0:360:2);
\node (A) at (0.75,8) [shape=circle, fill=black] {};
\node (A') at (-0.75,8) [shape=circle, fill=black] {};
\node (B) at (0,10) [shape=circle, fill=black] {};
\node (C) at (0.75,6.15) [shape=circle, fill=black] {};
\node (C') at (-0.75,6.15) [shape=circle, fill=black] {};
\node (D) at (0,12) [shape=circle, fill=black] {};

\draw [thick] (C) to (A);
\draw [thick] (C') to (A');
\draw [thick] (B) to (D);

\draw [thick] (2,-2) arc (0:360:2);
\node (a) at (0.75,-2) [shape=circle, fill=black] {};
\node (a') at (-0.75,-2) [shape=circle, fill=black] {};
\node (b) at (0,0) [shape=circle, fill=black] {};
\node (c) at (0.75,-3.85) [shape=circle, fill=black] {};
\node (c') at (-0.75,-3.85) [shape=circle, fill=black] {};
\node (d) at (0,2) [shape=circle, fill=black] {};

\draw [thick] (c) to (a);
\draw [thick] (c') to (a');
\draw [thick] (b) to (d);
\draw [thick, rounded corners] (D) to (-3,12) to (-3,2) to (3,2) to (3,12) to (D);

\node at (0,-6) {$\mathcal A(1)\mathcal A$};


\draw [thick] (14,10) arc (0:360:3);
\node (A1) at (8,10) [shape=circle, fill=black] {};
\node (A1') at (14,10) [shape=circle, fill=black] {};
\node (B1) at (11,7) [shape=circle, fill=black] {};
\node (C1) at (11,13) [shape=circle, fill=black] {};
\node (D1) at (11,11) [shape=circle, fill=black] {};

\draw [thick] (B1) to (11,5.5);
\draw [thick] (C1) to (D1);
\draw [thick] (12,10) arc (0:360:1);

\draw [thick] (14,0) arc (0:360:3);
\node (a1) at (8,0) [shape=circle, fill=black] {};
\node (a1') at (14,0) [shape=circle, fill=black] {};
\node (b1) at (11,-3) [shape=circle, fill=black] {};
\node (c1) at (11,3) [shape=circle, fill=black] {};
\node (d1) at (11,1) [shape=circle, fill=black] {};

\draw [thick] (b1) to (11,-4.5);
\draw [thick] (c1) to (d1);
\draw [thick] (12,0) arc (0:360:1);

\draw [thick, rounded corners] (A1) to (6,10) to (6,0) to (a1);
\draw [thick, rounded corners] (A1') to (16,10) to (16,0) to (a1');

\node at (11,-6) {$\mathcal B(3)\mathcal B$};


\draw [thick] (24,8) arc (0:360:2);
\node (A2) at (24,8) [shape=circle, fill=black] {};
\node (A2') at (20,8) [shape=circle, fill=black] {};
\node (B2) at (22,10) [shape=circle, fill=black] {};
\node (C2) at (23,9.7) [shape=circle, fill=black] {};
\node (C2') at (21,9.7) [shape=circle, fill=black] {};
 
\draw [thick] (C2) to (23,8);
\draw [thick] (C2') to (21,8);
\draw [thick] (B2) to (22,11.5);

\node (D2) at (22,13) [shape=circle, fill=black] {};
\draw [thick, rounded corners] (A2) to (25,8) to (25,13) to (19,13) to (19,8) to (A2');

\node (E2) at (22,15) [shape=circle, fill=black] {};
\draw [thick] (D2) to (E2);

\draw [thick] (24,-2) arc (0:360:2);
\node (a2) at (24,-2) [shape=circle, fill=black] {};
\node (a2') at (20,-2) [shape=circle, fill=black] {};
\node (b2) at (22,0) [shape=circle, fill=black] {};
\node (c2) at (23,-0.3) [shape=circle, fill=black] {};
\node (c2') at (21,-0.3) [shape=circle, fill=black] {};
 
\draw [thick] (c2) to (23,-2);
\draw [thick] (c2') to (21,-2);
\draw [thick] (b2) to (22,1.5);

\node (d2) at (22,3) [shape=circle, fill=black] {};
\draw [thick, rounded corners] (a2) to (25,-2) to (25,3) to (19,3) to (19,-2) to (a2');

\node (e2) at (22,5) [shape=circle, fill=black] {};
\draw [thick] (d2) to (e2);

\draw [thick, rounded corners] (E2) to (18,15) to (18,5) to (26,5) to (26,15) to (E2);

\node at (22,-6)  {$\mathcal C(1)\mathcal C$};

\end{tikzpicture}

\end{center}
\caption{The dessins $\mathcal A(1)\mathcal A$, $\mathcal B(3)\mathcal B$ and $\mathcal C(1)\mathcal C$.}
\label{dessinsA*B*C*}
\end{figure}
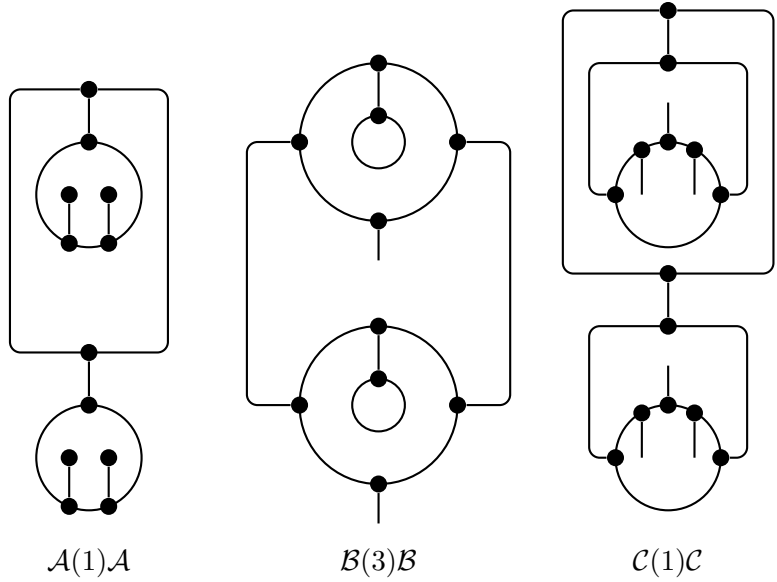


\medskip

\noindent{\bf Example} Figure~\ref{dessinsA*B*C*} shows the dessins $\mathcal A^*=\mathcal A(1)\mathcal A$, $\mathcal B^*=\mathcal B(3)\mathcal B$ and $\mathcal C^*=\mathcal C(1)\mathcal C$, of degrees $28$, $30$ and $42$. Since $\mathcal A$ is planar and has $\beta=2$, the above argument shows that $\mathcal A^*$ has the same monodromy group ${\rm PSL}_2(13)$ as $\mathcal A$, but now in its imprimitive representation of degree 28 on the cosets of the unique subgroup $H^*\cong {\rm C}_{13}\rtimes{\rm C}_3$ of index $2$ in $H\cong{\rm C}_{13}\rtimes{\rm C}_6$. On the other hand, since the monodromy group of $\mathcal B$ is ${\rm A}_{15}$, and a point stabiliser $H\cong{\rm A}_{14}$ in ${\rm A}_{15}$ has no subgroup of index $2$, the monodromy group of $\mathcal B^*$ must be a proper covering group of ${\rm A}_{15}$; according to GAP it is an extension of an elementary abelian group of order $2^{14}$ by ${\rm A}_{15}$, which must be the imprimitive group $({\rm S}_2\wr{
\rm A}_{15})\cap{\rm A}_{30}$. Similarly, although the monodromy group of $\mathcal C$ is ${\rm PGL}_3(2)\cong{\rm PSL}_2(7)$, that of $\mathcal C^*$ is a covering of this group by an elementary abelian group of order $2^6$; as an abstract group, and as the automorphism group of the regular cover of $\mathcal C^*$, this is the unique Hurwitz group of genus $257$, corresponding to the normal subgroup $K'K^2$ of $\Delta$, where $K$ corresponds to the regular Hurwitz dessin of genus $3$. These dessins $\mathcal A^*$, $\mathcal B^*$ and $\mathcal C^*$ are all reflexible, because $\mathcal A$, $\mathcal B$ and $\mathcal C$ are reflexible, invariant under a reflection which preserves the chosen handle (as an unordered pair of fixed points).

\begin{figure}[h!]
\begin{center}
\begin{tikzpicture}[scale=0.3, inner sep=0.8mm]

\draw [thick] (14,30) arc (0:360:3);
\node (A1) at (8,30) [shape=circle, fill=black] {};
\node (A1') at (14,30) [shape=circle, fill=black] {};
\node (B1) at (9.4,27.5) [shape=circle, fill=black] {};
\node (B1') at (12.6,27.5) [shape=circle, fill=black] {};
\node (C1) at (10.45,29.15) [shape=circle, fill=black] {};
\draw [thick] (B1) to (C1);
\draw [thick] (B1') to (13.6,26);
\draw [thick] (12,30) arc (0:360:1);

\draw [thick] (14,20) arc (0:360:3);
\node (D1) at (8,20) [shape=circle, fill=black] {};
\node (D1') at (14,20) [shape=circle, fill=black] {};
\node (E1) at (9.4,17.5) [shape=circle, fill=black] {};
\node (E1') at (12.6,17.5) [shape=circle, fill=black] {};
\node (F1) at (10.45,19.15) [shape=circle, fill=black] {};
\draw [thick] (E1) to (F1);
\draw [thick] (E1') to (13.6,16);
\draw [thick] (12,20) arc (0:360:1);

\draw [thick, rounded corners] (A1) to (6,30) to (6,20) to (D1);
\draw [thick, rounded corners] (A1') to (16,30) to (16,20) to (D1');


\draw [thick] (29,10) arc (0:360:3);
\node (A2) at (23,10) [shape=circle, fill=black] {};
\node (A2') at (29,10) [shape=circle, fill=black] {};
\node (B2) at (24.4,7.5) [shape=circle, fill=black] {};
\node (B2') at (27.6,7.5) [shape=circle, fill=black] {};
\node (C2) at (26.55,9.15) [shape=circle, fill=black] {};
\draw [thick] (B2') to (C2);
\draw [thick] (B2) to (23.4,6);
\draw [thick] (27,10) arc (0:360:1);

\draw [thick] (29,0) arc (0:360:3);
\node (D2) at (23,0) [shape=circle, fill=black] {};
\node (D2') at (29,0) [shape=circle, fill=black] {};
\node (E2) at (24.4,-2.5) [shape=circle, fill=black] {};
\node (E2') at (27.6,-2.5) [shape=circle, fill=black] {};
\node (F2) at (25.45,-0.85) [shape=circle, fill=black] {};
\draw [thick] (E2) to (F2);
\draw [thick] (E2') to (28.6,-4);
\draw [thick] (27,0) arc (0:360:1);

\draw [thick, rounded corners] (A2) to (21,10) to (21,0) to (D2);
\draw [thick, rounded corners] (A2') to (31,10) to (31,0) to (D2');


\draw [thick] (14,10) arc (0:360:3);
\node (a1) at (8,10) [shape=circle, fill=black] {};
\node (a1') at (14,10) [shape=circle, fill=black] {};
\node (b1) at (9.4,7.5) [shape=circle, fill=black] {};
\node (b1') at (12.6,7.5) [shape=circle, fill=black] {};
\node (c1) at (10.45,9.15) [shape=circle, fill=black] {};
\draw [thick] (b1) to (c1);
\draw [thick] (b1') to (13.6,6);
\draw [thick] (12,10) arc (0:360:1);

\draw [thick] (14,0) arc (0:360:3);
\node (d1) at (8,0) [shape=circle, fill=black] {};
\node (d1') at (14,0) [shape=circle, fill=black] {};
\node (e1) at (9.4,-2.5) [shape=circle, fill=black] {};
\node (e1') at (12.6,-2.5) [shape=circle, fill=black] {};
\node (f1) at (11.55,-0.85) [shape=circle, fill=black] {};
\draw [thick] (e1') to (f1);
\draw [thick] (e1) to (8.4,-4);
\draw [thick] (12,0) arc (0:360:1);

\draw [thick, rounded corners] (a1) to (6,10) to (6,0) to (d1);
\draw [thick, rounded corners] (a1') to (16,10) to (16,0) to (d1');


\draw [thick] (29,30) arc (0:360:3);
\node (a2) at (23,30) [shape=circle, fill=black] {};
\node (a2') at (29,30) [shape=circle, fill=black] {};
\node (b2) at (24.4,27.5) [shape=circle, fill=black] {};
\node (b2') at (27.6,27.5) [shape=circle, fill=black] {};
\node (c2) at (26.55,29.15) [shape=circle, fill=black] {};
\draw [thick] (b2') to (c2);
\draw [thick] (b2) to (23.4,26);
\draw [thick] (27,30) arc (0:360:1);

\draw [thick] (29,20) arc (0:360:3);
\node (d2) at (23,20) [shape=circle, fill=black] {};
\node (d2') at (29,20) [shape=circle, fill=black] {};
\node (e2) at (24.4,17.5) [shape=circle, fill=black] {};
\node (e2') at (27.6,17.5) [shape=circle, fill=black] {};
\node (f2) at (26.55,19.15) [shape=circle, fill=black] {};
\draw [thick] (e2') to (f2);
\draw [thick] (e2) to (23.4,16);
\draw [thick] (27,20) arc (0:360:1);

\draw [thick, rounded corners] (a2) to (21,30) to (21,20) to (d2);
\draw [thick, rounded corners] (a2') to (31,30) to (31,20) to (d2');

\end{tikzpicture}

\end{center}
\caption{The dessins $\mathcal B(2)\mathcal B$.}
\label{dessinsB2B}
\end{figure}
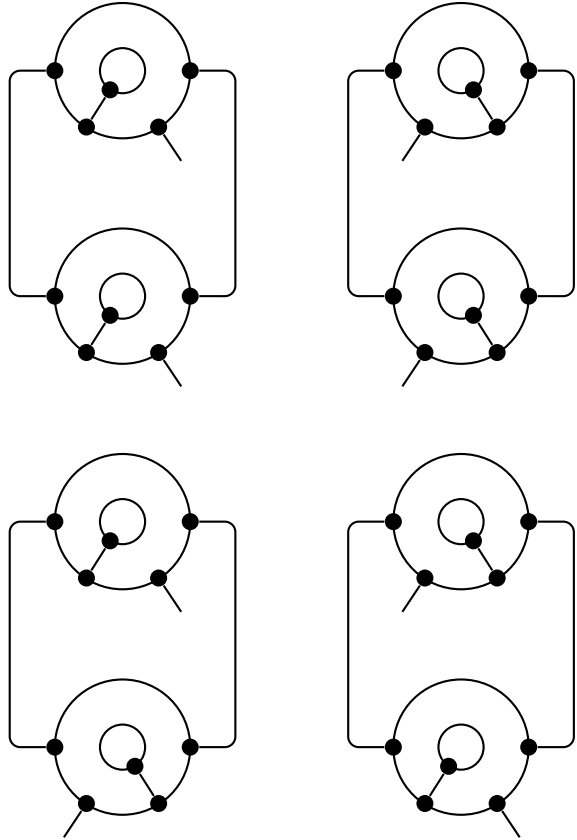

\medskip

\noindent{\bf Example} The dessin $\mathcal B$ also has two $(2)$-handles, transposed by its reflection. We can therefore form four dessins $\mathcal B(2)\mathcal B$ of degree $30$, shown in Figure~\ref{dessinsB2B}, by choosing one of these handles  in each copy of $\mathcal B$. The two dessins in the top row, where the same handle has been chosen from each copy, clearly form a chiral pair; GAP shows that their monodromy group is $({\rm S}_2\wr{\rm A}_{15})\cap{\rm A}_{30}$. The two dessins in the bottom row, where different handles have been chosen from each copy, are also mirror images of each other, but they are in fact isomorphic, by a half-turn of the sphere transposing the two copies of $\mathcal B$, so that this dessin is invariant under the antipodal symmetry of the sphere. In this case GAP shows that the monodromy group is ${\rm A}_{30}$.

\medskip

As a variation on the concept of a double covering $\mathcal D(k)\mathcal D$, if $\mathcal D$ has a $(k)$-handle then so has $\overline\mathcal D$, and even if $\mathcal D\not\cong \overline\mathcal D$ we can form the dessin $\mathcal D(k)\overline\mathcal D\cong\overline\mathcal D(k)\mathcal D$.
If $\mathcal D$ has degree $n$ and signature $(g;3^{[\alpha]},2^{[\beta]},7^{[\gamma]})$, then so has $\overline\mathcal D$, so as in the case of $\mathcal D^*$ it follows from Theorem~\ref{additivity} that $\mathcal D(k)\overline\mathcal D$ has degree $n^*=2n$ and signature $(2g;3^{[2\alpha]},2^{[2\beta-4]},7^{[2\gamma]})$.

\medskip

\noindent{\bf Example} If we join the dessins $\mathcal S$ and $\overline\mathcal S$ in Figure~{\ref{dessinsSSbar}, we see a phenomenon similar to that in the preceding example, as shown in Figure~\ref{dessinSkSbar}. The dessins $\mathcal S(1)\overline\mathcal S$ and $\overline\mathcal S(1)\mathcal S$ in the first row, clearly mirror images of each other, are also isomorphic under a rotation of order $2$, and are each invariant under the antipodal isometry of the sphere. Although the monodromy group of $\mathcal S$ and $\overline\mathcal S$ is ${\rm PGL}_3(2)\cong {\rm PSL}_2(7)$, and for $\mathcal S(1)\mathcal S$ and $\overline\mathcal S(1)\overline\mathcal S$ it is a covering of this group, this dessin has monodromy group ${\rm PSL}_2(13)$, in its natural representation. It is, in fact, one of a Galois orbit of three dessins corresponding to this action; the others are $\mathcal A$, shown in Figure~\ref{dessinsABC}, and the dessin $\mathcal S(2)\overline\mathcal S\cong\overline\mathcal S(2)\mathcal S$, also invariant under the antipodal isometry, shown in the second row of Figure~\ref{dessinSkSbar}.

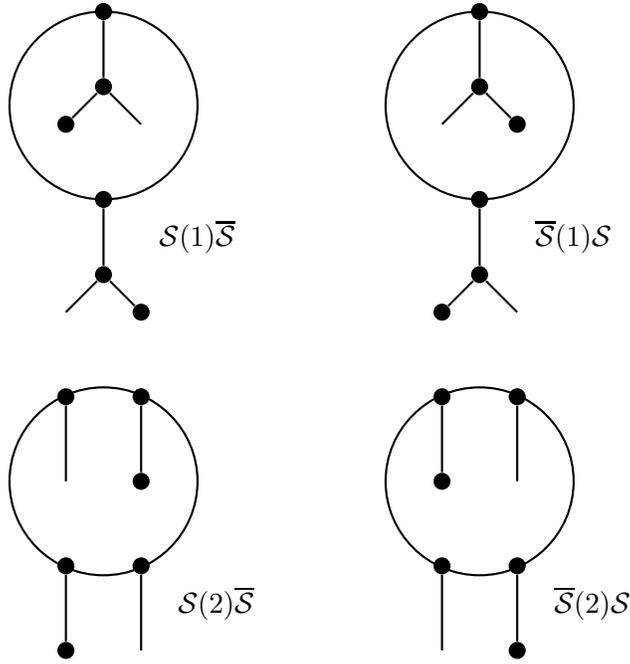
\begin{figure}[h!]
\begin{center}
\begin{tikzpicture}[scale=0.25, inner sep=0.8mm]

\draw [thick] (-5,5) arc (0:360:5);
\node (A) at (-10,10) [shape=circle, fill=black] {};
\node (B) at (-10,6) [shape=circle, fill=black] {};
\node (C) at (-10,0) [shape=circle, fill=black] {};
\node (D) at (-10,-4) [shape=circle, fill=black] {};
\node (E) at (-12,4) [shape=circle, fill=black] {};
\node (F) at (-8,-6) [shape=circle, fill=black] {};

\draw [thick] (A) to (B);
\draw [thick] (C) to (D);
\draw [thick] (E) to (B) to (-8,4);
\draw [thick] (F) to (D) to (-12,-6);

\node at (-5,-2) {$\mathcal S(1)\overline\mathcal S$};


\draw [thick] (15,5) arc (0:360:5);
\node (A') at (10,10) [shape=circle, fill=black] {};
\node (B') at (10,6) [shape=circle, fill=black] {};
\node (C') at (10,0) [shape=circle, fill=black] {};
\node (D') at (10,-4) [shape=circle, fill=black] {};
\node (E') at (12,4) [shape=circle, fill=black] {};
\node (F') at (8,-6) [shape=circle, fill=black] {};

\draw [thick] (A') to (B');
\draw [thick] (C') to (D');
\draw [thick] (E') to (B') to (8,4);
\draw [thick] (F') to (D') to (12,-6);

\node at (15,-2) {$\overline\mathcal S(1)\mathcal S$};


\draw [thick] (-5,-15) arc (0:360:5);
\node (a) at (-12,-10.5) [shape=circle, fill=black] {};
\node (b) at (-8,-10.5) [shape=circle, fill=black] {};
\node (c) at (-8,-15) [shape=circle, fill=black] {};
\node (d) at (-12,-19.5) [shape=circle, fill=black] {};
\node (e) at (-8,-19.5) [shape=circle, fill=black] {};
\node (f) at (-12,-24) [shape=circle, fill=black] {};

\draw [thick] (b) to (c);
\draw [thick] (d) to (f);
\draw [thick] (a) to (-12,-15);
\draw [thick] (e) to (-8,-24);

\node at (-4,-21.5) {$\mathcal S(2)\overline\mathcal S$};


\draw [thick] (15,-15) arc (0:360:5);
\node (a') at (12,-10.5) [shape=circle, fill=black] {};
\node (b') at (8,-10.5) [shape=circle, fill=black] {};
\node (c') at (8,-15) [shape=circle, fill=black] {};
\node (d') at (12,-19.5) [shape=circle, fill=black] {};
\node (e') at (8,-19.5) [shape=circle, fill=black] {};
\node (f') at (12,-24) [shape=circle, fill=black] {};

\draw [thick] (b') to (c');
\draw [thick] (d') to (f');
\draw [thick] (a') to (12,-15);
\draw [thick] (e') to (8,-24);

\node at (16,-21.5) {$\overline\mathcal S(2)\mathcal S$};

\end{tikzpicture}

\end{center}
\caption{The dessins $\mathcal S(k)\overline\mathcal S$ and $\overline\mathcal S(k)\mathcal S$ for $k=1,2$. }
\label{dessinSkSbar}
\end{figure}


\subsection{Multiple joins}

If dessins $\mathcal D_i\;(i=1,2)$ each have several mutually disjoint handles, then one can use these to make a multiple join $\mathcal D$ by joining $m_k$ $(k)$-handles in $\mathcal D_1$ to the same number in $\mathcal D_2$ (provided that many exist) for $k=1, 2, 3$. If each $\mathcal D_i$ has degree $n_i$ and signature $(g_i; 2^{[\alpha_i]}, 3^{[\beta_i]}, 7^{[\gamma_i]})$, then $\mathcal D$  has degree $n_1+n_2$ and signature $(g; 2^{[\alpha]}, 3^{[\beta]}, 7^{[\gamma]})$ where $g=g_1+g_2+m-1$ ($m:=\sum_km_k$), $\alpha=\alpha_1+\alpha_2$, $\beta=\beta_1+\beta_2-4m$ and $\gamma=\gamma_1+\gamma_2$.

In this way, Hurwitz dessins of arbitrary genus $g$ can be formed from a pair of Hurwitz dessins of genus $0$ with at least $g+1$ compatible handles. For instance, since the dessin $\mathcal G$ in Figure~\ref{dessinG} (based on Conder's diagram $G$ in~\cite{Conder-80}) has genus $0$ and three disjoint $(1)$-handles, we could take each $\mathcal D_i$ to be the composition $\mathcal G(1)\cdots(1)\mathcal G$ of $g-1$ copies of $\mathcal G$. (Using such a composition as a `stem', on which to attach further Hurwitz dessins, is a basic idea in~\cite{Conder-80} and several subsequent papers.)


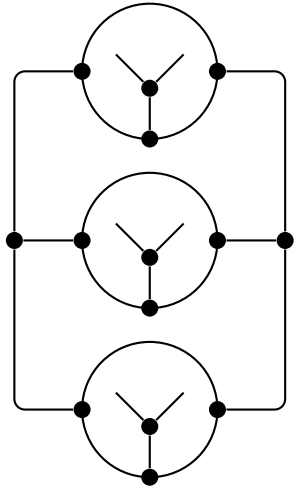
\begin{figure}[h!]
\begin{center}
\begin{tikzpicture}[scale=0.45, inner sep=0.8mm]

\draw [thick] (2,5) arc (0:360:2);
\node (a1) at (2,5) [shape=circle, fill=black] {};
\node (a1') at (-2,5) [shape=circle, fill=black] {};
\node (b1) at (0,3) [shape=circle, fill=black] {};
\node (c1) at (0,4.5) [shape=circle, fill=black] {};

\draw [thick] (b1) to (c1);
\draw [thick] (1,5.5) to (c1) to (-1,5.5);

\draw [thick] (2,0) arc (0:360:2);
\node (a2) at (2,0) [shape=circle, fill=black] {};
\node (a2') at (-2,0) [shape=circle, fill=black] {};
\node (b2) at (0,-2) [shape=circle, fill=black] {};
\node (c2) at (0,-0.5) [shape=circle, fill=black] {};

\draw [thick] (b2) to (c2);
\draw [thick] (1,0.5) to (c2) to (-1,0.5);

\draw [thick] (2,-5) arc (0:360:2);
\node (a3) at (2,-5) [shape=circle, fill=black] {};
\node (a3') at (-2,-5) [shape=circle, fill=black] {};
\node (b3) at (0,-7) [shape=circle, fill=black] {};
\node (c3) at (0,-5.5) [shape=circle, fill=black] {};

\draw [thick] (b3) to (c3);
\draw [thick] (1,-4.5) to (c3) to (-1,-4.5);

\node (d) at (4,0) [shape=circle, fill=black] {};
\node (d') at (-4,0) [shape=circle, fill=black] {};
\draw [thick, rounded corners] (a1) to (4,5) to (d) to (4,-5) to (a3);
\draw [thick, rounded corners] (a1') to (-4,5) to (d') to (-4,-5) to (a3');
\draw [thick] (a2) to (d);
\draw [thick] (a2') to (d');

\end{tikzpicture}

\end{center}
\caption{The dessin $\mathcal G$.}
\label{dessinG}
\end{figure}


\subsection{$y$-handles from ${\rm PSL}_2(q)$}

Motivated by the example of the dessin $\mathcal A$, it is natural to look for other instances of $y$-handles arising from the Hurwitz groups $G={\rm PSL}_2(q)$ in their natural representation. Based on Macbeath's classification~\cite{Macb} of the Hurwitz groups of this form, as explained in \S\ref{Macbeath}, we have the following:

\begin{theorem}\label{PSL2qnathandles} A Hurwitz dessin corresponding to the natural representation of the group $G={\rm PSL}_2(q)$, $q=p^e$, has a $(k)$-handle if and only if either
\begin{itemize}
\item $q=p=13$, with ${\rm tr}(z)=\pm 5$ and $k=1$, or
\item $q=p=29$, with ${\rm tr}(z)=\pm 3$ and $k=2$, or
\item $q=p=41$, with ${\rm tr}(z)=\pm 11$ and $k=3$.
\end{itemize}

\end{theorem}

\noindent{\sl Proof.} If a $(k)$-handle exists then $y$ has two fixed points, so $q$ must be odd and $y$ must have order dividing $(q-1)/2$ where $q\equiv 1$  mod~$(4)$. By applying an automorphism of $G$ we may assume that $y$ has the form $t\mapsto -t$, fixing $a=0$ and $b=\infty$, or equivalently, that $y$ corresponds to the pair of elements
\[\pm\left( \begin{array}{cc}
i & 0 \\ 
0 & -i \\
\end{array}\right)\]
of ${\rm SL}_2(q)$, where $i^2=-1$. An element of $G$ has order $3$ if and only if it has trace $\pm 1$, so the most general form for $x$ is
\[\pm\left( \begin{array}{cc}
a & b \\ 
c & 1-a \\
\end{array}\right)\]
where $a(1-a)-bc=1$. It then follows that $xy$ has trace
\[t=\pm (2a-1)i.\]

Now any element of order $7$, such as $xy$, has trace $\pm\tau_j$ where $\tau_j=\lambda^j+\lambda^{-j}$ for $j=1, 2, 3$ and $\lambda$ is a primitive 7th root of $1$ (possibly in an extension field). One easily checks that $\tau_1+\tau_2+\tau_3=-1$, $\tau_1\tau_2+\tau_2\tau_3+\tau_3\tau_1=-2$ and $\tau_1\tau_2\tau_3=1$, so the elements $\pm\tau_j$ are respectively the roots of the polynomials
\[p_{\pm}(t):=t(t^2-2)\pm(t^2-1).\]

First we look for $(1)$-handles. These require $x:0\mapsto\infty$, so $a=0$ and $xy$ has trace $t=\pm i$.
Thus $i$ must be a root of one of the polynomials $p_{\pm}(t)$, so $3i=\pm 2$ and hence $p=13$ with $t=\pm 5$. By Macbeath's classification, it follows that $q=13$. The dessin is that shown earlier as $\mathcal A$ in Figure~\ref{dessinsABC}.

In the case of a $(2)$-handle we require $xyx:0\mapsto\infty$. Writing $x$ and $y$ as above we find that
\[xyx: 0\mapsto\frac{(1-2a)b}{2a^2-3a+2},\]
so this is equivalent to
\[2a^2-3a+2=0.\]
Putting $t=\delta(2a-1)i$ with $\delta^2=1$ gives $t^2=-4a^2+4a-1=-2a+3$,
so
\[\begin{array}{cl}
0=p_{\pm}(t)&=\delta(2a-1)i(-2a+1)\pm(-2a+2)\\
&=\delta i(-4a^2+4a-1)\pm2(1-a)\\
&=\delta i(-2a+3)\pm2(1-a),\\
\end{array}
\]
and hence $0=(2a-3)^2+4(a-1)^2=8a^2-20a+13=-8a+5$, that is, $a=5/8$. This gives $0=50-120+128=58$ and hence $q=p=29$. The roots of $2a^2-3a+2$ in $\mathbb F_{29}$ are $a=-3$ and $-10$, and we have $i=\pm12$, so $xy$ has trace $t=\delta(2a-1)i=\pm 3$ or $\pm 9$. Now the roots of $p_+(t)$ are $3, 12$ and 13, so $a=-3$, giving $t=\pm 3$.

It is now straightforward to check that the resulting dessin is the planar dessin $\mathcal F$ in Figure~\ref{dessinF}, corresponding to Conder's diagram $F$ in~\cite{Conder-80}.

\begin{figure}[h!]
\begin{center}
\begin{tikzpicture}[scale=0.5, inner sep=0.8mm]

\draw [thick] (2,0) arc (0:360:2);
\node (a) at (2,0) [shape=circle, fill=black] {};
\node (a') at (-2,0) [shape=circle, fill=black] {};
\node (b) at (0,2) [shape=circle, fill=black] {};
\node (c) at (0.75,-1.85) [shape=circle, fill=black] {};
\node (c') at (-0.75,-1.85) [shape=circle, fill=black] {};

\draw [thick] (c) to (0.75,0);
\draw [thick] (c') to (-0.75,0);

\node (d) at (0,4) [shape=circle, fill=black] {};
\draw [thick] (b) to (d);
\node (e) at (4,4) [shape=circle, fill=black] {};
\node (e') at (-4,4) [shape=circle, fill=black] {};
\node (f) at (4,2) [shape=circle, fill=black] {};
\node (f') at (-4,2) [shape=circle, fill=black] {};

\draw [thick] (4.5,1.5) arc (0:360:0.5);
\draw [thick] (-3.5,1.5) arc (0:360:0.5);

\draw [thick] [rounded corners] (a) to (6,0) to (6,4) to (-6,4) to (-6,0) to (a');
\draw [thick] (e) to (f);
\draw [thick] (e') to (f');

\end{tikzpicture}

\end{center}
\caption{The dessin $\mathcal F$.}
\label{dessinF}
\end{figure}
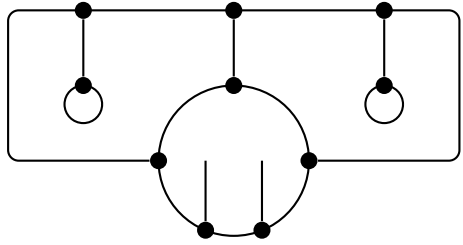


In the case of a $(3)$-handle we require $(xy)^2x:0\mapsto\infty$. Writing $x$ and $y$ as above we find that this is equivalent to
\[4a^3-8a^2+8a-3=0.\]
Putting $t=\delta(2a-1)i$ with $\delta^2=1$ gives $t^2=-4a^2+4a-1$,
so
\[\begin{array}{cl}
0=p_{\pm}(t)&=t(t^2-2)\pm(t^2-1)\\
&=\delta(2a-1)i(-4a^2+4a-3)\pm(-4a^2+4a-2)\\
&=\delta i(-8a^3+12a^2-10a+3)\pm(-4a^2+4a-2)\\
&=\delta i(-4a^2+6a-3)\pm(-4a^2+4a-2)\\
\end{array}
\]
so that
\[\begin{array}{cl}
0&=(-4a^2+6a-3)^2+(-4a^2+4a-2)^2\\
&=32a^4-80a^3+92a^2-52a+13.\\
\end{array}
\]
Using the cubic to eliminate the leading term gives
\[-16a^3+28a^2-28a+13=0=0\]
Repeating this gives
\[-4a^2+4a+1=0\]
Thus $4a^2=4a+1$, so $4a^3=4a^2+a=5a+1$, and hence the cubic becomes $5a-4=0$, that is, $a=4/5$. Substituting this in the above quadratic gives $64=105$, so $q=p=41$, with $a=9$. Here $i=\pm 9$ so $t=\pm 11$. As a check, $p_+(-11)=0$, so $-11=\lambda+\lambda^{-1}$ for some primitive 7th root of 1 in $\mathbb F_{41^3}$.

The corresponding dessin $\mathcal T$ has genus $g=1$. It is shown in Figure~\ref{L2(41)} with opposite sides of the outer hexagon identified to form a torus. \hfill$\square$

\begin{figure}[h!]
\begin{center}
\begin{tikzpicture}[scale=0.25, inner sep=0.8mm]

\node (a) at (-10,5) [shape=circle, fill=black] {};
\node (a') at (10,5) [shape=circle, fill=black] {};
\node (b) at (-5,5) [shape=circle, fill=black] {};
\node (b') at (5,5) [shape=circle, fill=black] {};
\node (c) at (0,0) [shape=circle, fill=black] {};
\node (d) at (0,-5) [shape=circle, fill=black] {};
\node (e) at (-15,0) [shape=circle, fill=black] {};
\node (e') at (15,0) [shape=circle, fill=black] {};
\node (f) at (-10,0) [shape=circle, fill=black] {};
\node (f') at (10,0) [shape=circle, fill=black] {};
\node (g) at (-10,-5) [shape=circle, fill=black] {};
\node (g') at (10,-5) [shape=circle, fill=black] {};
\node (h) at (-2,0.4) [shape=circle, fill=black] {};
\node (h') at (2,0.4) [shape=circle, fill=black] {};

\draw [thick] (5,5) arc (0:360:5);
\draw [thick] (a) to (b);
\draw [thick] (a') to (b');
\draw [thick] (c) to (d);
\draw [thick] (e) to (f);
\draw [thick] (e') to (f');
\draw [thick] (g) to (g');
\draw [thick] (-10,12) to (-10,-12);
\draw [thick] (10,12) to (10,-12);
\draw [thick] (h) to (-2,5);
\draw [thick] (h') to (2,5);
\draw [thick] (-16.5,6) to (e) to (-16.5,-6);
\draw [thick] (16.5,6) to (e') to (16.5,-6);

\draw [thick, dashed] (-13,12) to (13,12) to (20,0) to (13,-12) to (-13,-12) to (-20,0) to (-13,12);

\end{tikzpicture}

\end{center}
\caption{A torus dessin $\mathcal T$ with monodromy group ${\rm PSL}_2(41)$.}
\label{L2(41)}
\end{figure}
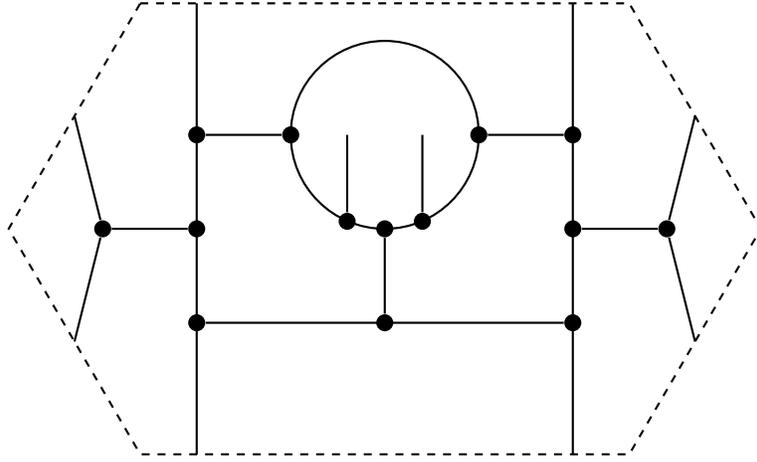


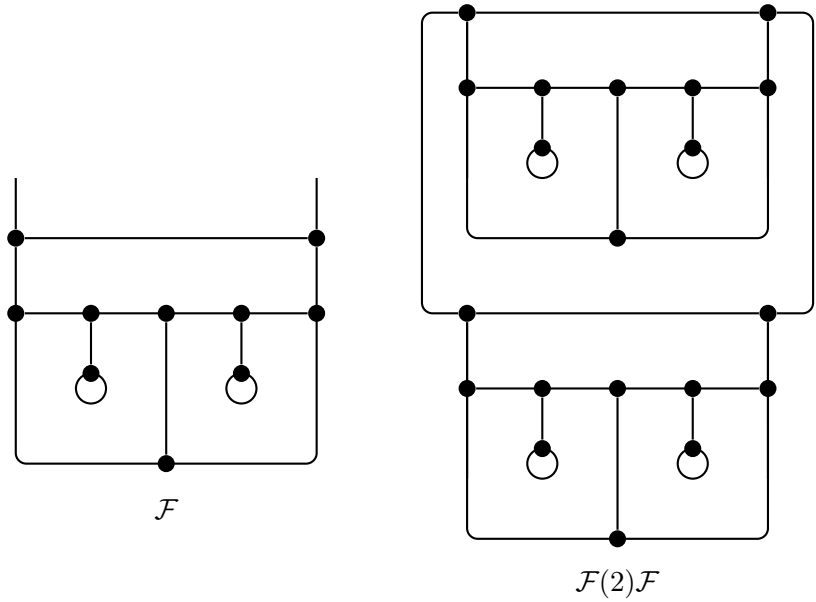
\begin{figure}[h!]
\begin{center}
\begin{tikzpicture}[scale=0.2, inner sep=0.8mm]

\node (a) at (-10,5) [shape=circle, fill=black] {};
\node (a') at (10,5) [shape=circle, fill=black] {};
\node (b) at (-10,0) [shape=circle, fill=black] {};
\node (b') at (10,0) [shape=circle, fill=black] {};
\node (c) at (-5,0) [shape=circle, fill=black] {};
\node (c') at (5,0) [shape=circle, fill=black] {};
\node (d) at (0,0) [shape=circle, fill=black] {};

\draw [thick] (a) to (a');
\draw [thick] (b) to (b');
\draw [thick] (a) to (-10, 9);
\draw [thick] (a') to (10,9);

\node (e) at (-5,-4) [shape=circle, fill=black] {};
\node (e') at (5,-4) [shape=circle, fill=black] {};
\node (f) at (0,-10) [shape=circle, fill=black] {};

\draw [thick, rounded corners] (a) to (-10,-10) to (10,-10) to (a');
\draw [thick] (c) to (e);
\draw [thick] (c') to (e');
\draw [thick] (d) to (f);

\draw [thick] (-4,-5) arc (0:360:1);
\draw [thick] (6,-5) arc (0:360:1);

\node at (0,-13) {$\mathcal F$};


\node (A) at (20,20) [shape=circle, fill=black] {};
\node (A') at (40,20) [shape=circle, fill=black] {};
\node (B) at (20,15) [shape=circle, fill=black] {};
\node (B') at (40,15) [shape=circle, fill=black] {};
\node (C) at (25,15) [shape=circle, fill=black] {};
\node (C') at (35,15) [shape=circle, fill=black] {};
\node (D) at (30,15) [shape=circle, fill=black] {};

\draw [thick] (A) to (A');
\draw [thick] (B) to (B');
\draw [thick] (A) to (20, 9);
\draw [thick] (A') to (40,9);

\node (E) at (25,11) [shape=circle, fill=black] {};
\node (E') at (35,11) [shape=circle, fill=black] {};
\node (F) at (30,5) [shape=circle, fill=black] {};

\draw [thick, rounded corners] (A) to (20,5) to (40,5) to (A');
\draw [thick] (C) to (E);
\draw [thick] (C') to (E');
\draw [thick] (D) to (F);

\draw [thick] (26,10) arc (0:360:1);
\draw [thick] (36,10) arc (0:360:1);


\node (A1) at (20,0) [shape=circle, fill=black] {};
\node (A1') at (40,0) [shape=circle, fill=black] {};
\node (B1) at (20,-5) [shape=circle, fill=black] {};
\node (B1') at (40,-5) [shape=circle, fill=black] {};
\node (C1) at (25,-5) [shape=circle, fill=black] {};
\node (C1') at (35,-5) [shape=circle, fill=black] {};
\node (D1) at (30,-5) [shape=circle, fill=black] {};

\draw [thick] (A1) to (A1');
\draw [thick] (B1) to (B1');
\draw [thick] (A1) to (20, -11);
\draw [thick] (A1') to (40,-11);

\node (E1) at (25,-9) [shape=circle, fill=black] {};
\node (E1') at (35,-9) [shape=circle, fill=black] {};
\node (F1) at (30,-15) [shape=circle, fill=black] {};

\draw [thick, rounded corners] (A1) to (20,-15) to (40,-15) to (A1');
\draw [thick] (C1) to (E1);
\draw [thick] (C1') to (E1');
\draw [thick] (D1) to (F1);

\draw [thick] (26,-10) arc (0:360:1);
\draw [thick] (36,-10) arc (0:360:1);

\draw [thick, rounded corners] (A) to (17,20) to (17,0) to (A1);
\draw [thick, rounded corners] (A') to (43,20) to (43,0) to (A1');

\node at (30,-18) {$\mathcal F(2)\mathcal F$};

\end{tikzpicture}

\end{center}
\caption{$\mathcal F$ and $\mathcal F(2)\mathcal F$.}
\label{F,F2F}
\end{figure}

\medskip

The dessin $\mathcal F$ with monodromy group ${\rm PSL}_2(29)$ which appears in this theorem corresponds to a conjugacy class of subgroups $M$ of index $30$ in $\Delta$ with signature $(0; 2,2,7,7)$. As in the case $q=13$, it follows from Lemma~\ref{uniqueD*} that the dessin $\mathcal F^*=\mathcal F(2)\mathcal F$, shown on the right in Figure~\ref{F,F2F}, also has monodromy group ${\rm PSL}_2(29)$, acting on the cosets of the unique subgroup $H^*=H^2\cong{\rm C}_{29}\rtimes{\rm C}_7$ of index $2$ in a natural point stabiliser $H\cong{\rm C}_{29}\rtimes{\rm C}_{14}$. 

However, Lemma~\ref{uniqueD*} does not apply to the dessin $\mathcal T$ with monodromy group ${\rm PSL}_2(41)$, since this has genus $1$. Indeed, the dessin $\mathcal T^*=\mathcal T(3)\mathcal T$, which has genus $2$ and has no fixed points for $y$, cannot be that corresponding to the unique subgroup $H^2\cong{\rm C}_{29}\rtimes{\rm C}_{10}$ of index $2$ in $H\cong{\rm C}_{29}\rtimes{\rm C}_{20}$, since the latter dessin is an unbranched covering of $\mathcal T$, of genus $1$, with four fixed points for $y$ forming two $(3)$-handles.

\medskip



Theorem~\ref{PSL2qnathandles} shows that although the number $\beta$ of fixed points of $y$ is invariant under the absolute Galois group ${\mathbb G}={\rm Gal}\,\overline{\mathbb Q}/{\mathbb Q}$ (as indeed are the cycle structures of $x, y$ and $z$ for any dessin, see~\cite{JS}), the property of having a $(k)$-handle is not invariant. As Macbeath showed in Theorem~\ref{Macbthm}, for each prime $p\equiv\pm 1$ mod~$(7)$ there are three Hurwitz dessins with the natural representation of ${\rm PSL}_2(p)$ as their monodromy group, distinguished by the values of $\pm{\rm tr}(z)$. In each case, as shown by Streit~\cite{Streit} they form an orbit of $\mathbb G$, but for $p=13$, $29$ and $41$ only one of them has a handle, by Theorem~\ref{PSL2qnathandles}. For the other two dessins in each orbit the two fixed points of $y$ lie in different faces, as shown in Figure~\ref{orbitA} for $p=13$.

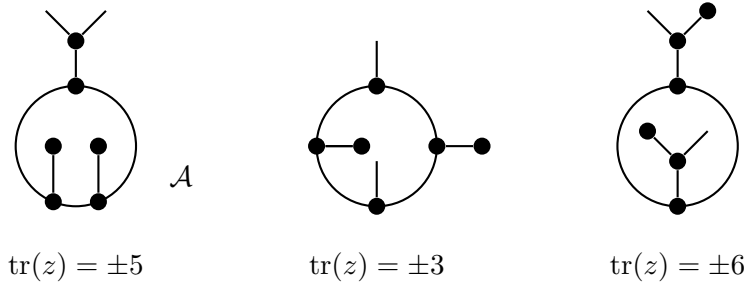
\begin{figure}[h!]
\begin{center}
\begin{tikzpicture}[scale=0.4, inner sep=0.8mm]

\draw [thick] (2,0) arc (0:360:2);
\node (a) at (0.75,0) [shape=circle, fill=black] {};
\node (a') at (-0.75,0) [shape=circle, fill=black] {};
\node (b) at (0,2) [shape=circle, fill=black] {};
\node (c) at (0.75,-1.85) [shape=circle, fill=black] {};
\node (c') at (-0.75,-1.85) [shape=circle, fill=black] {};
\node (d) at (0,3.5) [shape=circle, fill=black] {};

\draw [thick] (c) to (a);
\draw [thick] (c') to (a');
\draw [thick] (b) to (d);
\draw [thick] (-1,4.5) to (d) to (1,4.5);

\node at (3.5,-1) {$\mathcal A$};
\node at (0,-4) {${\rm tr}(z)=\pm 5$};


\draw [thick] (12,0) arc (0:360:2);
\node (A) at (12,0) [shape=circle, fill=black] {};
\node (B) at (10,2) [shape=circle, fill=black] {};
\node (C) at (8,0) [shape=circle, fill=black] {};
\node (D) at (10,-2) [shape=circle, fill=black] {};
\node (A') at (13.5,0) [shape=circle, fill=black] {};
\node (C') at (9.5,0) [shape=circle, fill=black] {};

\draw [thick] (A') to (A);
\draw [thick] (C') to (C);
\draw [thick] (B) to (10,3.5);
\draw [thick] (C') to (C);
\draw [thick] (D) to (10,-0.5);

\node at (10,-4) {${\rm tr}(z)=\pm 3$};


\draw [thick] (22,0) arc (0:360:2);
\node (E) at (20,2) [shape=circle, fill=black] {};
\node (F) at (20,-2) [shape=circle, fill=black] {};
\node (E') at (20,3.5) [shape=circle, fill=black] {};
\node (E'') at (21,4.5) [shape=circle, fill=black] {};
\node (F') at (20,-0.5) [shape=circle, fill=black] {};
\node (F'') at (19,0.5) [shape=circle, fill=black] {};
\draw [thick] (E') to (E);
\draw [thick] (F') to (F);
\draw [thick] (F'') to (F') to (21,0.5);
\draw [thick] (19,4.5) to (E') to (E'');

\node at (20,-4) {${\rm tr}(z)=\pm 6$};

\end{tikzpicture}

\end{center}
\caption{The Galois orbit containing $\mathcal A$.}
\label{orbitA}
\end{figure}

\begin{corollary}\label{PSLycorollary}
Suppose that a Hurwitz dessin $\mathcal D$ has a monodromy group $G={\rm PSL}_2(q)$ in a representation which covers the natural representation. Then $\mathcal D$ has a $(k)$-handle if and only if
\begin{itemize}
\item one of the three conclusions of Theorem~\ref{PSL2qnathandles} holds, and
\item the point stabilisers for $\mathcal D$ have even order.
\end{itemize}
\end{corollary}

\noindent{\sl Proof.} If $\mathcal D$ has a $(k)$-handle, then by Lemma~\ref{handlecover} so does its quotient dessin $\overline\mathcal D$ corresponding to the natural representation of $G$, so one of the three conclusions of Theorem~\ref{PSL2qnathandles} applies, with $q=p=13$, $29$ or $41$. The point stabilisers for $\mathcal D$ are subgroups $H_0$ of natural point stabilisers $H\cong {\rm C}_p\rtimes {\rm C}_{(p-1)/2}$ in $G$. If $H_0$ has odd order then $\mathcal D$ has no fixed points for $y$ and hence no handles, whereas if $H_0$ has odd index in $H$, then the handle in $\overline\mathcal D$ must lift to at least one handle in $\mathcal D$. This deals with the cases $q=13$ or $29$, since every subgroup of $H$ has odd order or odd index. In the case $q=41$ there are also subgroups $H_0$ of even order and even index; however, these are all subgroups of odd index in $H^2$, which correspoinds to a  dessin with two handles, so their corresponding dessins also have handles. \hfill$\square$

\medskip

In each of the three cases in Theorem~\ref{PSL2qnathandles} the natural point stabilisers $H$ are Frobenius groups $H\cong {\rm C}_p\rtimes{\rm C}_{(p-1)/2}$. The subgroups $H_0$ of $H$ therefore consist of one subgroup ${\rm C}_p\rtimes{\rm C}_f$ and a conjugacy class of $p$ subgroups ${\rm C}_f$ for each factor $f$ of $(p-1)/2$. By Corollary~\ref{PSLycorollary}, the subgroups corresponding to dessins with $y$-handles are those for which $f$ is even. When $p=13$, $29$ or $41$ these values of $f$ are $2$ and $6$, or $2$ and $14$, or $2$, $4$, $10$ and $20$. Of course, we need to exclude $H_0=H$ in order to obtain non-identity coverings.

\medskip.

\noindent{\bf Example.} When $q=p=13$ in Corollary~\ref{PSLycorollary} the point stabilisers $H_0<H$ giving $y$-handles consist of a normal subgroup ${\rm C}_{13}\rtimes {\rm C}_2\cong{\rm D}_{13}$, and two conjugacy classes of $13$ subgroups ${\rm C}_2$ and ${\rm C}_6$. The first case corresponds to the $3$-sheeted regular covering $\mathcal G$ of $\mathcal A$ of degree $42$ shown in Figure~\ref{dessinG}, branched over the two fixed points of $x$ and considered later in \S\ref{xexample} as an example of an $x$-join. In this case, and also when $H_0={\rm C}_2$, the numbers $\alpha$, $\beta$ and $\gamma$ of fixed points of $x$, $y$ and $z$ are $0$, $6$ and $0$, with three $(1)$-handles; the dessins have genus $g=0$ and $6$ respectively. When $H_0={\rm C}_6$ we have $\alpha=\beta=2$ and $\gamma=0$; the dessin has genus $2$, and there is one $(1)$-handle. The last two dessins, of degrees $546$ and $182$, are too large for us to draw.

When $q=p=29$ the subgroups $H_0<H$ providing handles consist of a normal subgroup ${\rm C}_{29}\rtimes {\rm C}_2\cong{\rm D}_{29}$, and two conjugacy classes of $29$ subgroups ${\rm C}_2$ and ${\rm C}_{14}$. The first corresponds to a $7$-sheeted regular covering of $\mathcal F$ of degree $210$, branched over the two fixed points of $z$. Here, and also when $H_0={\rm C}_2$, we have $\alpha=\gamma=0$ and $\beta=14$, with seven $(2)$-handles, but with $g=0$ and $70$ respectively. When $H_0={\rm C}_{14}$ we have $\alpha=0$ and $\beta=\gamma=2$, so there is one $(2)$-handle and the genus is $10$. 
 
When $q=p=41$ the relevant subgroups $H_0<H$ are three normal subgroups ${\rm C}_{41}\rtimes {\rm C}_f$ for $f=2$, $4$ or $10$, and four conjugacy classes of $41$ subgroups  ${\rm C}_f$ for $f=2, 4, 10$ or $20$. The corresponding dessins have $\alpha=\gamma=0$ and $\beta=40/f$, with $20/f$ $(3)$-handles. The normal subgroups correspond to $20/f$-sheeted unbranched regular coverings of the dessin in Figure~\ref{L2(41)}}, of genus $1$. The dessins corresponding to subgroups ${\rm C}_f$ have genus $1+400/f$.

\medskip

Of course, this leaves open the question of what handles can appear in other representations of Hurwitz groups ${\rm PSL}_2(q)$. As above, Lemma~\ref{handlecover} shows that it is sufficient to classify those handles arising from primitive representations, where the point-stabilisers are maximal subgroups. 


\section{$x$-joins}\label{xjoins}

\begin{figure}[h!]
\begin{center}
\begin{tikzpicture}[scale=0.5, inner sep=0.8mm]

\node (a) at (-6,5) [shape=circle, fill=black] {};
\node (b) at (-10,5) [shape=circle, fill=black] {};
\node (c) at (-8,3) [shape=circle, fill=black] {};
\draw [thick] (c) to (-8,1);
\draw [thick, dashed] (-8,1) to (-8,-1);
\draw [thick] (a) to (c) to (b);
\node at (-6,6) {$a$};
\node at (-10,6) {$b$};
\node at (-6,0) {$k=1$};


\node (b2) at (-2,5) [shape=circle, fill=black] {};
\node (a2) at (2,5) [shape=circle, fill=black] {};
\node (d2) at (-2,2) [shape=circle, fill=black] {};
\node (c2) at (2,2) [shape=circle, fill=black] {};
\draw [thick] (b2) to (d2);
\draw [thick] (a2) to (c2);
\draw [thick] (-3,2) to (3,2);
\draw [thick, dashed] (-4,2) to (-3,2);
\draw [thick, dashed] (4,2) to (3,2);
\node at (2,6) {$a$};
\node at (-2,6) {$b$};
\node at (0,0) {$k=2$};


\node (b2) at (8,5) [shape=circle, fill=black] {};
\node (a2) at (12,5) [shape=circle, fill=black] {};
\node (d2) at (8,2) [shape=circle, fill=black] {};
\node (c2) at (12,2) [shape=circle, fill=black] {};
\node (e) at (10,2) [shape=circle, fill=black] {};
\draw [thick] (a2) to (c2);
\draw [thick] (b2) to (d2);
\draw [thick] (7,2) to (13,2);
\draw [thick, dashed] (6,2) to (7,2);
\draw [thick, dashed] (14,2) to (13,2);
\draw [thick] (e) to (10,0);
\draw [thick, dashed] (10,0) to (10,-1);
\node at (12,6) {$a$};
\node at (8,6) {$b$};
\node at (8,0) {$k=3$};

\end{tikzpicture}

\end{center}
\caption{$(k)$-handles for $x$, with $k=1, 2, 3$.}
\label{x-handles}
\end{figure}
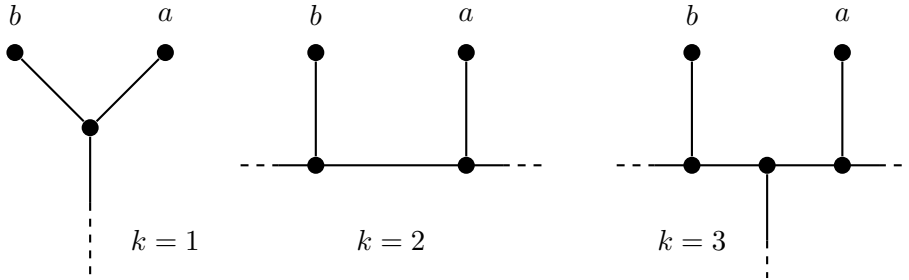

Instead of using fixed points of $y$, one can also define handles, and a similar joining operation, based on fixed points $a$ and $b$ of the standard generator $x$ of order $3$. In this case, we will say that they form a $(k)$-{\sl handle for\/} $x$ if $b=a(yx)^ky$. 
We may assume that $k=1, 2$ or $3$, as illustrated in Figure~\ref{x-handles}.

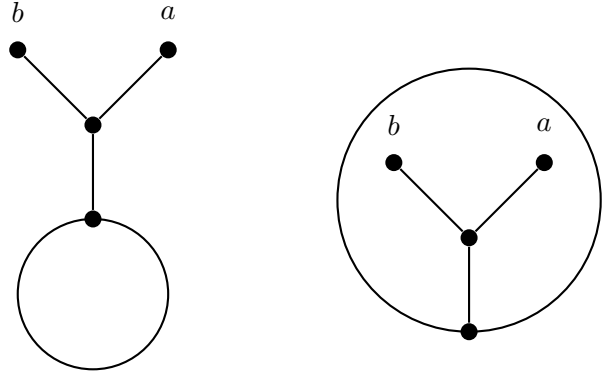
\begin{figure}[h!]
\begin{center}
\begin{tikzpicture}[scale=0.5, inner sep=0.8mm]

\node (a) at (2,5) [shape=circle, fill=black] {};
\node (b) at (-2,5) [shape=circle, fill=black] {};
\node (c) at (0,3) [shape=circle, fill=black] {};
\node (d) at (0,0.5) [shape=circle, fill=black] {};
\draw [thick] (c) to (d);
\draw [thick] (a) to (c) to (b);
\node at (2,6) {$a$};
\node at (-2,6) {$b$};
\draw [thick] (2,-1.5) arc (0:360:2);

\node (A) at (12,2) [shape=circle, fill=black] {};
\node (B) at (8,2) [shape=circle, fill=black] {};
\node (C) at (10,0) [shape=circle, fill=black] {};
\node (D) at (10,-2.5) [shape=circle, fill=black] {};
\draw [thick] (C) to (D);
\draw [thick] (A) to (C) to (B);
\node at (12,3) {$a$};
\node at (8,3) {$b$};
\draw [thick] (13.5,1) arc (0:360:3.5);

\end{tikzpicture}

\end{center}
\caption{The unique dessin with a $(1)$-handle for $x$.}
\label{1x-handle}
\end{figure}

If $k=1$, it is easy to see that the dessin is unique: it must be as in Figure~\ref{1x-handle}, shown with the fixed points in the outer or inner face, since no further vertices, edges or faces can be added. This dessin, of degree $8$ and genus $0$, corresponds to the natural action of ${\rm PSL}_2(7)$ on the projective line over $\mathbb F_7$, with point stabiliser $H\cong{\rm C}_7\rtimes{\rm C}_3$.

Figure~\ref{2,3-handles} shows that the cases $k=2$ and $k=3$ are equivalent: a $(2)$-handle for $x$ must be as shown on the left, so by transposing the roles of $a$ and $b$ it can also be regarded as a $(3)$-handle as on the right, and vice versa. We will therefore always assume that $k=1$, as in Figure~\ref{1x-handle}, or $k=2$; we will call these trivial and non-trivial $x$-handles respectively.

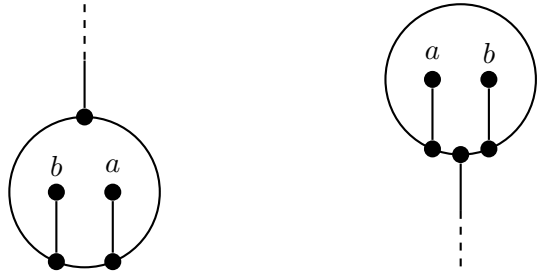
\begin{figure}[h!]
\begin{center}
\begin{tikzpicture}[scale=0.5, inner sep=0.8mm]

\draw [thick] (2,0) arc (0:360:2);
\node (a) at (0.75,0) [shape=circle, fill=black] {};
\node (a') at (-0.75,0) [shape=circle, fill=black] {};
\node (b) at (0,2) [shape=circle, fill=black] {};
\node (c) at (0.75,-1.85) [shape=circle, fill=black] {};
\node (c') at (-0.75,-1.85) [shape=circle, fill=black] {};

\draw [thick] (c) to (a);
\draw [thick] (c') to (a');
\draw [thick] (b) to (0,3.5);
\draw [thick, dashed] (0,5) to (0,3.5);

\node at (0.75,0.7) {$a$};
\node at (-0.75,0.7) {$b$};


\draw [thick] (12,3) arc (0:360:2);
\node (A) at (10.75,3) [shape=circle, fill=black] {};
\node (A') at (9.25,3) [shape=circle, fill=black] {};
\node (B) at (10,1) [shape=circle, fill=black] {};
\node (C) at (10.75,1.15) [shape=circle, fill=black] {};
\node (C') at (9.25,1.15) [shape=circle, fill=black] {};

\draw [thick] (C) to (A);
\draw [thick] (C') to (A');
\draw [thick] (B) to (10,-0.5);
\draw [thick, dashed] (10,-0.5) to (10,-2);

\node at (10.75,3.7) {$b$};
\node at (9.25,3.7) {$a$};

\end{tikzpicture}

\end{center}
\caption{Equivalence of $(2)$- and $(3)$-handles for $x$.}
\label{2,3-handles}
\end{figure}

If Hurwitz dessins $\mathcal D_i\;(i=1, 2, 3)$ have $(k)$-handles $(a_i, b_i)$ for $x$, with the same value of $k$, we can form an $x$-{\sl join} $\mathcal D_1(x)\mathcal D_2(x)\mathcal D_3$, by defining the cycles of $x$ and $y$ to be those they have on $\mathcal D_1$, $\mathcal D_2$ and $\mathcal D_3$, except that the six fixed points $a_i, b_i$ of $x$ become two 3-cycles $a=(a_1,a_2, a_3)$ and $b=(b_3,b_2, b_1)$ (note the reverse cyclic ordering of subscripts). The result is a connected Hurwitz dessin, and again this can be regarded as a connected sum operation, joining surfaces across cuts. 

The analogue of Theorem~\ref{additivity} for $x$-joins is as follows:

\begin{theorem}\label{3-additivity}
If Hurwitz dessins $\mathcal D_i\;(i=1,2,3)$ of degree $n_i$ and signature $(g_i;3^{[\alpha_i]},2^{[\beta_i]},7^{[\gamma_i]})$ have an $x$-join $\mathcal D$, then $\mathcal D$ has degree $n$ and signature  $(g;3^{[\alpha]},2^{[\beta]},7^{[\gamma]})$ where
\[n=\sum n_i,\; g=\sum g_i, \; \alpha=\sum\alpha_i-6, \; \beta=\sum\beta_i\;\hbox{and}\; \gamma=\sum\gamma_i.\]
\end{theorem}

\noindent{\sl Proof.} The proof is similar to that for Theorem~\ref{additivity}. \hfill$\square$

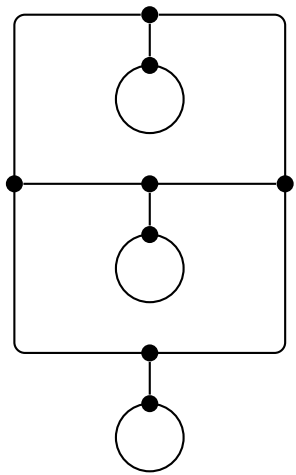
\begin{figure}[h!]
\begin{center}
\begin{tikzpicture}[scale=0.45, inner sep=0.8mm]

\draw [thick] (1,2.5) arc (0:360:1);
\node (a1) at (0,5) [shape=circle, fill=black] {};
\node (b1) at (0,3.5) [shape=circle, fill=black] {};
\draw [thick] (a1) to (b1);


\draw [thick] (1,-2.5) arc (0:360:1);
\node (a2) at (0,0) [shape=circle, fill=black] {};
\node (b2) at (0,-1.5) [shape=circle, fill=black] {};
\draw [thick] (a2) to (b2);


\draw [thick] (1,-7.5) arc (0:360:1);
\node (a3) at (0,-5) [shape=circle, fill=black] {};
\node (b3) at (0,-6.5) [shape=circle, fill=black] {};
\draw [thick] (a3) to (b3);


\node (d) at (-4,0) [shape=circle, fill=black] {};
\node (d') at (4,0) [shape=circle, fill=black] {};
\draw [thick, rounded corners] (a1) to (-4,5) to (-4,-5) to (4,-5) to (4,5) to (a1);
\draw [thick] (d) to (d');

\end{tikzpicture}

\end{center}
\caption{The trivial $x$-join, with $k=1$.}
\label{24-dessin}
\end{figure}

In the trivial case $k=1$, the only possible $x$-join involves three copies ${\mathcal D}_i$ of the dessin $\mathcal D$ of degree $8$ in Figure~\ref{1x-handle}. The resulting dessin $\mathcal D^*=\mathcal D(x)\mathcal D(x)\mathcal D$ is shown in Figure~\ref{24-dessin}. It has degree $24$ and genus $0$, and is a $3$-sheeted regular covering of $\mathcal D$. The monodromy group of $\mathcal D^*$ is described in the Example in \S5.1. Since neither $x$ nor $y$ has fixed points, this dessin cannot be used to form further joins.
 
\begin{figure}[h!]
\begin{center}
\begin{tikzpicture}[scale=0.45, inner sep=0.8mm]

\draw [thick] (2,5) arc (0:360:2);
\node (a1) at (2,5) [shape=circle, fill=black] {};
\node (a1') at (-2,5) [shape=circle, fill=black] {};
\node (b1) at (0,3) [shape=circle, fill=black] {};

\draw [thick] (b1) to (0,4);
\draw [thick, dashed] (0,5) to (0,4);


\draw [thick] (2,0) arc (0:360:2);
\node (a2) at (2,0) [shape=circle, fill=black] {};
\node (a2') at (-2,0) [shape=circle, fill=black] {};
\node (b2) at (0,-2) [shape=circle, fill=black] {};

\draw [thick] (b2) to (0,-1);
\draw [thick, dashed] (0,0) to (0,-1);


\draw [thick] (2,-5) arc (0:360:2);
\node (a3) at (2,-5) [shape=circle, fill=black] {};
\node (a3') at (-2,-5) [shape=circle, fill=black] {};
\node (b3) at (0,-7) [shape=circle, fill=black] {};

\draw [thick] (b3) to (0,-6);
\draw [thick, dashed] (0,-5) to (0,-6);


\node (d) at (4,0) [shape=circle, fill=black] {};
\node (d') at (-4,0) [shape=circle, fill=black] {};
\draw [thick, rounded corners] (a1) to (4,5) to (d) to (4,-5) to (a3);
\draw [thick, rounded corners] (a1') to (-4,5) to (d') to (-4,-5) to (a3');
\draw [thick] (a2) to (d);
\draw [thick] (a2') to (d');

\node at (7,5) {$\mathcal D_1$};
\node at (7,0) {$\mathcal D_2$};
\node at (7,-5) {$\mathcal D_3$};

\node at (5,0) {$a$};
\node at (-5,0) {$b$};

\end{tikzpicture}

\end{center}
\caption{The general dessin $\mathcal D_1(x)\mathcal D_2(x)\mathcal D_3$, with $k=2$.}
\label{x-join}
\end{figure}
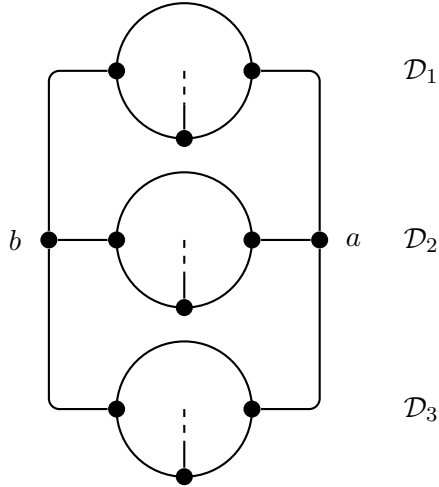

In the more general case $k=2$, a typical $x$-join is shown in Figure~\ref{x-join}, with the handles of each dessin drawn in the outer face. If we restrict attention to the faces containing the points $a_i$ and $b_i$, we have a 3-sheeted covering, branched over the points $a$ and $b$. Whereas in the case of $y$-joins we had the same monodromy permutation, a 2-cycle, at each of the two the branch-points, here we use two mutually inverse 3-cycles; the important point is that in both cases the product of the two monodromy permutations is the identity. The resulting dessin depends only on the chosen cyclic order of the three dessins $\mathcal D_i$; reversing this gives the mirror image of $\overline\mathcal D_1(x)\overline\mathcal D_2(x)\overline\mathcal D_3$, where $\overline\mathcal D_i$ denotes the mirror image of $\mathcal D_i$. It follows that if each $\mathcal D_i$ is reflexible the result is completely independent of the order of the three dessins. If the three dessins $\mathcal D_i$ are mutually isomorphic, by isomorphisms which match up the selected handles, then  $\mathcal D_1(x)\mathcal D_2(x)\mathcal D_3$ has an automorphism of order 3, permuting the dessins $\mathcal D_i$ cyclically and fixing $a$ and $b$.



\subsection{The triple cover $\mathcal D(x)\mathcal D(x)\mathcal D$}

If a dessin $\mathcal D$ has an $x$-handle, one can form the join $\mathcal D^*=\mathcal D(x)\mathcal D(x)\mathcal D$, 
a triple cover of $\mathcal D$ branched over two points, corresponding to the fixed points $a$ and $b$ of $x$ in the chosen handle. If $\mathcal D$ has degree $n$ and signature $(g;3^{[\alpha]},2^{[\beta]},7^{[\gamma]})$, then $\mathcal D^*$ has degree $3n$ and signature $(3g;\alpha^{[3\alpha-6]},3^{[3\beta]},7^{[3\gamma]})$.

The covering $\mathcal D^*\to\mathcal D$ is regular, induced by the obvious automorphism of order $3$ of $\mathcal D^*$, so $M^*$ is a normal subgroup of $M$. (It is, in fact, the normal closure in $M$ of all its standard generators apart from the two generators $X_j$ corresponding to the fixed points of $X$ in the chosen handle.) Thus $M^*$ contains the subgroup $M'M^3$ generated by the commutators and cubes of the elements of $M$. Now the signature of $M$ shows that $M/M'M^3$ is an elementary abelian $3$-group of rank $r=2g+\alpha-1$, so $M$ has $(3^r-1)/2$ subgroups of index $3$, each corresponding to a regular triple covering of $\mathcal D$. This leads to the following analogue of Lemma~\ref{uniqueD*}:

\begin{lemma}\label{uniquetriplecover}
Suppose that a planar Hurwitz dessin $\mathcal D$ has exactly two fixed points for $X$, forming a $(k)$-handle; then $\mathcal D^*=\mathcal D(x)\mathcal D(x)\mathcal D$ is the unique Hurwitz dessin which is a triple covering of $\mathcal D$. If, in addition, the point stabiliser $H$ in the monodromy group $G$ of $\mathcal D$ has a subgroup $H^*$ of index $3$, then $\mathcal D^*$ also has monodromy group $G$, with point stabiliser $H^*$. \hfill$\square$
\end{lemma}

\noindent{\bf Example}. This applies to the dessin $\mathcal D$ in Figure~\ref{1x-handle}, showing that the monodromy group of $\mathcal D^*$ (shown in Figure~\ref{24-dessin}) is ${\rm PSL}_2(7)$, like that of $\mathcal D$, but now acting on the cosets of a Sylow $7$-subgroup $H^*\cong {\rm C}_7$.

\medskip

However, if $g>0$ or $\alpha>2$ then $\mathcal D^*$ is one of several triple coverings of $\mathcal D$, and its imprimitive monodromy group could be isomorphic to $G$ or a proper covering group of $G$ contained in ${\rm S}_3\wr G$. The latter must happen if $H$ has no subgroup of index $3$, for example if $G$ is ${\rm A}_n$ acting naturally.



\subsection{A non-trivial example of an $x$-join}\label{xexample}

The triple cover $\mathcal A^*=\mathcal A(x)\mathcal A(x)\mathcal A$ of the dessin $\mathcal A$ of degree $14$ in Figures~\ref{AandA} and \ref{dessinsABC} is a Hurwitz dessin of degree 42 and genus 0. This is the dessin $\mathcal G$ based on Conder's diagram $G$ in~\cite{Conder-80} and shown in Figure~\ref{dessinG}. It is a 3-sheeted regular covering of $\mathcal A$, branched over the two points corresponding to the fixed points $a$ and $b$ of $x$, so that $\mathcal A\cong \mathcal G/{\rm C}_3$ where the group ${\rm C}_3$ of automorphisms of $\mathcal G$ permutes the three copies of $\mathcal A$ cyclically. It follows from Lemma~\ref{uniquetriplecover} that, like $\mathcal A$, this dessin $\mathcal G$ also has monodromy group ${\rm PSL}_2(13)$, but now in its imprimitive representation of degree 42 on the cosets of a dihedral subgroup $H^*\cong {\rm D}_{13}={\rm C}_{13}\rtimes {\rm C}_2$.

This example illustrates an important general point about $x$-joins, analogous to a point made earlier about $y$-joins (see Figure~\ref{wrongSS} and the accompanying text). In Figure~\ref{dessinG} there appear to be several different ways of joining three copies of the dessin $\mathcal A$ together: some of them could first be rotated through a half-turn, or equivalently reflected in a horizontal axis, so that the `Y' containing the $y$-handle is inverted. However, unless we do this to all three copies, this would break the rule that vertices $a_i$ in the $x$-handles must be identified with each other to form a new vertex $a$, with the same applying to the vertices $b_i$: there cannot be any `mixing' of vertices $a_i$ and $b_i$. Equivalently, there must be combinatorial isomorphisms between the oriented faces containing the $x$-handles in the three dessins, so that $xy$ still has order 7 after the joining operation. Incorrect identifications here could result in faces of valency other than 1 or 7, giving a transitive representation of $\Delta(3,2,\infty)$, isomorphic to the modular group ${\rm PSL}_2(\mathbb Z)$, rather than $\Delta$ (see~\S\ref{modular}).


\subsection{$x$-handles from ${\rm PSL}_2(q)$}

The following is the analogue of Theorem~\ref{PSL2qnathandles} for $x$-handles:

\begin{theorem} A dessin corresponding to the natural representation of the group $G={\rm PSL}_2(q)$ has an $x$-handle if and only if $q=13$, with ${\rm tr}(z)=\pm 5$.
\end{theorem}

\noindent{\sl Proof.} If an $x$-handle exists then $x$ has two fixed points, so it must be an elliptic element, of order dividing $(q-1)/2$ where $q\equiv 1$  mod~$(3)$. By applying an automorphism of $G$ we may assume that $x$ has the form $t\mapsto \omega t$, fixing $a=0$ and $b=\infty$, where $\omega^2+\omega+1=0$, or equivalently, that $x$ corresponds to the pair of elements
\[\pm\left( \begin{array}{cc}
\omega^2 & 0 \\ 
0 & \omega \\
\end{array}\right)\]
of ${\rm SL}_2(q)$. An element of $G$ has order $2$ if and only if it has trace $0$, so the most general form for $y$ is
\[\pm\left( \begin{array}{cc}
a & b \\ 
c & -a \\
\end{array}\right)\]
where $a^2+bc+1=0$. It then follows that $xy$ has trace
\[\pm t=\pm a(\omega^2-\omega)=\pm a\sqrt{-3}.\]

As before, $t=\pm\tau_j$ where $\tau_j=\lambda^j+\lambda^{-j}$ where $\lambda$ is a primitive 7th root of $1$, and the elements $\pm\tau_j$ are respectively the roots of the polynomials
\[p_{\pm}(t):=t(t^2-2)\pm(t^2-1).\]

Now $0$ and $\infty$ form an $x$-handle if and only if $yxyxy:0\mapsto\infty$. Writing $x$ and $y$ as above we find that this is equivalent to
\[\frac{abc(\omega-1)}{bc+a^2\omega^2}=a.\]
Now $a\ne 0$ (since otherwise $\langle x, y\rangle\cong {\rm S}_3$), so this simplifies to
\[3a^2=\omega-2.\]
Thus $t^2=-3a^2=2-\omega$, so 
\[0=-p_{\pm}(t)=t(t^2-2)\pm(t^2-1)=-t\omega\pm(1-\omega)\]
and hence
\[t^2\omega^2=(1-\omega)^2.\]
This simplifies to $\omega=3$, so that $3^2+3+1=0$. Thus $p=13$, so $q=13$ with $t=\pm 5$ and $a=\pm 3$. \hfill$\square$

\medskip

The dessin characterised by this theorem is $\mathcal A$, shown in Figures~\ref{AandA} and \ref{dessinsABC}. 

\begin{corollary}
Suppose that a Hurwitz dessin $\mathcal D$ has a monodromy group $G={\rm PSL}_2(q)$ in a representation which covers the natural representation. Then $\mathcal D$ has an $x$-handle if and only if
\begin{itemize}
\item $q=13$, with ${\rm tr}(z)=\pm 5$, and
\item the point stabilisers have order divisible by $3$.
\end{itemize}
\end{corollary}

\noindent{\sl Proof.} The proof is similar to that of Corollary~\ref{PSLycorollary} for $y$-handles, with the prime $3$ replacing $2$. \hfill$\square$

\medskip

The possible point-stabilisers in ${\rm PSL}_2(13)$ satisfying these conditions are therefore those isomorphic to  ${\rm C}_{13}\rtimes{\rm C}_6$,  ${\rm C}_{13}\rtimes{\rm C}_3$, ${\rm C}_6$, or ${\rm C}_3$. The first two correspond to the planar dessins $\mathcal A$ and $\mathcal A(1)\mathcal A$ considered earlier (see Figures~\ref{AandA}, \ref{dessinsABC} and~\ref{dessinsA*B*C*}); these have one $x$-handle and two, respectively, as have the dessins corresponding to the last two, which are $13$-sheeted branched coverings of these, of genus $2$ and $4$.

As with $y$-handles, this leaves open the question of which other representation of Hurwitz groups ${\rm PSL}_2(q)$ provide $x$-handles. As before, it is sufficient to consider those handles arising from primitive representations.


\subsection{More dessins with $x$-handles.}

It seems to be rather unusual for a Hurwitz dessin to contain an $x$-handle, at least among those dessins of relatively low degree. For instance, of the 14 dessins $\mathcal A,\ldots, \mathcal N$ based on Conder's diagrams $A,\ldots, N$ in~\cite{Conder-80}, only $\mathcal A$ has one: indeed, in each of the remaining 13 cases, $x$ has at most one fixed point (see Table~\ref{dessinsAtoN} in the Appendix).

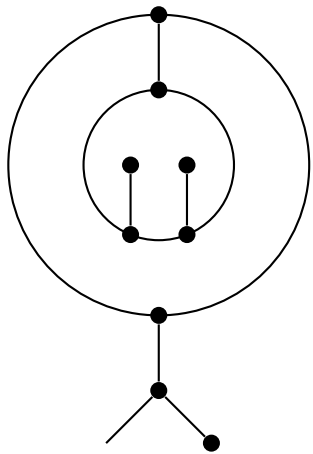
\begin{figure}[h!]
\begin{center}
\begin{tikzpicture}[scale=0.5, inner sep=0.8mm]

\draw [thick] (2,0) arc (0:360:2);
\node (a) at (0.75,0) [shape=circle, fill=black] {};
\node (a') at (-0.75,0) [shape=circle, fill=black] {};
\node (b) at (0,2) [shape=circle, fill=black] {};
\node (c) at (0.75,-1.85) [shape=circle, fill=black] {};
\node (c') at (-0.75,-1.85) [shape=circle, fill=black] {};

\draw [thick] (c) to (a);
\draw [thick] (c') to (a');

\draw [thick] (4,0) arc (0:360:4);
\node (d) at (0,4) [shape=circle, fill=black] {};
\node (e) at (0,-4) [shape=circle, fill=black] {};
\node (f) at (0,-6) [shape=circle, fill=black] {};
\node (g) at (1.4,-7.4) [shape=circle, fill=black] {};

\draw [thick] (b) to (d);
\draw [thick] (e) to (f) to (g);
\draw [thick] (f) to (-1.4,-7.4);

\end{tikzpicture}

\end{center}
\caption{A dessin of degree 21 with an $x$-handle.}
\label{21dessin}
\end{figure}

One simple technique for finding further examples is to start by drawing a single $x$-handle in a face, and systematically work outwards, adding one face at a time until the dessin closes up. The constraint that each face should have seven sides or one greatly reduces the number of possibilities at each stage. In addition to $\mathcal A$, one soon finds the following examples, all planar.

Figure~\ref{21dessin} and its mirror image illustrate a chiral pair of dessins of degree 21, each with one $x$-handle. One can number the half-edges, starting on the right of the top vertex, so that the standard generators for the dessin are
\[x=(1,2,3)(4,5,6)(7,8,9)(11,12,13)(15,16,17)(18,19,20)(10)(14)(21),\]
\[y=(1,17)(2,15)(3,4)(5,7)(6,11)(8,12)(9,10)(12,14)(16,18)(20,21)(19),\]
\[z=(1,16,20,21,19,18,15)(2,17,3,6,13,7,4)(5,9,10,8,12,14,11).\]
We then find from GAP that the monodromy group is ${\rm A}_{21}$. Alternatively, one can prove this by hand, using the fact that the commutator $[x,y]=x^{-1}yxy$ is a 19-cycle
\[(1,7,9,6,13,12,3,18,20,21,16,4,14,8,11,10,5,17,19);\]
this clearly implies that the group is primitive, and by~\cite{Jones-14} the only primitive permutation groups of degree $n$ containing an $(n-2)$-cycle are the symmetric group ${\rm S}_n$, the alternating group ${\rm A}_n$ for $n$ odd, and subgroups of ${\rm P\Gamma L}_2(q)$ containing ${\rm PGL}_2(q)$ where $n=q+1$ for some prime power $q$. Since Hurwitz groups are perfect, and 20 is not a prime power, we have the required result.

\begin{figure}[h!]
\begin{center}
\begin{tikzpicture}[scale=0.5, inner sep=0.8mm]

\draw [thick] (2,0) arc (0:360:2);
\node (a) at (0.75,0) [shape=circle, fill=black] {};
\node (a') at (-0.75,0) [shape=circle, fill=black] {};
\node (b) at (0,2) [shape=circle, fill=black] {};
\node (c) at (0.75,-1.85) [shape=circle, fill=black] {};
\node (c') at (-0.75,-1.85) [shape=circle, fill=black] {};

\draw [thick] (c) to (a);
\draw [thick] (c') to (a');

\draw [thick] (4,0) arc (0:360:4);
\node (d) at (0,4) [shape=circle, fill=black] {};
\node (e) at (0,-4) [shape=circle, fill=black] {};
\node (f) at (0,-6) [shape=circle, fill=black] {};

\draw [thick] (b) to (d);
\draw [thick] (e) to (f);

\draw [thick] (6,0) arc (0:360:6);
\node (g) at (6,0) [shape=circle, fill=black] {};
\node (g') at (-6,0) [shape=circle, fill=black] {};
\node (h) at (8,0) [shape=circle, fill=black] {};
\node (h') at (-8,0) [shape=circle, fill=black] {};

\draw [thick] (g) to (h);
\draw [thick] (g') to (h');

\end{tikzpicture}

\end{center}
\caption{A dessin of degree 28 with two $x$-handles.}
\label{28dessin}
\end{figure}
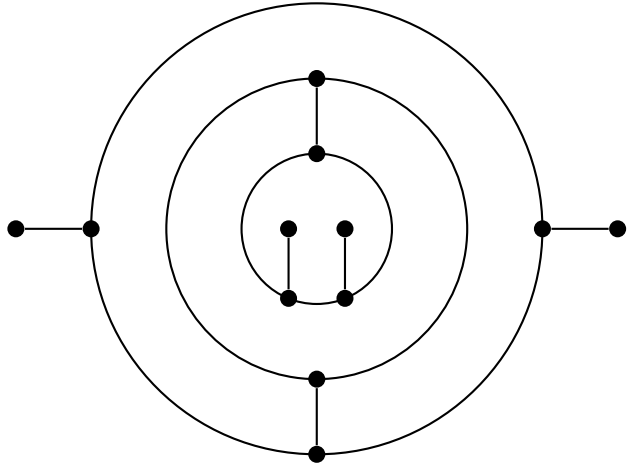

Figure~\ref{28dessin} illustrates a reflexible dessin of degree 28 with two $x$-handles. It is a redrawing of the $y$-join $\mathcal A(1)\mathcal A$ discussed earlier (see Figure~\ref{dessinsA*B*C*})
. Its monodromy group is ${\rm PSL}_2(13)$, in its imprimitive representation of degree 28 on the cosets of a subgroup ${\rm C}_{13}\rtimes{\rm C}_3$.

\begin{figure}[h!]
\begin{center}
\begin{tikzpicture}[scale=0.5, inner sep=0.8mm]

\draw [thick] (2,0) arc (0:360:2);
\node (a) at (0.75,0) [shape=circle, fill=black] {};
\node (a') at (-0.75,0) [shape=circle, fill=black] {};
\node (b) at (0,2) [shape=circle, fill=black] {};
\node (c) at (0.75,-1.85) [shape=circle, fill=black] {};
\node (c') at (-0.75,-1.85) [shape=circle, fill=black] {};

\draw [thick] (c) to (a);
\draw [thick] (c') to (a');

\draw [thick] (4,0) arc (0:360:4);
\node (d) at (0,4) [shape=circle, fill=black] {};
\node (e) at (0,-4) [shape=circle, fill=black] {};
\node (f) at (0,-6) [shape=circle, fill=black] {};

\draw [thick] (b) to (d);
\draw [thick] (e) to (f);

\draw [thick] (6,0) arc (0:360:6);
\node (g) at (6,0) [shape=circle, fill=black] {};
\node (g') at (-6,0) [shape=circle, fill=black] {};
\node (h) at (8,0) [shape=circle, fill=black] {};

\draw [thick] (g) to (h);
\draw [thick] (g') to (-8,0);
\draw [thick] (10,0) arc (0:360:1);

\end{tikzpicture}

\end{center}
\caption{A dessin of degree 29 with an $x$-handle.}
\label{29dessin}
\end{figure}
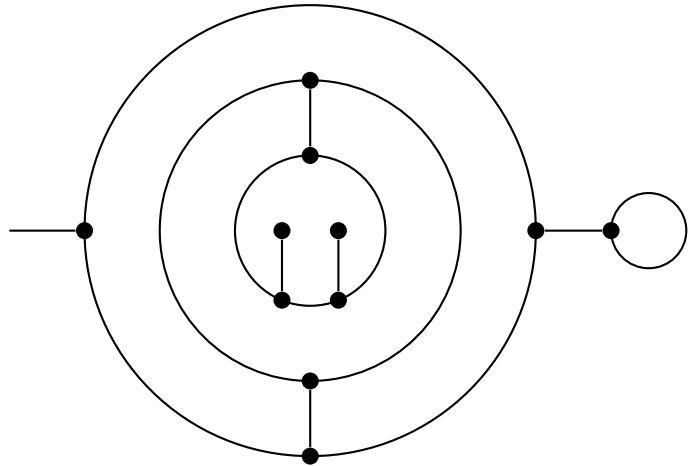

Figure~\ref{29dessin} and its mirror image illustrate a chiral pair of dessins of degree $n=29$, each with one $x$-handle. As in the case $n=21$, GAP tells us that the monodromy group is ${\rm A}_n$, and again the fact that $[x,y]$ is an $(n-2)$-cycle allows us to confirm this by hand.

\begin{figure}[h!]
\begin{center}
\begin{tikzpicture}[scale=0.5, inner sep=0.8mm]

\draw [thick] (2,0) arc (0:360:2);
\node (a) at (0.75,0) [shape=circle, fill=black] {};
\node (a') at (-0.75,0) [shape=circle, fill=black] {};
\node (b) at (0,2) [shape=circle, fill=black] {};
\node (c) at (0.75,-1.85) [shape=circle, fill=black] {};
\node (c') at (-0.75,-1.85) [shape=circle, fill=black] {};

\draw [thick] (c) to (a);
\draw [thick] (c') to (a');

\draw [thick] (4,0) arc (0:360:4);
\node (d) at (0,4) [shape=circle, fill=black] {};
\node (e) at (0,-4) [shape=circle, fill=black] {};
\node (f) at (0,-6) [shape=circle, fill=black] {};

\draw [thick] (b) to (d);
\draw [thick] (e) to (f);

\draw [thick] (6,0) arc (0:360:6);
\node (g) at (6,0) [shape=circle, fill=black] {};
\node (g') at (-6,0) [shape=circle, fill=black] {};
\node (h) at (8,0) [shape=circle, fill=black] {};
\node (h') at (-8,0) [shape=circle, fill=black] {};

\draw [thick] (g) to (h);
\draw [thick] (g') to (h');

\draw [thick] (8,0) arc (0:360:8);
\node (i) at (0,8) [shape=circle, fill=black] {};
\node (j) at (0,7) [shape=circle, fill=black] {};
\draw [thick] (i) to (j);
\node (k) at (1,-7.9) [shape=circle, fill=black] {};
\node (k') at (-1,-7.9) [shape=circle, fill=black] {};
\draw [thick] (k) to (1,-9.5);
\draw [thick] (k') to (-1,-9.5);

\end{tikzpicture}

\end{center}
\caption{A dessin of degree 42 with an $x$-handle.}
\label{42dessin}
\end{figure}
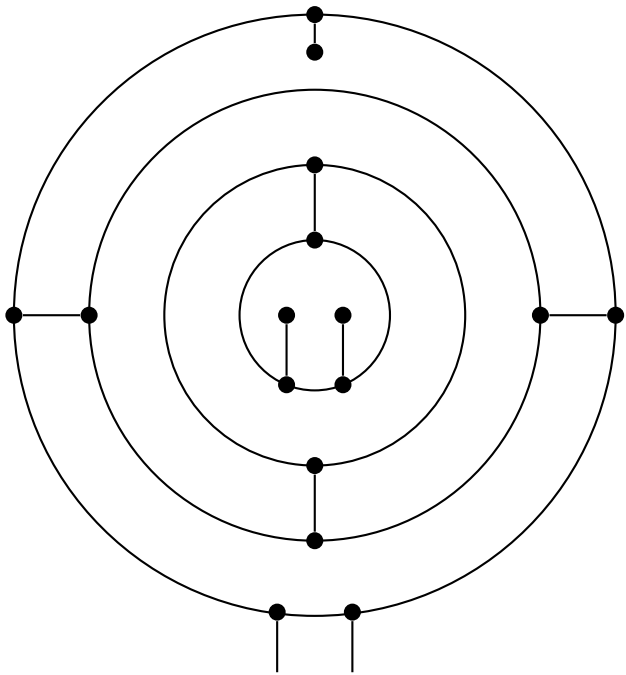

Figure~\ref{42dessin} illustrates a reflexible dessin $\mathcal D$ of degree 42 with an $x$-handle. According to GAP its monodromy group is ${\rm A}_{42}$. Using the $(2)$-handle in the outer face, we can form a dessin $\mathcal D_1=\mathcal D(1)\mathcal D$ of degree 84 with two $x$-handles. Using one of these handles we can form the triple cover $\mathcal D_2=\mathcal D_1(x)\mathcal D_1(x)\mathcal D_1$ of degree 252 with three $x$-handles. Iteration gives a sequence of dessins $\mathcal D_{n+1}=\mathcal D_n(x)\mathcal D_n(x)\mathcal D_n$ of degree $28.3^n$ with $3(3^{n-2}+1)/2$ $x$-handles, for all $n\ge 1$; they are all reflexible and planar. Many other similar iterations are possible, since if $\mathcal D$ is any dessin with at least one $x$-handle, and  $\mathcal D'$ and $\mathcal D''$ are any dessins with at least two each, then $\mathcal D(x)\mathcal D'(x)\mathcal D''$ also has at least two. Indeed, by taking them to have $k-1$, two and two $x$-handles, one can construct a dessin with $k$ $x$-handles for each $k\ge 2$.


\section{Handles and joins for the modular group}\label{modular}

In~\cite{Conder-81}, Conder extended the technique of joining coset diagrams which he and Higman had developed for $\Delta(3,2,7)$ to the triangle groups $\Delta(3,2,r)$ for all integers $r\ge 7$, obtaining similar results to those in the case $r=7$.

One can obtain a further extension to the modular group $\Gamma={\rm PSL}_2({\mathbb Z})=\Delta(3,2,\infty)\cong {\rm C}_3*{\rm C}_2$, by omitting the relation $Z^7=1$, so that faces of any valency are allowed, while vertices all have valency $3$ or $1$. In this case, any ordered pair of fixed points of $x$ or of $y$ in the same face form an $x$- or $y$-handle, so that three or two such dessins with handles can be joined as described earlier to produce another, but this time with no requirement that the handles should be in isomorphic faces.

In any such trivalent dessin, the order of $z$ is the level of the corresponding subgroup $M$ of $\Gamma$; in the case of a congruence subgroup, this coincides with the level as defined number-theoretically, by a result of Wohlfahrt~\cite{Woh}.

As an example of dessins rising from congruence subgroups, there is a natural action of $\Gamma$ by M\"obius transformations on $\mathbb P^1(\mathbb Q)$. For any prime $p$, reducing this action mod~$(p)$ gives the natural action of ${\rm PSL}_2(p)$ on $\mathbb P^1(\mathbb F_p)$. The generating set
\[X: t\mapsto\frac{1}{1-t},\quad Y:t\mapsto \frac{-1}{t},\quad Z: t\mapsto t+1\]
for $\Gamma$, with defining relations
\[X^3=Y^2=XYZ=1,\]
reduces mod~$(p)$ to give a generating triple $(x,y,z)$ of type $(3,2,p)$ for $G(p):={\rm PSL}_2(p)$, so the natural action of $G(p)$ on the projective line $\mathbb P^1(p)$, with stabilisers $H\cong{\rm C}_p\rtimes{\rm C}_{(p-1)/2}$ for odd $p$, gives a dessin $\mathcal P(p)$ of that type and of degree $p+1$ with monodromy group $G(p)$. There are two faces, of valencies $p$ and $1$; there are two free edges (giving a $y$-handle) or none as $p\equiv\pm 1$ mod~$(4)$, and for $p>3$ there are two $1$-valent vertices (giving an $x$-handle) or none as $p\equiv \pm 1$ mod~$(3)$.

\medskip

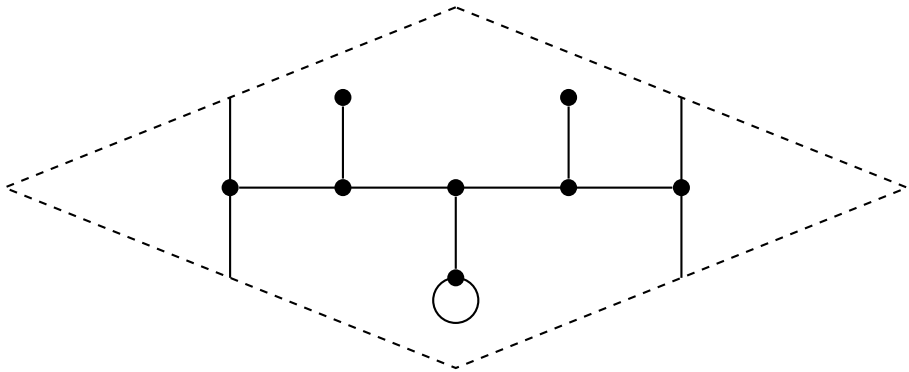
\begin{figure}[h!]
\begin{center}
\begin{tikzpicture}[scale=0.3, inner sep=0.8mm]

\node (a) at (-10,0) [shape=circle, fill=black] {};
\node (b) at (-5,0) [shape=circle, fill=black] {};
\node (c) at (0,0) [shape=circle, fill=black] {};
\node (b') at (5,0) [shape=circle, fill=black] {};
\node (a') at (10,0) [shape=circle, fill=black] {};
\node (f) at (0,-4) [shape=circle, fill=black] {};
\node (g) at (-5,4) [shape=circle, fill=black] {};
\node (g') at (5,4) [shape=circle, fill=black] {};

\draw [thick] (a) to (a');
\draw [thick] (c) to (f);
\draw [thick] (b) to (g);
\draw [thick] (b') to (g');
\draw [thick] (-10,4) to (-10,-4);
\draw [thick] (10,4) to (10,-4);

\draw [thick, dashed] (0,8) to (-20,0) to (0,-8) to (20,0) to (0,8);

\draw [thick] (1,-5) arc (0:360:1);

\end{tikzpicture}

\end{center}
\caption{The dessin $\mathcal P(19)$.}
\label{PSL2(19)}
\end{figure}

\noindent{\bf Example.} By Dirichlet's theorem there are infinitely many primes $p\equiv 7$ mod~$(12)$. For such primes, $\mathcal P(p)$ has $(p+5)/3$ vertices (two of valency $1$, corresponding to the fixed points of $x$, the rest of valency $3$), $(p+1)/2$ edges (none of them free), and two faces, so it has genus $g=(p-7)/12$. (See Figure~\ref{1x-handle} for $\mathcal P(7)$, and Figure~\ref{PSL2(19)} for $\mathcal P(19)$, with opposite sides of the outer parallelogram identified to form a torus.) The two vertices of valency $1$ both lie in the face of valency $p$, so they form a handle. Given three primes $p_i\equiv 7$ mod~$(12)$, the join $\mathcal D=\mathcal P(p_1)(x)\mathcal P(p_2)(x)\mathcal P(p_3)$ has degree $s+3$, where $s=p_1+p_2+p_3$, and has genus $(s-21)/12$.

If we take each $p_i=p$ we obtain a regular triple cover $\mathcal P(p)^*$ of $\mathcal P(p)$, with $G(p)$ as its monodromy group, now acting on the cosets of a subgroup $H^*\cong{\rm C}_p\rtimes{\rm C}_{(p-1)/6}$. However, if the primes $p_i$ are not all equal, then the monodromy group can be very different from the groups $G(p_i)$. For instance, the torus dessin $\mathcal P(7)(x)\mathcal P(7)(x)\mathcal P(19)$ shown in Figure~\ref{dessinP(7,7,19)} has monodromy group $G\cong {\rm A}_{36}$; this fact, confirmed by GAP, can be proved as follows, using the cycle structure $1^{[3]},7,9,17$ of $z$.

If $G$ is imprimitive it has $a$ blocks of size $b$, where $ab=36$ and $a, b>1$, so it can be embedded in the wreath product ${\rm S}_b\wr {\rm S}_a$. This has order $(b!)^aa!$ with $a, b\le 18$, so we require $\{a, b\}=\{18, 2\}$ for it to have order divisible by $17$. Since $G$ is generated by elements $x$ and $z$ of odd order it cannot have two blocks, so $a=18$ and $b=2$. However, if the group ${\rm S}_2\wr {\rm S}_{18}=({\rm S}_2)^{18}\rtimes {\rm S}_{18}$ has an element of order $7.9.17$ then so has its quotient ${\rm S}_{18}$; this is clearly false, so $G$ is primitive. Applying Jordan's Theorem (see~\S\ref{Jordan}) to the $7$-cycle $z^{153}$ gives $G\cong {\rm A}_{36}$ or ${\rm S}_{36}$, and as before the latter is impossible since $G=\langle x, z\rangle$.

\begin{figure}[h!]
\begin{center}
\begin{tikzpicture}[scale=0.3, inner sep=0.8mm]

\node (a) at (-10,-3) [shape=circle, fill=black] {};
\node (b) at (-5,-3) [shape=circle, fill=black] {};
\node (c) at (0,-3) [shape=circle, fill=black] {};
\node (b') at (5,-3) [shape=circle, fill=black] {};
\node (a') at (10,-3) [shape=circle, fill=black] {};
\node (f) at (0,-5) [shape=circle, fill=black] {};

\node (h) at (-5,2) [shape=circle, fill=black] {};
\node (h') at (5,2) [shape=circle, fill=black] {};
\node (i) at (0,7) [shape=circle, fill=black] {};
\node (j) at (0,2) [shape=circle, fill=black] {};
\node (k) at (0,5) [shape=circle, fill=black] {};
\node (l) at (0,0) [shape=circle, fill=black] {};

\draw [thick] (a) to (a');
\draw [thick] (c) to (f);
\draw [thick] (b) to (g);
\draw [thick] (b') to (g');
\draw [thick, rounded corners] (h) to (-5,7) to (5,7) to (h');
\draw [thick] (h) to (h');
\draw [thick] (i) to (k);
\draw [thick] (j) to (l);
\draw [thick] (-10,5.25) to (-10,-4.25);
\draw [thick] (10,5.25) to (10,-4.25);

\draw [thick, dashed] (0,10) to (-20,0.5) to (0,-9) to (20,0.5) to (0,10);

\draw [thick] (1,-6) arc (0:360:1);
\draw [thick] (1,4) arc (0:360:1);
\draw [thick] (1,-1) arc (0:360:1);

\end{tikzpicture}

\end{center}
\caption{The dessin $\mathcal P(7)(x)\mathcal P(7)(x)\mathcal P(19)$.}
\label{dessinP(7,7,19)}
\end{figure}


\section{Handles and joins for dessins of all types}

Here we briefly sketch how to define handles and joins for dessins of any type $(p,q,r)$. It is sufficient to define handles and joins only for $x$, generalising the construction given earlier for them in the case $p=3$, $q=2$.

We will now represent a dessin $\mathcal D$ in the usual way as a bipartite map, with $x$ and $y$ rotating edges around their incident black and white vertices, and with faces corresponding to the cycles of $z$. (Up to now, with $q=2$, we have omitted the white vertices, but now, with arbitrary $q$, we need them.) As before, a $(k)$-handle for $x$ is a pair $a, b$ of fixed points for $x$ with $b=a(yx)^ky$, so that the corresponding black vertices lie in the same face.

\begin{figure}[h!]
\begin{center}
\begin{tikzpicture}[scale=0.4, inner sep=0.8mm]

\draw [thick, dotted] (2,5) arc (0:360:2);
\node (a1) at (1.2,5) [shape=circle, draw] {};
\node (a1') at (-1.2,5) [shape=circle, draw] {};


\draw [thick, dotted] (2,0) arc (0:360:2);
\node (a2) at (1.2,0) [shape=circle, draw] {};
\node (a2') at (-1.2,0) [shape=circle, draw] {};


\draw [thick, dotted] (2,-8) arc (0:360:2);
\node (a3) at (1.2,-8) [shape=circle, draw] {};
\node (a3') at (-1.2,-8) [shape=circle, draw] {};


\node (d) at (4,0) [shape=circle, fill=black] {};
\node (d') at (-4,0) [shape=circle, fill=black] {};
\draw [thick, rounded corners] (a1) to (4,5) to (d) to (4,-8) to (a3);
\draw [thick, rounded corners] (a1') to (-4,5) to (d') to (-4,-8) to (a3');
\draw [thick] (a2) to (d);
\draw [thick] (a2') to (d');
\draw [thick, dotted] (0,-3) to (0,-5);
\draw [thick, dotted] (7,-3) to (7,-5);

\node at (7,5) {$\mathcal D_1$};
\node at (7,0) {$\mathcal D_2$};
\node at (7,-8) {$\mathcal D_d$};

\node at (5,0) {$a$};
\node at (-5,0) {$b$};

\end{tikzpicture}

\end{center}
\caption{The $x$-join $\mathcal D_1(x)\mathcal D_2(x)\cdots(x)\mathcal D_d$.}
\label{generalx-join}
\end{figure}
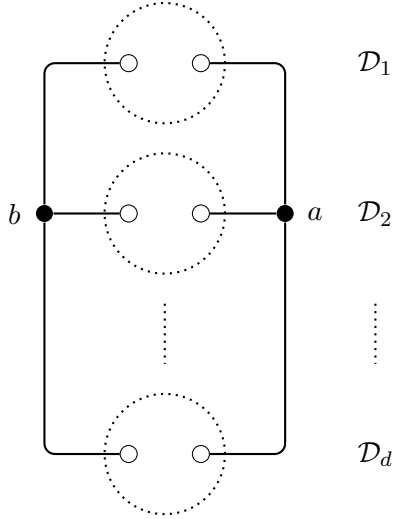

Let dessins $\mathcal D_i\;(i=1, 2, \ldots, d)$ of type $(p,q,r)$ have $(k)$-handles $(a_i, b_i)$ for $x$, with the same $k$, where $d$ divides $p$. An $x$-{\sl join} $\mathcal D_1(x)\mathcal D_2(x)\cdots(x)\mathcal D_d$ is formed by defining the cycles of $x$ and $y$ to be those they have in the dessins $\mathcal D_i$, except that the fixed points $a_i, b_i$ of $x$ become two cycles $a=(a_1,a_2, \ldots, a_d)$ and $b=(b_d, b_{d-1}\ldots, b_1)$. The result is a dessin $\mathcal D$ of type $(p,q,r)$, formed by joining the dessins $\mathcal D_i$  across cuts (see Figure~\ref{generalx-join}).  As before, this operation is additive in the degrees and genera of the dessins $\mathcal D_i$. The elliptic periods in the signature of $\mathcal D$ are those for the dessins $\mathcal D_i$, except that $2d$ periods $p$ corresponding to the fixed points $a_i, b_i$ are replaced with two periods $p/d$ corresponding to $a$ and $b$ (these can be omitted if $d=p$).

Everitt~\cite{Everitt} used a similar joining operation with $k=1$ for coset diagrams for Dyck groups (cocompact groups of genus $0$, of signature $(0;m_1,\ldots, m_r)$ for integers $m_i\ge 2$). This was an important ingredient in proving that each non-elementary Fuchsian group $\Gamma$ has almost all alternating groups as quotients. It is hoped to use a combination of these joining operations to extend the results in~\cite{Jones-18a, Jones-18b} on maximal subgroups from certain triangle groups to all such groups $\Gamma$.


\section{Appendix: dessins based on Conder's diagrams}


In~\cite{Conder-80}, Conder used $14$ coset diagrams $A,\ldots, N$ for subgroups of finite index $n$ in $\Delta$. As explained earlier, his diagrams can be interpreted as planar Hurwitz dessins $\mathcal A,\ldots, \mathcal N$ of degree $n$. We have already used the dessins $\mathcal A$, $\mathcal B$, $\mathcal C$, $\mathcal F$ and $\mathcal G$ as examples of Hurwitz dessins: see Figure~\ref{dessinsABC} for $\mathcal A$, $\mathcal B$ and $\mathcal C$, and Figures~\ref{dessinF} and \ref{dessinG} for $\mathcal F$ and $\mathcal G$; for completeness, we give all the dessins here, together with their basic properties (see also~\cite{JP}).

\begin{table}[ht]
\centering
\begin{tabular}{| p{1cm} | p{1.4cm} | p{3cm} | p{1.9cm} | p{1.2cm} |}
\hline
Dessin & degree $n$ & monodromy group $G$ & number of $(k)$-handles & $\alpha, \beta, \gamma$\\
\hline\hline
$\mathcal A$ & $14$ & ${\rm PSL}_2(13)$ & $1, 0, 0$ & $2, 2, 0$ \\
\hline
$\mathcal B$ & $15$ & ${\rm A}_{15}$ & $0, 1, 1$ & $0, 3, 1$ \\
\hline
$\mathcal C$ & $21$ & ${\rm PGL}_3(2)>{\rm D}_4$ & $1, 0, 1$ & $0, 5, 0$ \\
\hline
$\mathcal D$ & $22$ & ${\rm A}_{22}$ & $0, 1, 0$ & $1, 2, 1$ \\
\hline
$\mathcal E$ & $28$ & ${\rm PSL}_2(8)>{\rm D}_9$ & $1, 1, 0$ & $1, 4, 0$ \\
\hline
$\mathcal F$ & $30$ & ${\rm PSL}_2(29)$ & $0, 1, 0$ & $0, 2, 2$ \\
\hline
$\mathcal G$ & $42$ & ${\rm PSL}_2(13)>{\rm D}_{13}$ & $3, 0, 0$ & $0, 6, 0$ \\
\hline
$\mathcal H$ & $42$ & ${\rm A}_{42}$ & $1, 0, 1$ & $0, 6, 0$ \\
\hline
$\mathcal I$ & $57$ & ${\rm A}_{57}$ & $0, 2, 0$ & $0, 5, 1$ \\
\hline
$\mathcal J$ & $72$ & $({\rm S}_2\wr{\rm A}_{36})\cap{\rm A}_{72}$ & $2, 0, 0$ & $0, 4, 2$ \\
\hline
$\mathcal K$ & $72$ & ${\rm A}_{72}$ & $1, 0, 0$ & $0, 4, 2$ \\
\hline
$\mathcal L$ & $102$ & ${\rm A}_{102}$ & $0, 1, 0$ & $0, 2, 4$ \\
\hline
$\mathcal M$ & $108$ & ${\rm A}_{108}$ & $1, 1, 0$ & $0, 4, 3$ \\
\hline
$\mathcal N$ & $108$ & ${\rm A}_{108}$ & $1, 0, 1$ & $0, 4, 3$ \\
\hline

\end{tabular}
\caption{Dessins $\mathcal A,\ldots, \mathcal N$ based on Conder's diagrams $A,\ldots, N$}
\label{dessinsAtoN}
\end{table}

 \medskip

In Table~\ref{dessinsAtoN} we give the degree $n$ and the monodromy group $G$ (in most cases obtained by using GAP) for each dessin. When $G$ is ${\rm A}_n$ or ${\rm PSL}_2(n-1)$, it always acts naturally; such an action is primitive, so the automorphism group of the dessin is trivial. For $\mathcal C, \mathcal E$ and $\mathcal G$, where $G$ does not act naturally, we give a point-stabiliser $H$ in the form $G>H$. The orientation-preserving automorphism group of the dessin is trivial for $\mathcal C$ and $\mathcal E$, and isomorphic to ${\rm C}_3$ and ${\rm C}_2$ for $\mathcal G$ and $\mathcal J$. We also give the number of disjoint $(k)$-handles for $k=1,2,3$ and the numbers $\alpha$, $\beta$ and $\gamma$ of fixed points of $x, y$ and $z$, which determine the signature of the dessin (they all have genus $0$).

\bigskip

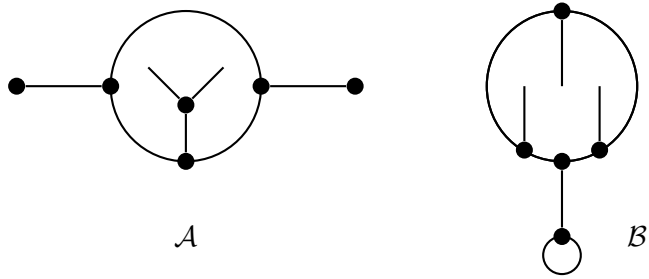
\begin{figure}[h]
\begin{center}
\begin{tikzpicture}[scale=0.5, inner sep=0.8mm]

\draw [thick] (2,0) arc (0:360:2);
\node (a) at (2,0) [shape=circle, fill=black] {};
\node (a') at (-2,0) [shape=circle, fill=black] {};
\node (b) at (0,-2) [shape=circle, fill=black] {};
\node (c) at (0,-0.5) [shape=circle, fill=black] {};
\node (d) at (4.5,0) [shape=circle, fill=black] {};
\node (d') at (-4.5,0) [shape=circle, fill=black] {};

\draw [thick] (b) to (c);
\draw [thick] (1,0.5) to (c) to (-1,0.5);
\draw [thick] (a) to (d);
\draw [thick] (a') to (d');

\node at (0,-4)  {$\mathcal A$};


\draw [thick] (12,0) arc (0:360:2);
\node (b2) at (10,-2) [shape=circle, fill=black] {};
\node (c2) at (11,-1.7) [shape=circle, fill=black] {};
\node (c'2) at (9,-1.7) [shape=circle, fill=black] {};
\node (d2) at (10,2) [shape=circle, fill=black] {};

\draw [thick] (c2) to (11,0);
\draw [thick] (c'2) to (9,0);
\draw [thick] (d2) to (10,0);

\node (e2) at (10,-4) [shape=circle, fill=black] {};
\draw [thick] (b2) to (e2);
\draw [thick] (12,0) arc (0:360:2);
\draw [thick] (10.5,-4.5) arc (0:360:0.5);

\node at (12,-4)  {$\mathcal B$};

\end{tikzpicture}

\end{center}
\caption{The dessins $\mathcal A$ and $\mathcal B$.}
\label{dessinsAB}
\end{figure}


\begin{figure}[h!]
\begin{center}
\begin{tikzpicture}[scale=0.5, inner sep=0.8mm]

\draw [thick] (2,0) arc (0:360:2);
\node (a) at (2,0) [shape=circle, fill=black] {};
\node (a') at (-2,0) [shape=circle, fill=black] {};
\node (b) at (0,2) [shape=circle, fill=black] {};
\node (c) at (1,1.7) [shape=circle, fill=black] {};
\node (c') at (-1,1.7) [shape=circle, fill=black] {};

\draw [thick] (c) to (1,0);
\draw [thick] (c') to (-1,0);
\draw [thick] (b) to (0,3.5);

\node (d) at (0,5) [shape=circle, fill=black] {};
\draw [rounded corners, thick] (a) to (3,0) to (3,5) to (-3,5) to (-3,0) to (a');

\node (e) at (0,7) [shape=circle, fill=black] {};
\draw [thick] (d) to (e);
\draw [thick] (1,8) to (e) to (-1,8);

\node at (3,-2)  {$\mathcal C$};


\draw [thick] (12,4) arc (0:360:2);
\node (a4) at (12,4) [shape=circle, fill=black] {};
\node (a'4) at (8,4) [shape=circle, fill=black] {};
\node (b4) at (10,6) [shape=circle, fill=black] {};
\node (c4) at (10.75,2.15) [shape=circle, fill=black] {};
\node (c'4) at (9.25,2.15) [shape=circle, fill=black] {};
\node (d4) at (10,8) [shape=circle, fill=black] {};
\node (e4) at (10,0) [shape=circle, fill=black] {};
\node (f4) at (10,-2) [shape=circle, fill=black] {};

\draw [thick] (c4) to (10.75,4);
\draw [thick] (c'4) to (9.25,4);
\draw [thick] (b4) to (d4);
\draw [thick] (e4) to (f4);
\draw [rounded corners, thick] (a4) to (13,4) to (13,-2) to (7,-2) to (7,4) to (a'4);

\draw [thick] (10.5,8.5) arc (0:360:0.5);

\node at (15,-2)  {$\mathcal D$};

\end{tikzpicture}

\end{center}
\caption{The dessins $\mathcal C$ and $\mathcal D$.}
\label{dessinsCD}
\end{figure}
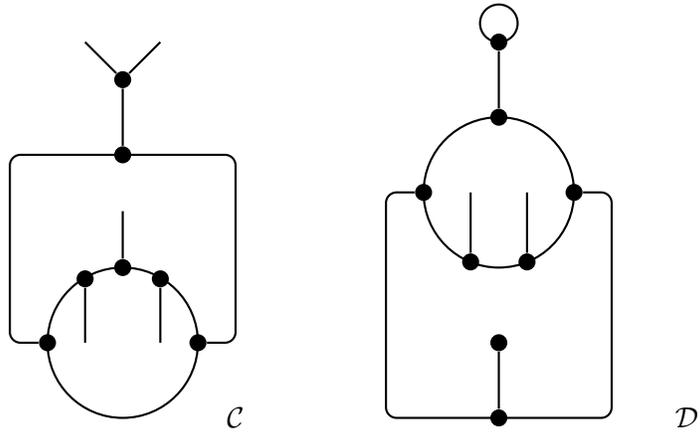

\newpage


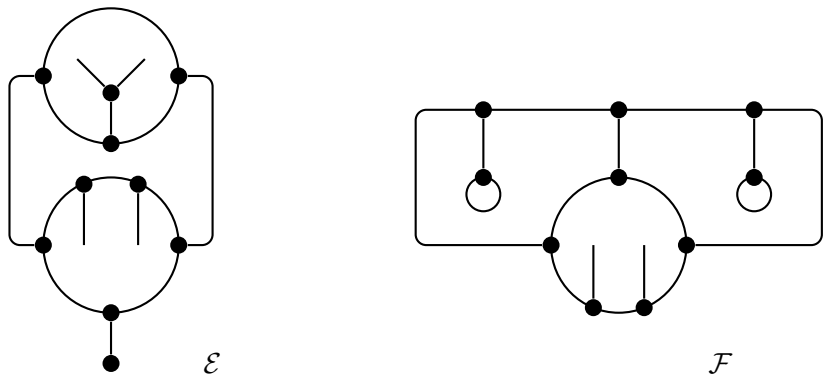
\begin{figure}[h!]
\begin{center}
\begin{tikzpicture}[scale=0.45, inner sep=0.8mm]

\draw [thick] (2,5) arc (0:360:2);
\draw [thick] (2,0) arc (0:360:2);
\node (a) at (2,5) [shape=circle, fill=black] {};
\node (a') at (-2,5) [shape=circle, fill=black] {};
\node (b) at (0,3) [shape=circle, fill=black] {};
\node (c) at (0,4.5) [shape=circle, fill=black] {};
\node (d) at (0.8,1.8) [shape=circle, fill=black] {};
\node (d') at (-0.8,1.8) [shape=circle, fill=black] {};
\node (e) at (2,0) [shape=circle, fill=black] {};
\node (e') at (-2,0) [shape=circle, fill=black] {};
\node (f) at (0,-2) [shape=circle, fill=black] {};
\node (g) at (0,-3.5) [shape=circle, fill=black] {};

\draw [rounded corners, thick] (a) to (3,5) to (3,0) to (e);
\draw [rounded corners, thick] (a') to (-3,5) to (-3,0) to (e');
\draw [thick] (b) to (c);
\draw [thick] (1,5.5) to (c) to (-1,5.5);
\draw [thick] (d) to (0.8,0);
\draw [thick] (d') to (-0.8,0);
\draw [thick] (f) to (g);

\node at (3,-3.5) {$\mathcal E$};


\draw [thick] (17,0) arc (0:360:2);
\node (a6) at (17,0) [shape=circle, fill=black] {};
\node (a'6) at (13,0) [shape=circle, fill=black] {};
\node (b6) at (15,2) [shape=circle, fill=black] {};
\node (c6) at (15.75,-1.85) [shape=circle, fill=black] {};
\node (c'6) at (14.25,-1.85) [shape=circle, fill=black] {};

\draw [thick] (c6) to (15.75,0);
\draw [thick] (c'6) to (14.25,0);

\node (d6) at (15,4) [shape=circle, fill=black] {};
\draw [thick] (b6) to (d6);
\node (e6) at (19,4) [shape=circle, fill=black] {};
\node (e'6) at (11,4) [shape=circle, fill=black] {};
\node (f6) at (19,2) [shape=circle, fill=black] {};
\node (f'6) at (11,2) [shape=circle, fill=black] {};

\draw [thick] (19.5,1.5) arc (0:360:0.5);
\draw [thick] (11.5,1.5) arc (0:360:0.5);

\draw [rounded corners, thick] (a6) to (21,0) to (21,4) to (9,4) to (9,0) to (a'6);
\draw [thick] (e6) to (f6);
\draw [thick] (e'6) to (f'6);

\node at (18,-3.5) {$\mathcal F$};

\end{tikzpicture}

\end{center}
\caption{The dessins $\mathcal E$ and $\mathcal F$.}
\label{DessinsEF}
\end{figure}


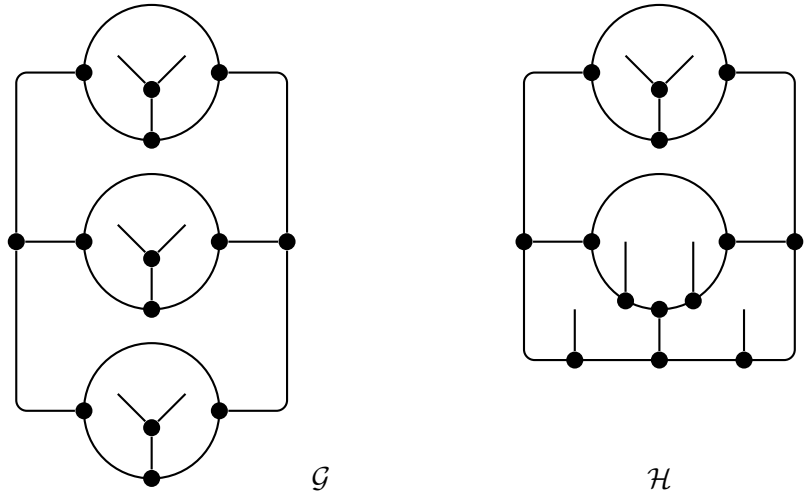
\begin{figure}[h!]
\begin{center}
\begin{tikzpicture}[scale=0.45, inner sep=0.8mm]

\draw [thick] (2,5) arc (0:360:2);
\node (a1) at (2,5) [shape=circle, fill=black] {};
\node (a1') at (-2,5) [shape=circle, fill=black] {};
\node (b1) at (0,3) [shape=circle, fill=black] {};
\node (c1) at (0,4.5) [shape=circle, fill=black] {};

\draw [thick] (b1) to (c1);
\draw [thick] (1,5.5) to (c1) to (-1,5.5);

\draw [thick] (2,0) arc (0:360:2);
\node (a2) at (2,0) [shape=circle, fill=black] {};
\node (a2') at (-2,0) [shape=circle, fill=black] {};
\node (b2) at (0,-2) [shape=circle, fill=black] {};
\node (c2) at (0,-0.5) [shape=circle, fill=black] {};

\draw [thick] (b2) to (c2);
\draw [thick] (1,0.5) to (c2) to (-1,0.5);

\draw [thick] (2,-5) arc (0:360:2);
\node (a3) at (2,-5) [shape=circle, fill=black] {};
\node (a3') at (-2,-5) [shape=circle, fill=black] {};
\node (b3) at (0,-7) [shape=circle, fill=black] {};
\node (c3) at (0,-5.5) [shape=circle, fill=black] {};

\draw [thick] (b3) to (c3);
\draw [thick] (1,-4.5) to (c3) to (-1,-4.5);

\node (d) at (4,0) [shape=circle, fill=black] {};
\node (d') at (-4,0) [shape=circle, fill=black] {};
\draw [rounded corners, thick] (a1) to (4,5) to (d) to (4,-5) to (a3);
\draw [rounded corners, thick] (a1') to (-4,5) to (d') to (-4,-5) to (a3');
\draw [thick] (a2) to (d);
\draw [thick] (a2') to (d');

\node at (5, -7) {$\mathcal G$};


\draw [thick] (17,5) arc (0:360:2);
\node (a1) at (17,5) [shape=circle, fill=black] {};
\node (a1') at (13,5) [shape=circle, fill=black] {};
\node (b1) at (15,3) [shape=circle, fill=black] {};
\node (c1) at (15,4.5) [shape=circle, fill=black] {};

\draw [thick] (b1) to (c1);
\draw [thick] (16,5.5) to (c1) to (14,5.5);

\draw [thick] (17,0) arc (0:360:2);
\node (a2) at (17,0) [shape=circle, fill=black] {};
\node (a2') at (13,0) [shape=circle, fill=black] {};
\node (c2) at (16,-1.75) [shape=circle, fill=black] {};
\node (c2') at (14,-1.75) [shape=circle, fill=black] {};
\node (d2) at (15,-2) [shape=circle, fill=black] {};

\draw [thick] (c2) to (16,0);
\draw [thick] (c2') to (14,0);

\node (e) at (19,0) [shape=circle, fill=black] {};
\node (e') at (11,0) [shape=circle, fill=black] {};
\node (f) at (17.5,-3.5) [shape=circle, fill=black] {};
\node (f') at (12.5,-3.5) [shape=circle, fill=black] {};
\node (g) at (15,-3.5) [shape=circle, fill=black] {};
\draw [thick] (g) to (d2);
\draw [rounded corners, thick] (a1) to (19,5) to (19,-3.5) to (11,-3.5) to (11,5) to (a1');
\draw [thick] (a2) to (e);
\draw [thick] (a2') to (e');
\draw [thick] (f) to (17.5,-2);
\draw [thick] (f') to (12.5,-2);

\node at (15, -7) {$\mathcal H$};

\end{tikzpicture}

\end{center}
\caption{The dessins $\mathcal G$ and $\mathcal H$.}
\label{DessinsGH}
\end{figure}

\newpage


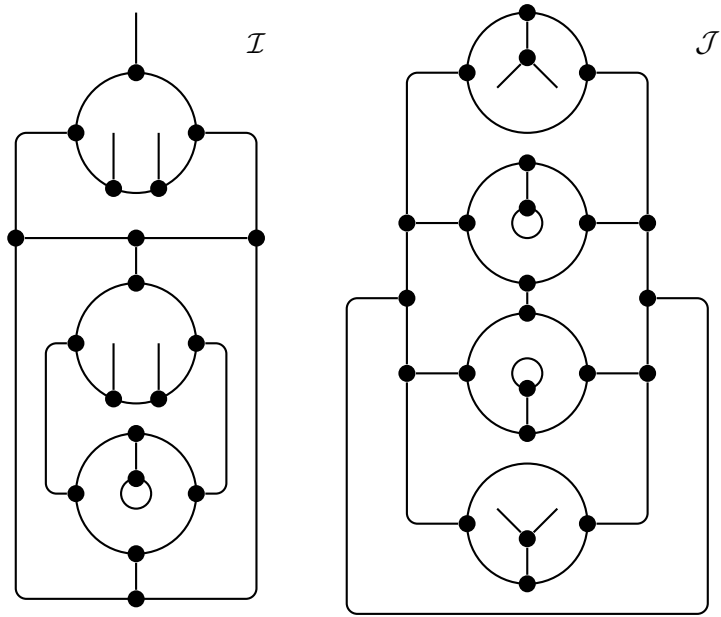
\begin{figure}[h!]
\begin{center}
\begin{tikzpicture}[scale=0.4, inner sep=0.8mm]

\draw [thick] (-11,8) arc (0:360:2);
\node (A1) at (-11,8) [shape=circle, fill=black] {};
\node (A1') at (-15,8) [shape=circle, fill=black] {};
\node (B1) at (-13,10) [shape=circle, fill=black] {};
\node (C1) at (-12.25,6.15) [shape=circle, fill=black] {};
\node (C1') at (-13.75,6.15) [shape=circle, fill=black] {};

\draw [thick] (C1) to (-12.25,8);
\draw [thick] (C1') to (-13.75,8);
\draw [thick] (B1) to (-13,12);


\draw [thick] (-11,1) arc (0:360:2);
\node (A2) at (-11,1) [shape=circle, fill=black] {};
\node (A2') at (-15,1) [shape=circle, fill=black] {};
\node (B2) at (-13,3) [shape=circle, fill=black] {};
\node (C2) at (-12.25,-0.85) [shape=circle, fill=black] {};
\node (C2') at (-13.75,-0.85) [shape=circle, fill=black] {};

\draw [thick] (C2) to (-12.25,1);
\draw [thick] (C2') to (-13.75,1);


\draw [thick] (-11,-4) arc (0:360:2);
\node (A3) at (-11,-4) [shape=circle, fill=black] {};
\node (A3') at (-15,-4) [shape=circle, fill=black] {};
\node (B3) at (-13,-2) [shape=circle, fill=black] {};
\node (C3) at (-13,-3.5) [shape=circle, fill=black] {};
\node (D3) at (-13,-6) [shape=circle, fill=black] {};

\draw [thick] (B3) to (C3);
\draw [thick] (-12.5,-4) arc (0:360:0.5);


\node (E) at (-13,4.5) [shape=circle, fill=black] {};
\node (F) at (-9,4.5) [shape=circle, fill=black] {};
\node (F') at (-17,4.5) [shape=circle, fill=black] {};
\node (G) at (-13,-7.5) [shape=circle, fill=black] {};
\draw [thick] (F) to (F');
\draw [rounded corners, thick] (A1) to (-9,8) to (-9,-7.5) to (-17,-7.5) to (-17,8) to (A1');
\draw [thick] (E) to (B2);
\draw [thick] (G) to (D3);
\draw [rounded corners, thick] (A2) to (-10,1) to (-10,-4) to (A3);
\draw [rounded corners, thick] (A2') to (-16,1) to (-16,-4) to (A3');

\node at (-9,11) {$\mathcal I$};


\draw [thick] (2,10) arc (0:360:2);
\node (a1) at (2,10) [shape=circle, fill=black] {};
\node (a1') at (-2,10) [shape=circle, fill=black] {};
\node (b1) at (0,12) [shape=circle, fill=black] {};
\node (c1) at (0,10.5) [shape=circle, fill=black] {};

\draw [thick] (b1) to (c1);
\draw [thick] (1,9.5) to (c1) to (-1,9.5);


\draw [thick] (2,5) arc (0:360:2);
\node (a2) at (2,5) [shape=circle, fill=black] {};
\node (a2') at (-2,5) [shape=circle, fill=black] {};
\node (b2) at (0,7) [shape=circle, fill=black] {};
\node (c2) at (0,5.5) [shape=circle, fill=black] {};
\node (d2) at (0,3) [shape=circle, fill=black] {};

\draw [thick] (b2) to (c2);
\draw [thick] (0.5,5) arc (0:360:0.5);


\draw [thick] (2,0) arc (0:360:2);
\node (a3) at (2,0) [shape=circle, fill=black] {};
\node (a3') at (-2,0) [shape=circle, fill=black] {};
\node (b3) at (0,-2) [shape=circle, fill=black] {};
\node (c3) at (0,-0.5) [shape=circle, fill=black] {};
\node (d3) at (0,2) [shape=circle, fill=black] {};

\draw [thick] (b3) to (c3);
\draw [thick] (0.5,0) arc (0:360:0.5);


\draw [thick] (2,-5) arc (0:360:2);
\node (a4) at (2,-5) [shape=circle, fill=black] {};
\node (a4') at (-2,-5) [shape=circle, fill=black] {};
\node (b4) at (0,-7) [shape=circle, fill=black] {};
\node (c4) at (0,-5.5) [shape=circle, fill=black] {};

\draw [thick] (b4) to (c4);
\draw [thick] (1,-4.5) to (c4) to (-1,-4.5);


\node (e2) at (4,5) [shape=circle, fill=black] {};
\node (e2') at (-4,5) [shape=circle, fill=black] {};
\node (e3) at (4,0) [shape=circle, fill=black] {};
\node (e3') at (-4,0) [shape=circle, fill=black] {};
\draw [thick] (a2) to (e2);
\draw [thick] (a2') to (e2');
\draw [thick] (a3) to (e3);
\draw [thick] (a3') to (e3');
\draw [rounded corners, thick] (a1) to (4,10) to (e2) to (e3) to (4,-5) to (a4);
\draw [rounded corners, thick] (a1') to (-4,10) to (e2') to (e3') to (-4,-5) to (a4');


\node (f) at (4,2.5) [shape=circle, fill=black] {};
\node (f') at (-4,2.5) [shape=circle, fill=black] {};
\draw [rounded corners, thick] (f) to (6,2.5) to (6,-8) to (-6,-8) to (-6,2.5) to (f');
\draw [thick] (d2) to (d3);

\node at (6,11) {$\mathcal J$};

\end{tikzpicture}

\end{center}
\caption{The dessins $\mathcal I$ and $\mathcal J$.}
\label{DessinsIJ}
\end{figure}


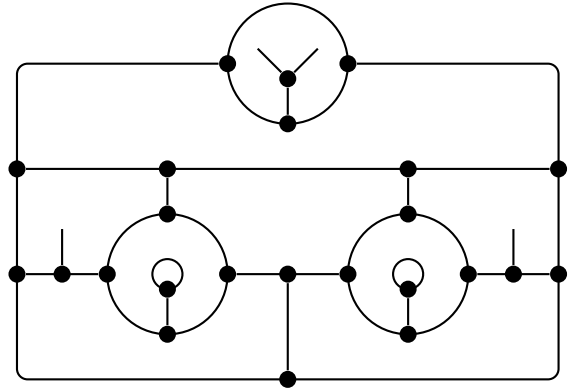
\begin{figure}[h!]
\begin{center}
\begin{tikzpicture}[scale=0.4, inner sep=0.8mm]

\node (x) at (0,11) {};

\draw [thick] (2,7) arc (0:360:2);
\node (a) at (2,7) [shape=circle, fill=black] {};
\node (a') at (-2,7) [shape=circle, fill=black] {};
\node (b) at (0,5) [shape=circle, fill=black] {};
\node (c) at (0,6.5) [shape=circle, fill=black] {};

\draw [thick] (b) to (c);
\draw [thick] (1,7.5) to (c) to (-1,7.5);


\draw [thick] (-2,0) arc (0:360:2);
\node (a1) at (-2,0) [shape=circle, fill=black] {};
\node (a1') at (-6,0) [shape=circle, fill=black] {};
\node (b1) at (-4,-2) [shape=circle, fill=black] {};
\node (c1) at (-4,-0.5) [shape=circle, fill=black] {};
\node (d1) at (-4,2) [shape=circle, fill=black] {};

\draw [thick] (b1) to (c1);
\draw [thick] (-3.5,0) arc (0:360:0.5);

\draw [thick] (6,0) arc (0:360:2);
\node (a2) at (6,0) [shape=circle, fill=black] {};
\node (a2') at (2,0) [shape=circle, fill=black] {};
\node (b2) at (4,-2) [shape=circle, fill=black] {};
\node (c2) at (4,-0.5) [shape=circle, fill=black] {};
\node (d2) at (4,2) [shape=circle, fill=black] {};

\draw [thick] (b2) to (c2);
\draw [thick] (4.5,0) arc (0:360:0.5);


\node (e) at (4,3.5) [shape=circle, fill=black] {};
\node (f) at (9,3.5) [shape=circle, fill=black] {};
\node (g) at (9,0) [shape=circle, fill=black] {};
\node (h) at (7.5,0) [shape=circle, fill=black] {};
\node (e') at (-4,3.5) [shape=circle, fill=black] {};
\node (f') at (-9,3.5) [shape=circle, fill=black] {};
\node (g') at (-9,0) [shape=circle, fill=black] {};
\node (h') at (-7.5,0) [shape=circle, fill=black] {};
\node (i) at (0,-3.5) [shape=circle, fill=black] {};
\node (j) at (0,0) [shape=circle, fill=black] {};

\draw [thick] (f) to (f');
\draw [thick] (g) to (a2);
\draw [thick] (g') to (a1');
\draw [thick] (a1) to (a2');
\draw [thick] (e) to (d2);
\draw [thick] (e') to (d1);
\draw [thick] (i) to (j);
\draw [thick] (h) to (7.5,1.5);
\draw [thick] (h') to (-7.5,1.5);
\draw [rounded corners, thick] (a') to (-9,7) to (-9,-3.5) to (9,-3.5) to (9,7) to (a);

\end{tikzpicture}

\end{center}
\caption{The dessin $\mathcal K$.}
\label{dessinK}
\end{figure}


\begin{figure}[h!]
\begin{center}
\begin{tikzpicture}[scale=0.4, inner sep=0.8mm]

\draw [thick] (13,0) arc (0:360:2);
\node (a1) at (13,0) [shape=circle, fill=black] {};
\node (a'1) at (9,0) [shape=circle, fill=black] {};
\node (b1) at (11,-2) [shape=circle, fill=black] {};
\node (c1) at (11,-0.5) [shape=circle, fill=black] {};
\node (d1) at (11,2) [shape=circle, fill=black] {};

\draw [thick] (b1) to (c1);
\draw [thick] (11.5,0) arc (0:360:0.5);


\draw [thick] (-9,0) arc (0:360:2);
\node (a2) at (-9,0) [shape=circle, fill=black] {};
\node (a'2) at (-13,0) [shape=circle, fill=black] {};
\node (b2) at (-11,-2) [shape=circle, fill=black] {};
\node (c2) at (-11,-0.5) [shape=circle, fill=black] {};
\node (d2) at (-11,2) [shape=circle, fill=black] {};

\draw [thick] (b2) to (c2);
\draw [thick] (-10.5,0) arc (0:360:0.5);


\draw [thick] (2,0) arc (0:360:2);
\node (a) at (2,0) [shape=circle, fill=black] {};
\node (a') at (-2,0) [shape=circle, fill=black] {};
\node (b) at (0,2) [shape=circle, fill=black] {};
\node (c) at (0.75,-1.85) [shape=circle, fill=black] {};
\node (c') at (-0.75,-1.85) [shape=circle, fill=black] {};

\draw [thick] (c) to (0.75,0);
\draw [thick] (c') to (-0.75,0);


\draw [thick] (6,-8) arc (0:360:2);
\node (a3) at (6,-8) [shape=circle, fill=black] {};
\node (a'3) at (2,-8) [shape=circle, fill=black] {};
\node (b3) at (4,-6) [shape=circle, fill=black] {};
\node (c3) at (4,-7.5) [shape=circle, fill=black] {};
\node (d3) at (4,-10) [shape=circle, fill=black] {};

\draw [thick] (b3) to (c3);
\draw [thick] (4.5,-8) arc (0:360:0.5);


\draw [thick] (-2,-8) arc (0:360:2);
\node (a4) at (-2,-8) [shape=circle, fill=black] {};
\node (a'4) at (-6,-8) [shape=circle, fill=black] {};
\node (b4) at (-4,-6) [shape=circle, fill=black] {};
\node (c4) at (-4,-7.5) [shape=circle, fill=black] {};
\node (d4) at (-4,-10) [shape=circle, fill=black] {};

\draw [thick] (b4) to (c4);
\draw [thick] (-3.5,-8) arc (0:360:0.5);


\draw [thick] (a) to (a'1);
\draw [thick] (a') to (a2);

\node (e) at (4,0) [shape=circle, fill=black] {};
\node (e') at (-4,0) [shape=circle, fill=black] {};
\draw [rounded corners, thick] (e) to (4,-4) to (-4,-4) to (e');

\node (f) at (0,-4) [shape=circle, fill=black] {};
\node (g) at (0,-8) [shape=circle, fill=black] {};
\draw [thick] (f) to (g);
\draw [thick] (a'3) to (g) to (a4);

\node (h) at (7.5,0) [shape=circle, fill=black] {};
\node (h') at (-7.5,0) [shape=circle, fill=black] {};
\draw [rounded corners, thick] (h) to (7.5,-8) to (a3);
\draw [rounded corners, thick] (h') to (-7.5,-8) to (a'4);

\node (i) at (4,-12) [shape=circle, fill=black] {};
\node (i') at (-4,-12) [shape=circle, fill=black] {};
\draw [rounded corners, thick] (a1) to (14.5,0) to (14.5,-12) to (i) to (i') to (-14.5,-12) to (-14.5,0) to (a'2);
\draw [thick] (i) to (d3);
\draw [thick] (i') to (d4);

\node (j) at (0,4) [shape=circle, fill=black] {};
\draw [thick] (j) to (b);
\draw [rounded corners, thick] (d1) to (11,4) to (j) to (-11,4) to (d2);

\end{tikzpicture}

\end{center}
\caption{The dessin $\mathcal L$.}
\label{dessinL}
\end{figure}
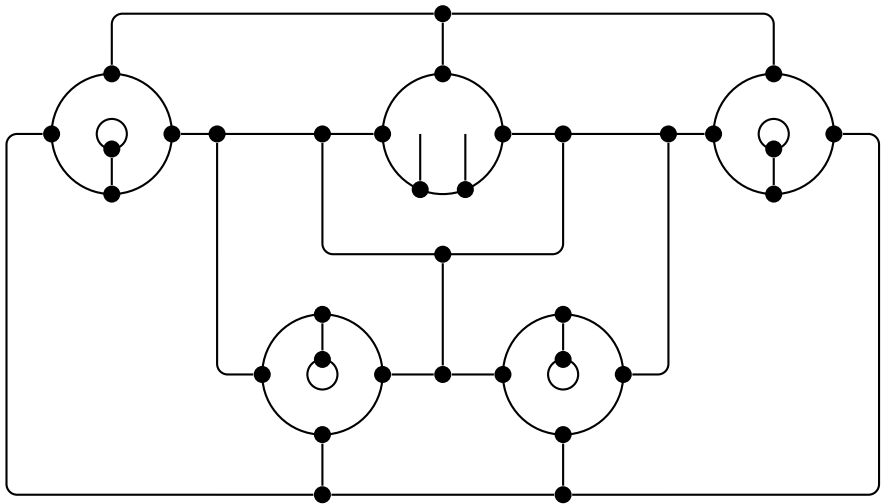

\newpage


\begin{figure}[h!]
\begin{center}
\begin{tikzpicture}[scale=0.37, inner sep=0.8mm]

\draw [thick] (10,0) arc (0:360:2);
\node (a1) at (10,0) [shape=circle, fill=black] {};
\node (a'1) at (6,0) [shape=circle, fill=black] {};
\node (b1) at (8,-2) [shape=circle, fill=black] {};
\node (c1) at (8,-0.5) [shape=circle, fill=black] {};
\node (d1) at (8,2) [shape=circle, fill=black] {};

\draw [thick] (b1) to (c1);
\draw [thick] (8.5,0) arc (0:360:0.5);


\draw [thick] (-6,0) arc (0:360:2);
\node (a2) at (-6,0) [shape=circle, fill=black] {};
\node (a'2) at (-10,0) [shape=circle, fill=black] {};
\node (b2) at (-8,-2) [shape=circle, fill=black] {};
\node (c2) at (-8,-0.5) [shape=circle, fill=black] {};
\node (d2) at (-8,2) [shape=circle, fill=black] {};

\draw [thick] (b2) to (c2);
\draw [thick] (-7.5,0) arc (0:360:0.5);


\draw [thick] (2,-3) arc (0:360:2);
\node (a3) at (2,-3) [shape=circle, fill=black] {};
\node (a'3) at (-2,-3) [shape=circle, fill=black] {};
\node (b3) at (0,-5) [shape=circle, fill=black] {};
\node (c3) at (0,-3.5) [shape=circle, fill=black] {};

\draw [thick] (b3) to (c3);
\draw [thick] (1,-2.5) to (c3) to (-1,-2.5);


\draw [thick] (2,4) arc (0:360:2);
\node (a4) at (2,4) [shape=circle, fill=black] {};
\node (a'4) at (-2,4) [shape=circle, fill=black] {};
\node (b4) at (0,2) [shape=circle, fill=black] {};
\node (c4) at (0.75,5.85) [shape=circle, fill=black] {};
\node (c'4) at (-0.75,5.85) [shape=circle, fill=black] {};

\draw [thick] (c4) to (0.75,4);
\draw [thick] (c'4) to (-0.75,4);


\draw [thick] (2,11) arc (0:360:2);
\node (a5) at (2,11) [shape=circle, fill=black] {};
\node (a'5) at (-2,11) [shape=circle, fill=black] {};
\node (b5) at (0,13) [shape=circle, fill=black] {};
\node (c5) at (0,11.5) [shape=circle, fill=black] {};
\node (d5) at (0,9) [shape=circle, fill=black] {};

\draw [thick] (b5) to (c5);
\draw [thick] (0.5,11) arc (0:360:0.5);


\draw [thick] (a'1) to (a2);
\node (e) at (0,0) [shape=circle, fill=black] {};
\node (f) at (4,0) [shape=circle, fill=black] {};
\node (f') at (-4,0) [shape=circle, fill=black] {};
\node (g) at (4,-3) [shape=circle, fill=black] {};
\node (g') at (-4,-3) [shape=circle, fill=black] {};
\node (h) at (4,-6) [shape=circle, fill=black] {};
\node (h') at (-4,-6) [shape=circle, fill=black] {};
\draw [thick] (e) to (b4);
\draw [thick] (f) to (h);
\draw [thick] (f') to (h');
\draw [thick] (a3) to (g);
\draw [thick] (a'3) to (g');


\node (i) at (12,0) [shape=circle, fill=black] {};
\node (i') at (-12,0) [shape=circle, fill=black] {};
\draw [thick] (a1) to (i);
\draw [thick] (a'2) to (i');
\draw [rounded corners, thick] (a5) to (12,11) to (12,-6) to (-12,-6) to (-12,11) to (a'5);


\node (j) at (0,7.5) [shape=circle, fill=black] {};
\node (k) at (8,4) [shape=circle, fill=black] {};
\node (k') at (-8,4) [shape=circle, fill=black] {};
\draw [rounded corners, thick] (d1) to (8,7.5) to (-8,7.5) to (d2);
\draw [thick] (k) to (a4);
\draw [thick] (k') to (a'4);
\draw [thick] (j) to (d5);

\end{tikzpicture}

\end{center}
\caption{The dessin $\mathcal M$.}
\label{dessinM}
\end{figure}
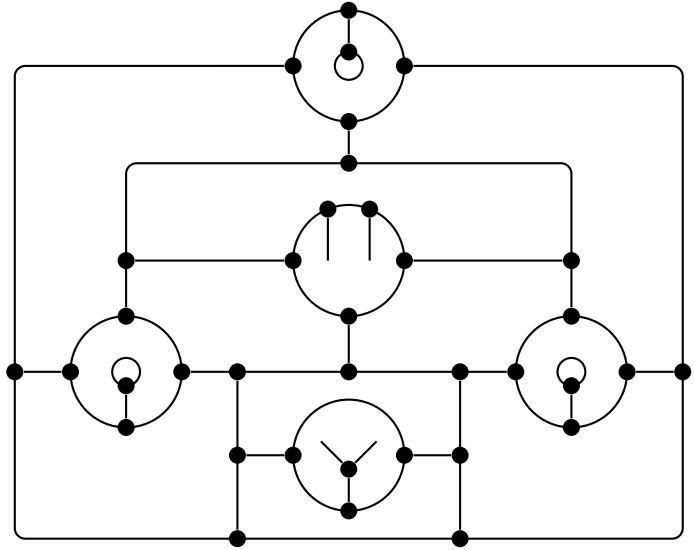


\begin{figure}[h!]
\begin{center}
\begin{tikzpicture}[scale=0.37, inner sep=0.8mm]

\draw [thick] (2,0) arc (0:360:2);
\node (a1) at (2,0) [shape=circle, fill=black] {};
\node (a'1) at (-2,0) [shape=circle, fill=black] {};
\node (b1) at (0,-2) [shape=circle, fill=black] {};
\node (c1) at (1,-1.7) [shape=circle, fill=black] {};
\node (c'1) at (-1,-1.7) [shape=circle, fill=black] {};

\draw [thick] (c1) to (1,0);
\draw [thick] (c'1) to (-1,0);


\draw [thick] (2,5) arc (0:360:2);
\node (a2) at (2,5) [shape=circle, fill=black] {};
\node (a'2) at (-2,5) [shape=circle, fill=black] {};
\node (b2) at (0,3) [shape=circle, fill=black] {};
\node (c2) at (0,4.5) [shape=circle, fill=black] {};
\node (d2) at (0,7) [shape=circle, fill=black] {};

\draw [thick] (b2) to (c2);
\draw [thick] (0.5,5) arc (0:360:0.5);


\draw [thick] (12,0) arc (0:360:2);
\node (a3) at (12,0) [shape=circle, fill=black] {};
\node (a'3) at (8,0) [shape=circle, fill=black] {};
\node (b3) at (10,-2) [shape=circle, fill=black] {};
\node (c3) at (10,-0.5) [shape=circle, fill=black] {};
\node (d3) at (10,2) [shape=circle, fill=black] {};

\draw [thick] (b3) to (c3);
\draw [thick] (10.5,0) arc (0:360:0.5);


\draw [thick] (-8,0) arc (0:360:2);
\node (a4) at (-8,0) [shape=circle, fill=black] {};
\node (a'4) at (-12,0) [shape=circle, fill=black] {};
\node (b4) at (-10,-2) [shape=circle, fill=black] {};
\node (c4) at (-10,-0.5) [shape=circle, fill=black] {};
\node (d4) at (-10,2) [shape=circle, fill=black] {};

\draw [thick] (b4) to (c4);
\draw [thick] (-9.5,0) arc (0:360:0.5);


\draw [thick] (2,-8) arc (0:360:2);
\node (a5) at (2,-8) [shape=circle, fill=black] {};
\node (a'5) at (-2,-8) [shape=circle, fill=black] {};
\node (b5) at (0,-10) [shape=circle, fill=black] {};
\node (c5) at (0,-8.5) [shape=circle, fill=black] {};

\draw [thick] (b5) to (c5);
\draw [thick] (1,-7.5) to (c5) to (-1,-7.5);


\draw [thick] (a1) to (a'3);
\draw [thick] (a'1) to (a4);
\node (e) at (4,0) [shape=circle, fill=black] {};
\node (e') at (-4,0) [shape=circle, fill=black] {};
\node (f) at (6,0) [shape=circle, fill=black] {};
\node (f') at (-6,0) [shape=circle, fill=black] {};
\node (g) at (6,-4) [shape=circle, fill=black] {};
\node (g') at (-6,-4) [shape=circle, fill=black] {};
\draw [thick] (f) to (g);
\draw [thick] (f') to (g');
\node (h) at (0,-4) [shape=circle, fill=black] {};
\node (i) at (13.5,-4) [shape=circle, fill=black] {};
\node (i') at (-13.5,-4) [shape=circle, fill=black] {};
\draw [thick] (i) to (i');
\draw [thick] (h) to (b1);


\draw [rounded corners, thick] (a2) to (4,5) to (e);
\draw [rounded corners, thick] (a'2) to (-4,5) to (e');
\node (j) at (0,9) [shape=circle, fill=black] {};
\draw [thick] (j) to (d2);
\draw [rounded corners, thick] (d3) to (10,9) to (-10,9) to (d4);
\node (k) at (13.5,0) [shape=circle, fill=black] {};
\node (k') at (-13.5,0) [shape=circle, fill=black] {};
\draw [thick] (a3) to (k);
\draw [thick] (a'4) to (k');
\draw [rounded corners, thick] (a5) to (13.5,-8) to (13.5,10) to (-13.5,10) to (-13.5,-8) to (a'5);

\end{tikzpicture}

\end{center}
\caption{The dessin $\mathcal N$.}
\label{dessinN}
\end{figure}
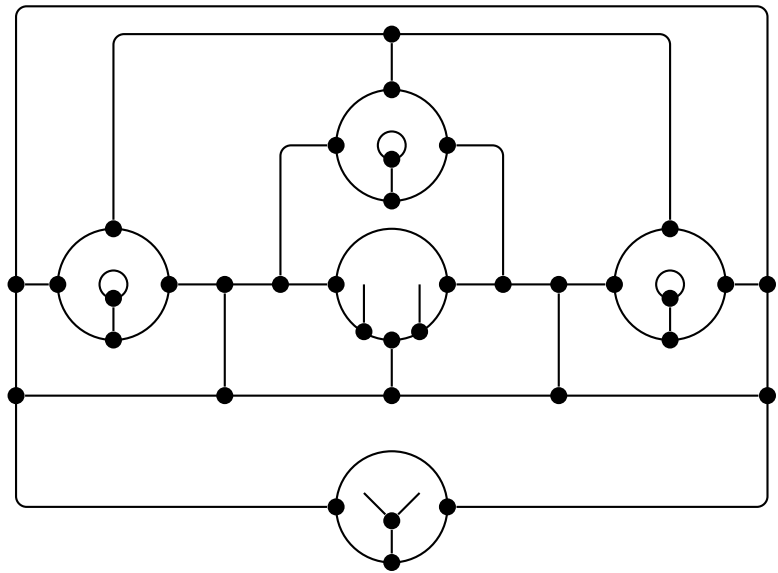


\end{document}